\documentclass[11pt, a4paper]{amsart}
\usepackage{etex}
\usepackage{amscd,amsmath,amsthm,amssymb}
\usepackage{mathtools}
\usepackage{mathrsfs}
\usepackage{comment}
\usepackage[all,cmtip]{xy}
\usepackage{xcolor}
\usepackage{enumitem} 
\usepackage{array}
\usepackage{caption}
\usepackage{mathabx}
\setlist[enumerate]{label=\textnormal{(\arabic*)}}

\usepackage{graphicx}

\usepackage{aurical}
\usepackage{amsbsy}
\usepackage{bm}

\usepackage{circuitikz}
\ctikzset{bipoles/resistor/height=0.15}
\ctikzset{bipoles/resistor/width=0.4}

\definecolor{cadmiumgreen}{rgb}{0.0, 0.42, 0.24}
\usepackage[
colorlinks, citecolor=cadmiumgreen,
pdfauthor={Omid Amini, Noema Nicolussi}, 
pdftitle={Higher rank Voronoi tilings},
pdfstartview ={FitV},
]{hyperref}

\definecolor{lblue}{rgb}{0.0, 0.55, 1.0}
\definecolor{dred}{rgb}{1.0, 0.0, 0.3}
\definecolor{byzant}{rgb}{0.74, 0.2, 0.64}
\definecolor{pomme}{rgb}{0.16, 0.67, 0.53}
\definecolor{aqua}{rgb}{0.0, 1.0, 1.0}
\definecolor{azure}{rgb}{0.0, 0.5, 1.0}
\definecolor{carnelian}{rgb}{0.7, 0.11, 0.11}
\definecolor{ceruleanblue}{rgb}{0.16, 0.32, 0.75}
\definecolor{pgreen}{rgb}{0.0, 0.65, 0.58}
 \definecolor{pblue}{rgb}{0.11, 0.22, 0.73}
\definecolor{pear}{rgb}{0.82, 0.89, 0.19}

\definecolor{pyellow}{rgb}{1.0, 0.97, 0.12}
\definecolor{lgreen}{rgb}{0.75, 1.0, 0.0}
\definecolor{brown}{rgb}{0.6, 0.4, 0.08}

\definecolor{ggray}{rgb}{0.27, 0.35, 0.27}

\usepackage{manfnt}

\usepackage[
alphabetic,
msc-links,
nobysame,
lite,
]{amsrefs} 

\usepackage{tikz, float} 
\usetikzlibrary {positioning}

\usepackage{pgfplots}
\usetikzlibrary{patterns}
\pgfplotsset{compat=1.10}
\usepgfplotslibrary{fillbetween}

\usetikzlibrary{matrix, arrows}
\usetikzlibrary{patterns}

\usetikzlibrary{calc,decorations.markings}
\usetikzlibrary{shapes,snakes}

\usepackage{a4wide}
\usepackage[T1]{fontenc}
\usepackage[OT1]{fontenc}
\usepackage[utf8]{inputenc}
\usepackage[english]{babel}
\usepackage{epsfig}
\usepackage[leqno]{amsmath}
\usepackage{tabularx}
\usepackage{subfigure}
\usepackage{url}
\usepackage{hyperref}
\usepackage{tikz-cd}
\usepackage{csquotes}

\usepackage{mathtools}
\usepackage{manfnt}
\usepackage{mathrsfs}
\usepackage{scalerel}

\DeclareRobustCommand{\SkipTocEntry}[5]{}

\usepackage{xcolor}
\definecolor{darkgreen}{HTML}{00AA00}

\usepackage{todonotes}

\newtheorem{thm}{Theorem}[section]

\newtheorem{lem}[thm]{Lemma}

\newtheorem{prop}[thm]{Proposition}

\newtheorem{claim}[thm]{Claim}
\newtheorem{cor}[thm]{Corollary}

\newtheorem{question}[thm]{Question}
\theoremstyle{definition}
\numberwithin{equation}{section}

\theoremstyle{definition}
\newtheorem{defii}[thm]{Definition}
\newenvironment{defi}
  {\pushQED{\qed}\defii}
  {\popQED\enddefii}
\newtheorem{remm}[thm]{Remark}
\newenvironment{remark}
  {\pushQED{\qed}\remm}
  {\popQED\endremm}
\newtheorem{exx}[thm]{Example}
\newenvironment{example}
  {\pushQED{\qed}\exx}
  {\popQED\endexx}
\newtheorem{convv}[thm]{Convention}
\newenvironment{convention}
  {\pushQED{\qed}\convv}
  {\popQED\endconvv}

\numberwithin{equation}{section}

\renewcommand{\~}{\widetilde}

\newcommand{\Q}{\mathbb{Q}}
\newcommand{\R}{\mathbb{R}}
\newcommand{\Z}{{\mathbb Z}}
\newcommand{\N}{{\mathbb N}}

\newcommand{\abs}[1]{\lvert #1\rvert}
\DeclareMathOperator{\rk}{rank} 
\newcommand{\rquot}[2]{#1/\,#2}       
\newcommand{\st}{\bigm|} 

\newcommand{\filter}{\mathrm{F}} 

\newcommand{\innone}[2]{{\langle#2\rangle_{_{\hspace{-.02cm}#1}}}}

\newcommand{\highinnoneorth}[2]{{\langleup{\vspace{.4cm}}#2\rangleup^{\perp}_{_{\hspace{-.02cm}#1}}}}
\newcommand{\innoneorth}[2]{{\langle#2\rangle^{\perp}_{_{\hspace{-.02cm}#1}}}}

\newcommand{\symbupl}{\tikz[scale=.1, baseline=-1.7]{\draw(.98,-.75)to[out=90,in=90](.98,.86) ;}} 
\newcommand{\langleup}{\langle \tikz[scale=.28, baseline=-0.7]{\symbupl}}
\newcommand{\symbupr}{\tikz[scale=.1, baseline=-1.7]{\draw(.91,-.77)to[out=90,in=90](.91,.86) ;}}
\newcommand{\rangleup}{\rangle\tikz[scale=.28, baseline=-0.7]{\symbupr}}
\newcommand{\highinn}[2]{{\langleup{\vspace{.4cm}}#2\rangleup_{_{\hspace{-.02cm}#1}}}}
\newcommand{\highinnprime}[2]{{\langleup{\vspace{.4cm}}#2\rangleup'_{_{\hspace{-.02cm}#1}}}}

\makeatletter
\let\@oldinfty\infty
\newcommand{\@sminfty}{{\hspace{-.02cm}\scaleto{\@oldinfty}{2.8pt}}} 
\newcommand{\@smsminfty}{_{\hspace{-.04cm}\scaleto{\@oldinfty}{2.8pt}}} 
\renewcommand{\infty}{{\mathchoice%
  {\displaystyle{\@oldinfty}}%
  {\textstyle{\@oldinfty}}%
  {\scriptstyle{\@sminfty}}%
  {\scriptscriptstyle{\@smsminfty}}}
}
\makeatother

\makeatletter
\let\@oldpi\pi
\newcommand{\@smsmpi}{_{\hspace{-.04cm}\scaleto{\@oldpi}{3.2pt}}} 
\renewcommand{\pi}{{\mathchoice%
  {\displaystyle{\@oldpi}}
  {\textstyle{\@oldpi}}%
  {\scriptstyle{\@oldpi}}%
  {\scriptscriptstyle{\@smsmpi}}}
}
\makeatother

\newcommand{\fin}{{\scaleto{\mathfrak f}{8pt}}}

\makeatletter
\let\@oldfin\fin
\newcommand{\@smfin}{{\hspace{-.01cm}\scaleto{\@oldfin}{5.4pt}}} 
\newcommand{\@smsmfin}{{\hspace{-.01cm}\scaleto{\@oldfin}{5pt}}} 
\renewcommand{\fin}{{\mathchoice%
  {\displaystyle{\@oldfin}}
  {\textstyle{\@oldfin}}%
  {\scriptstyle{\@smfin}}%
  {\scriptscriptstyle{\@smsmfin}}}
}
\makeatother

\newcommand{\sed}{{\scaleto{\mathrm{sed}}{7pt}}}

\makeatletter
\let\@oldsed\sed
\newcommand{\@smsed}{{\hspace{-.01cm}\scaleto{\@oldsed}{4.4pt}}} 
\newcommand{\@smsmsed}{{\hspace{-.01cm}\scaleto{\@oldsed}{4pt}}} 
\renewcommand{\sed}{{\mathchoice%
  {\displaystyle{\@oldsed}}
  {\textstyle{\@oldsed}}%
  {\scriptstyle{\@smsed}}%
  {\scriptscriptstyle{\@smsmsed}}}
}
\makeatother

\newcommand{\qf}{\mathfrak{q}}

\newcommand{\norm}[1]{\|#1\|}
\newcommand{\normsq}[1]{\|#1\|^{^{\scaleto{2}{3.5pt}}}}
\newcommand{\refnorm}[1]{|#1|}
\newcommand{\refnormsq}[1]{|#1|^{^{\scaleto{2}{3.5pt}}}}
\newcommand{\normorth}[1]{\|#1\|^{\perp}}

\newcommand{\dtor}{\zeta}

\newcommand{\Vor}{\mathrm{Vor}}
\newcommand{\interior}[1]{\mathrm{int}(#1)} 
\renewcommand{\L}{\mathbb L} 
\renewcommand{\S}{\mathbb S}
 
\newcommand{\vol}{\mathrm{vol}}

\newcommand{\T}{{\mathbb T}} 

\newcommand{\grind}[1]{{\scaleto{#1}{4pt}}}

\DeclareMathOperator{\Hom}{Hom}

\DeclareMathOperator{\id}{Id}

\renewcommand{\Re}{\operatorname{Re}}

\newcommand{\C}{{\mathbb C}}

\newcommand{\E}{{\mathbb E}}

\newcommand{\mg}{\mathscr M} 
\newcommand{\mgbar}{\comp{\mathscr M}} 
\newcommand{\mgg}[1]{\mg_{{\hspace{-.04cm}#1}}^{ }} 
\newcommand{\mgbarg}[1]{\mgbar_{{\hspace{-.09cm}#1}}^{ }} 

\makeatletter
\newsavebox\myboxA
\newsavebox\myboxB
\newlength\mylenA

\newcommand*\overbar[2][0.75]{%
    \sbox{\myboxA}{$\m@th#2$}%
    \setbox\myboxB\null
    \ht\myboxB=\ht\myboxA%
    \dp\myboxB=\dp\myboxA%
    \wd\myboxB=#1\wd\myboxA
    \sbox\myboxB{$\m@th\overline{\copy\myboxB}$}
    \setlength\mylenA{\the\wd\myboxA}
    \addtolength\mylenA{-\the\wd\myboxB}%
    \ifdim\wd\myboxB<\wd\myboxA%
       \rlap{\hskip 1\mylenA\usebox\myboxB}{\usebox\myboxA}%
    \else
        \hskip -0.5\mylenA\rlap{\usebox\myboxA}{\hskip 0.5\mylenA\usebox\myboxB}%
    \fi}

\newcommand{\comp}[1]{\overbar[.5]{#1}} 

\newcommand{\rsf}{\mathcal S} 
 
\newcommand{\graphgenus}{{h}} 
\newcommand{\mgr}{\mathcal G}  
\newcommand{\grm}[2]{\mathrm{gr}_{_{\hspace{-.08cm}#1}}^{#2}} 

\newcommand{\proj}{{\mathfrak p}}  
\newcommand{\projeuc}{\mathrm{p}}

\renewcommand{\curve}{\mathcal  C}

\renewcommand{\~}{\widetilde}

\newcommand{\Ical}{\mathcal I}
\newcommand{\dist}{\mathrm{d}}
\newcommand{\rest}[1]{\raisebox{-1pt}{$\vert$}_{#1}}

\newcommand{\suppp}[1]{{\mathrm{supp}(#1)^+}}
\newcommand{\suppm}[1]{{\mathrm{supp}(#1)^-}}
\newcommand{\supppm}[1]{{\mathrm{supp}(#1)^\pm}}

\newcommand{\smallcc}{{{\scaleto{\mathbb C}{4pt}}}}

\tikzstyle{Cwhite}=[scale = .8,circle, fill = white, minimum size=3mm] 
\tikzstyle{Cgray}=[scale = .4,circle, fill = gray, minimum size=3mm] 
\tikzstyle{Cblack2}=[scale = .4,circle, fill = black, minimum size=5mm] 
\tikzstyle{Cblack}=[scale = .7,circle, fill = black, minimum size=3mm]
\tikzstyle{C0}=[scale = .9,circle, fill = black!0, inner sep = 0pt, minimum size=3mm]
\tikzstyle{C1}=[scale = .7,circle, fill = black!0, inner sep = 0pt, minimum size=3mm]
\tikzstyle{Cred}=[scale = .4,circle, fill = red, minimum size=3mm]

\usepackage{stmaryrd}
\usepackage{trimclip}

\makeatletter
\DeclareRobustCommand{\shortto}{%
  \mathrel{\mathpalette\short@to\relax}%
}

\newcommand{\short@to}[2]{%
  \mkern2mu
  \clipbox{{.5\width} 0 0 0}{$\m@th#1\vphantom{+}{\shortrightarrow}$}%
  }
\makeatother

\newcommand{\transpose}{{\scaleto{\mathrm{T}}{4.5pt}}}

\newcommand{\lexeq}{\preceq_{{\scaleto{\mathrm{lex}}{3.5pt}}}}
\newcommand{\lexst}{\prec_{{\scaleto{\mathrm{lex}}{3.5pt}}}}
\newcommand{\gexeq}{\succeq_{{\scaleto{\mathrm{lex}}{3.5pt}}}}
\newcommand{\gex}{\succ_{{\scaleto{\mathrm{lex}}{3.5pt}}}}

\newcommand{\Jac}{\mathrm{Jac}}

\makeatletter
\newcommand{\orient}{{\mathchoice%
  {\displaystyle{\scaleto{\mathcal{O}}{5pt}}}%
  {\textstyle{\scaleto{\mathcal{O}}{5pt}}}%
  {\scriptstyle{{\scaleto{\mathcal{O}}{3.5pt}}}}%
  {\scriptscriptstyle{{\scaleto{\mathcal{O}}{3.5pt}}}}}
}
\makeatother

\renewcommand{\setminus}{\smallsetminus}

\newcommand{\symbuptm}{\tikz[scale=.28, baseline=-0.7]{\draw(.1,-.05)to[out=50,in=-130](.9,.05) (.5,0)--++(0,.5);}}
\newcommand{\almorth}{\tikz[scale=.28, baseline=-0.7]{\symbuptm \draw (-.44,.22) circle (14pt);}}
\newcommand{\aplus}{\hspace{-.2cm}\almorth\,}
\newcommand{\aperp}{\,\symbuptm\,}

\newcommand{\hdist}{d_{\mathrm{H}}}
\newcommand{\GH}{\scaleto{\mathrm{GH}}{4pt}}
\newcommand{\HA}{\scaleto{\mathrm{H}}{4pt}}
\renewcommand{\emptyset}{\varnothing}

\makeatletter
\@namedef{subjclassname@2020}{\textup{2020} Mathematics Subject Classification}
\makeatother

\newcommand{\msc}[1]{\href{http://www.ams.org/msc/msc2020.html?t=&s=#1}{#1}}

\begin{document}
\title[Higher rank inner products, Voronoi tilings and metric degenerations of tori]{Higher rank inner products, Voronoi tilings \\ and \\ metric degenerations of tori}

\author{Omid Amini}
\address{CNRS - Centre de math\'ematiques Laurent Schwartz, \'Ecole Polytechnique, Institut Polytechnique de Paris}
\email{\href{omid.amini@polytechnique.edu}{omid.amini@polytechnique.edu}}

\author{Noema Nicolussi}
\address{Faculty of Mathematics, Physics and Geodesy, Institute of Analysis and Number Theory, Graz University of Technology}
\email{\href{noema.nicolussi@univie.ac.at}{noema.nicolussi@univie.ac.at}}

\date{\today} 

\keywords{Higher rank inner products, multi-scale geometry, metric collapse, Voronoi tilings, asymptotic polyhedral geometry,  higher rank tropical geometry}
\subjclass[2020]{Primary  \msc{46C99}; \msc{52A23}; \msc{52C22}; Secondary \msc{05B45}; \msc{14T10}; \msc{14H15}; \msc{51F99}; \msc{51M20}; \msc{52B99}; \msc{58K99}; \msc{65D18}}

\begin{abstract}  
We introduce higher rank inner products on real and complex vector spaces and study their corresponding Voronoi tilings. We use the framework to describe metric degenerations of polarized tori and Hausdorff limits of Voronoi tilings of discrete sets. 

\end{abstract}

\maketitle

\setcounter{tocdepth}{1}

\tableofcontents

\section{Introduction}

The aim of this paper is to study degenerations of scalar products and their associated flat tori and Voronoi decompositions.

 Let $H$ be a finite dimensional real vector space. Consider a family $\innone{t}{\cdot\,, \cdot}$ of scalar products on $H$ parametrized by positive real numbers $t\in \R_+$. The question we address in this paper can be informally stated as follows.

 \begin{question} Is there a notion of ``\,limit" for the scalar products $\innone{t}{\cdot\,, \cdot}$ as $t$ tends to infinity?
\end{question}

This question is interesting in the case where the family is degenerating, which means that for some elements $x \in H$, the pairing $\innone{t}{x, x}$ tends either to $0$ or to $\infty$ as $t \to\infty$.

 Let $\T =\rquot{H}{\L}$ be the real torus associated to a full rank lattice $\L \subseteq H$.  A scalar product on $H$ induces a flat metric on $\T$ with associated distance function $\dist \colon \T\times \T \to [0, +\infty)$. We are particularly interested in the following related problem.

\begin{question} \label{Q:TorusIntro} Given a degenerating family of scalar products $\innone{t}{\cdot\,, \cdot}$ on $H$, describe the asymptotic behavior of the distance functions $\dist_t \colon \T \times \T \to [0, \infty)$ on the torus $\T$ for $t \to \infty$.
\end{question}

The notion of Voronoi decomposition allows to connect Question~\ref{Q:TorusIntro} to polyhedral geometry.
Recall that, given a discrete set $S \subseteq H$ and a scalar product $\innone{}{\cdot\,, \cdot}$ on $H$, we obtain a decomposition of $H$ into the so-called Voronoi cells $\Vor(\gamma)$, $\gamma\in S$, defined as the set of points in $H$ which have $\gamma$ as their nearest point in $S$.  For a full rank lattice $S = \L$, the distance function $\dist$ can be expressed in terms of the Voronoi decomposition.

We are thus led to the following question.

\begin{question} \label{Q:VoronoiIntro} 
Given a discrete subset $S\subseteq H$ and a degenerating family of scalar products $\innone{t}{\cdot\,, \cdot}$, $t \in \R_+$, describe the behavior of the induced Voronoi decompositions for $t \to \infty$.
\end{question}
 
 The above problems appear naturally in complex geometry and the work presented here has its origin in our series of papers~\cites{AN, AN2}. In these works, we develop a higher rank, multi-scale geometry in the setting of hybrid Riemann surfaces (which mix metric graphs and Riemann surfaces) and apply the framework to address questions on the asymptotic geometry of degenerating Riemann surfaces. In line with this approach, we can ask the following question. To any Riemann surface $S$ is associated its Jacobian torus $\T= \rquot{H_1(S, \R)}{H_1(S, \Z)}$. The complex structure of $S$ induces a scalar product on $H_1(S, \R)$ and a metric on $\T$, which is called polarization.  Let $\mgg{g}$ be the moduli space of Riemann surfaces of genus $g$ and let $\mgbarg{g}$ be its Deligne-Mumford compactification. 
  
\begin{question}\label{Q:Jacobians} Consider a family of smooth compact Riemann surfaces $S_t$ of genus $g$, $t \in \R_+$,  whose associated points $s_t$ in  $\mgg{g}$ converge to a point $s_\infty$ in $\mgbarg{g}$. For each $t$, let $\T_t$ be the polarized Jacobian of $S_t$.
 
Describe the metric behavior of the flat tori $\T_t$ as $t \to\infty$.
\end{question}

The main difficulties in these questions stem from multi-scale behavior appearing for degenerating scalar products. The pairing $\innone{t}{x,x}$ might tend to $+ \infty$ for some elements $x\in H$, while remaining bounded or going to zero for others. Questions~\ref{Q:TorusIntro} and~\ref{Q:VoronoiIntro} require to understand the effect of this multi-scale behavior on Voronoi decompositions and flat tori.

 \smallskip

 The aim of this paper is to answer the above questions by developing a limit theory of scalar products. This will be based on a generalization of the concept itself. More specifically...

  \smallskip
  
$\bullet$ We introduce \emph{inner products of higher rank}. These are symmetric bilinear forms with values in the vector space $\R^r$ which enjoy a suitable definiteness property with respect to the lexicographic order on $\R^r$. We interpret them as multi-scale limits of scalar products. The integer $r$ corresponds to the number of scales in a degenerating family of scalar products. 

  \smallskip
$\bullet$ We show that a higher rank inner product induces an \emph{almost orthogonal decomposition} of $H$ into subspaces $H_1, \dots, H_r$, on each of which, the behavior of a degenerating family of scalar products close to the limit is almost homogeneous.   
	
  \smallskip
$\bullet$ We introduce Voronoi decompositions with respect to inner products of higher rank. We prove that the cells in the Voronoi decomposition are the Hausdorff limits of Voronoi cells associated to degenerating scalar products. Voronoi cells of higher rank are not necessarily closed.  We show that the closure of each Voronoi cell decomposes into a sum of Voronoi cells of scalar products, each living in a homogeneous part of the almost orthogonal decomposition. 

 \smallskip
$\bullet$ We analyze multi-scale phenomena for the distance function on tori associated to scalar products degenerating to an inner product of higher rank. We single out subtori which exhibit different scales, and show that, upon proper renormalization, they converge to tori of possibly lower dimension in the Gromov--Hausdorff sense. We also provide an asymptotic expansion for the distance function in terms of the limiting inner product of higher rank.

 \smallskip
$\bullet$ As a direct application, we describe the asymptotic behavior of Jacobian tori and Voronoi decompositions associated to degenerating metric graphs. We also provide an informal overview of the application given in our forthcoming work~\cite{AN-AG-hybrid} to the study of the Jacobians of degenerating Riemann surfaces, which answers Question~\ref{Q:Jacobians}.

\smallskip

 In Section~\ref{ss:HermitianInnerProducts}, we explain how the set-up can be adapted to treat the same type of questions for Hermitian inner products on complex vector spaces, leading to a theory of higher rank Hermitian inner products. All the results discussed in this introduction and in the paper can be extended to complex vector spaces.
 
\subsection{Overview of results}

In the following, we state our main results.

\subsubsection{Definition of higher rank inner products}

 Let $r$ be a positive integer. Consider the real vector space $\R^r$ endowed with the lexicographic order $\lexeq$. Recall that for two points $a =(a_1, \dots, a_r)$ and $b=(b_1, \dots, b_r)$ in $\R^r$, we have $a \lexeq b$ if either $a=b$, or there exists $j\in [r]$ such that $a_i=b_i$ for all $i < j$ and $a_j<b_j$. Here and below, we let $[r] \coloneqq \{1, \dots, r\}$.

\smallskip

Consider a symmetric bilinear form $\highinn{}{\cdot \,, \cdot} \colon H \times H \to \R^r$ with coordinates $\highinn{j}{\cdot\,,\cdot}$, $j \in [r]$. Then $\highinn{}{\cdot \,, \cdot}$ defines a non-increasing filtration
\[
\filter^\bullet\colon \qquad \filter^1\coloneqq H \supseteq \filter^2 \supseteq \filter^3 \supseteq \dots \supseteq \filter^r \supseteq \filter^{r+1} \coloneqq (0)
\]
on $H$ by setting
\begin{align*}
\filter^j  \coloneqq \left\{ x \in H \st \,\highinn{}{x \,, y}_i = 0 \text{ for all $i <j$ and all $y \in H$}\right\}, \qquad j \in [r].
\end{align*}

We say that $\highinn{}{\cdot\,, \cdot}$ is an \emph{inner product},
if for each $j \in [r]$, the induced form $\highinn{j}{\cdot\,, \cdot}$ on $\filter^j / \filter^{j+1}$ is positive definite. Equivalently, if $\highinn{j}{\gamma, \gamma} >0$ for all $\gamma \in \filter^j \setminus \filter^{j+1}$. 

\smallskip

In case that $r=1$, we recover the classical definition. In this paper, we use the terminology  \emph{scalar product} for inner products with values in $\R$. If the image of an inner product $\highinn{}{\cdot\,, \cdot}$ is a one-dimensional subspace of $\R^r$, then  by fixing an isomorphism with $\R$, we obtain a scalar product. The above definition thus generalizes the notion of inner products to higher rank, in the sense that their image can have dimension higher than one.  By a slight abuse of notation, we refer to inner products with values in $\R^r$, for arbitrary positive integer $r$, as \emph{inner products of higher rank} or simply \emph{inner products} (although, strictly speaking, $r=1$ is allowed and, for larger $r$, their image could be a one-dimensional subspace of $\R^r$). 

\subsubsection{Higher rank inner products as limits of scalar products}
It turns out that higher rank inner products are closely related to multi-scale degenerations of scalar products.

\smallskip

Consider a symmetric bilinear form $\highinn{}{\cdot \,, \cdot} \colon H \times H \to \R^r$. Let $\R_+ = (0, +\infty)$. Given a vector $\underline L = (L_1, \dots, L_r)$ of reals $L_j \in \R_+$, we define the scalar bilinear form $\innone{\underline L}{\cdot\,, \cdot}\colon H\times H \to \R$ as
\[
\innone{\underline L}{x, y} \coloneqq \sum_{j=1}^r L_j \highinn{j}{x, y}, \qquad x,y \in H.
\]
We call $\innone{\underline L}{\cdot\,, \cdot}\colon H\times H \to \R$ the \emph{pullback} of $\highinn{}{\cdot \,, \cdot}$ by the vector $\underline L$. 

\smallskip

The following result provides an alternative, more analytic definition of inner products in terms of pullbacks (see Section~\ref{ss:PullbackFamilies}).
\begin{thm}\label{thm:pullback-intro} Let $\highinn{}{\cdot\,, \cdot} \colon H \times H \to \R^r$ be a symmetric bilinear form. The following statements are equivalent...

\smallskip

\noindent $(i)$ $\highinn{}{\cdot\,, \cdot}$ is an inner product.

\smallskip

\noindent $(ii)$ for any family $\underline L_t$,  $t \in \R_+$, of vectors $\underline L_t = (L_{t,1}, \dots, L_{t,r}) \in \R_+^r$,  with
\begin{equation} \label{eq:MultiScaleLimit-intro}
\lim_{t \to \infty} \frac{L_{t,j}}{L_{t, j+1}} =   \infty, \qquad j=1, \dots, r-1,
\end{equation}
the pullbacks $\innone{\underline L_t}{\cdot\,, \cdot}\colon H\times H \to \R$ are scalar products for all large enough $t \in \R_+$. 
\end{thm}
Given an inner product $\highinn{}{\cdot\,, \cdot} \colon H \times H \to \R^r$ and a family of vectors $\underline L_t$, $t \in \R_+$, satisfying \eqref{eq:MultiScaleLimit-intro}, we call the resulting family of scalar products $ \innone{t}{\cdot\,, \cdot} \coloneqq \innone{\underline L_t}{\cdot\,, \cdot}$, $t \in \R_+$, a \emph{pullback family} for $\highinn{}{\cdot\,, \cdot}$ with \emph{parameters} $\underline L_t = (L_{t,1}, \dots, L_{t,r})$, $t \in \R_+$.

\smallskip

By property \eqref{eq:MultiScaleLimit-intro}, the parameters $L_{t,j}$ exhibit different growth rates for $t \to \infty$. For instance, for the choice $L_{t,j} = t^{r-j}$, the first $(r-1)$ parameters tend to $+\infty$ with different speeds, while the last one remains constant. Note that for $j \in [r]$ and $x \in \filter^j\setminus \filter^{j+1}$, the scalar product behaves like
\[
\innone{t}{x,x} \sim L_{t,j} \, \highinn{j}{x,x}, \qquad t \to \infty.
\]
The scalar products $\innone{t}{\cdot \,, \cdot}$ thus show a multi-scale behavior with $r$ different scales in the limit $t \to \infty$. Since $L_{t,1} \gg L_{2,t} \gg \dots \gg L_{t,r}$ for $t \to \infty$ by \eqref{eq:MultiScaleLimit-intro}, the first parameter $L_{t,1}$ corresponds to the ``dominant scale", $L_{t,2}$ the ``second dominant scale", etc.

\smallskip

Thus, Theorem~\ref{thm:pullback-intro} relates higher rank inner products $\highinn{}{\cdot\,, \cdot}$ to families of scalar products with multi-scale behavior. On the other hand, the structure of the scalar products $\innone{t}{\cdot \,, \cdot} = \innone{\underline L_t}{\cdot \,, \cdot}$ in pullback families is rather specific. It is clearly desirable to include also other families of scalar products. For instance, in geometric applications to degenerations of metric graphs and Riemann surfaces (see Section~\ref{sec:tropical_curves}), the scalar products that appear are not necessarily pullbacks.

In order to remedy this, we introduce the notions of \emph{tamely} and $\omega$-\emph{tamely} (short for weak-tamely) degenerating families of scalar products. These are families of scalar products $\innone{t}{\cdot \,, \cdot}$, $t \in \R_+$,
that generalize pullback families, and which naturally appear when studying asymptotic geometry of complex algebraic varieties. Their properties allow to carry through similar arguments used in the proofs for pullback families. We refer to Section~\ref{ss:tameness} for the precise definition.

A natural idea is then to view  $\highinn{}{\cdot\,, \cdot}$ as a "limit" of its $(\omega$-)tamely degenerating families. 
As we outline briefly in Sections~\ref{ss:TameEquivalence} and~\ref{sec:tame-compactifications}, this idea leads to Hausdorff partial compactifications  of the cone of scalar products on $H$.

\subsubsection{Degeneration problem for tori and Voronoi decompositions}
Let $\innone{t}{\cdot \,, \cdot}$, $t \in \R_+$, be scalar products on $H$. Consider the real torus $\T= \rquot{H}{\L}$ associated to a full-dimensional lattice $\L\subseteq H$. Each scalar product $\innone{t}{\cdot \,, \cdot}$ defines a Riemannian metric on $\T$ and its universal cover $H$. Let $\dist_t \colon \T \times \T \to [0,+\infty)$ be the associated distance function. We are interested in the asymptotic behavior of the metric spaces $(\T, \dist_t)$ as $t \to \infty$. This question is closely related to the asymptotic behavior of the associated Voronoi decomposition. 

Recall that, for a discrete subset $S \subset H$ and a scalar product $\innone{}{\cdot \,, \cdot} \colon H \times H \to \R$, the space $H$ can be decomposed into the \emph{Voronoi cells} $\Vor(\gamma)$, $\gamma \in S$, given by
\begin{equation} \label{eq:ScalarVoronoiIntro}
\Vor(\gamma) = \big \{x \in H \st \innone{}{x-\gamma \,, x-\gamma} \leq \innone{}{x-\eta \,, x-\eta} \, \forall \eta \in S \big \}.
\end{equation}
The Voronoi cells $\Vor(\gamma)$ are convex (generalized) polyhedra with mutually disjoint interiors (see Section~\ref{ss:BasicNotations} for the notion of generalized polyhedra). Together, they form the \emph{Voronoi decomposition} $H = \bigcup_{\gamma \in S} \Vor(\gamma)$ of the space.
For a full rank lattice $S = \L$,
the Voronoi cell $\Vor(0)$ is a polyhedron and defines a fundamental domain in the universal cover $H \to \T$. 
The distance function $\dist \colon H \times H \to [0, \infty)$ can be written as 
\[
\dist(p,q) = \sqrt{\innone{}{x,x}}, \qquad p,q \in \T,
\]
where $x \in H$ is any representative of $p-q$ with $x \in \Vor(0)$. Given a family $\innone{t}{\cdot \,, \cdot}$, $t \in \R_+$, of scalar products on $H$ and a fixed discrete subset $S \subset H$, we are interested in the behavior of the associated Voronoi cells $\Vor_t(\gamma)$ as $t \to \infty$.

\subsubsection{Higher rank Voronoi decompositions}
In order to answer the above question, we introduce Voronoi decompositions for inner products of higher rank.

\smallskip

Let $\highinn{}{\cdot\,, \cdot} \colon H \times H \to \R^r$ be an inner product and $S \subset H$ discrete. Analogous to \eqref{eq:ScalarVoronoiIntro}, we define the (\emph{higher rank}) \emph{Voronoi cells} as
\begin{equation} \label{eq:HIgherRankVoronoiIntro}
	\Vor(\gamma) =  \big \{x \in H \st \highinn{}{x-\gamma \,, x-\gamma} \leq \highinn{}{x-\eta \,, x-\eta} \, \,\, \forall \eta \in S  \big \}, \qquad \gamma \in S.
\end{equation}

\smallskip

In general, higher rank Voronoi cells behave very differently from their classical counter-parts and certain pathologies may appear (see Section~\ref{ss:Pathologies}). In order to exclude this, we introduce the class of \emph{admissible discrete subsets} for the higher rank inner product $\highinn{}{\cdot\,, \cdot}$. 
A full rank lattice $\L \subset H$ is admissible exactly when the intersections $\L \cap \filter^j$ are of full rank for all $j \in [r]$, for the filtration $\filter^\bullet$ induced by $\highinn{}{\cdot\,, \cdot}$. For the general definition, we refer to Section~\ref{ss:AdmissibleDiscreteSetsAndVoronoiDecomposition}.
In Section~\ref{ss:HybridTopology}, we prove that admissible discrete sets coincide with the discrete sets for a certain topology on $H$ that combines features from the order topology on $(\R^r, \lexeq)$ and the Euclidean topology on $H$.

\smallskip

For an admissible discrete subset $S \subset H$, the Voronoi cells behave well (see Section~\ref{sec:HigherRankVoronoiDecompositions}). First of all, we show that they cover $H$, that is, we obtain a (\emph{higher rank}) \emph{Voronoi decomposition} $H = \bigcup_{\gamma \in S} \Vor(\gamma)$ (see Theorem~\ref{thm:voronoi-discrete}). Moreover,  we show that the closures of Voronoi cells are full dimensional convex (generalized) polyhedra with mutually disjoint interiors (see Corollary~\ref{cor:covering-discrete}). In fact, these cells enjoy a special structure and can be decomposed into a sum of Voronoi cells associated to scalar products. More precisely, we show that the space $H$ admits a canonical \emph{almost orthogonal decomposition}
\[
H = H_1 \aplus H_2 \aplus \dots \aplus H_r
\]
into a direct sum of subspaces $H_1, \dots, H_r$ such that $H_j$, $j\in[r]$, is isomorphic to the $j$-th graded piece $\grm{}{j} H = \rquot{\filter^j}{\filter^{j+1}}$ of the filtration $\filter^\bullet$, and the subspaces $H_j$, $j \in [r]$, are pairwise \emph{almost orthogonal} with respect to $\highinn{}{\cdot \,, \cdot}$. We refer to Section~\ref{ss:AlmostOrthogonalDecomposition} for the precision definition and more details.  
 The notation $\aplus$ that is used here, instead of $\oplus$, refers to the almost orthogonal property.  
 
The $j$-th component $\highinn{j}{\cdot \,, \cdot}$ of $\highinn{}{\cdot \,, \cdot}$ defines a scalar product on $H_j$.   
Using this decomposition, we then show that the closure of a higher rank Voronoi cell is a sum
	\[
\overline \Vor(\gamma) \cong V_{\gamma,1} \aplus V_{\gamma,2} \aplus \dots \aplus V_{\gamma,r}
\]
of Voronoi cells $V_{\gamma,j} \subset H_j$ for the scalar products $\highinn{j}{\cdot \,, \cdot}$ on $H_j$ and a suitable projection $S_{j,\gamma} \subset H_j$ of $S$ to $H_j$ (see Theorem~\ref{thm:VoronoiCellDecomposition-discrete}).

\smallskip

 We prove that the higher rank Voronoi cells are the limits of Voronoi cells for degenerating scalar products (see Theorem~\ref{thm:Hausdorff_voronoi}).

\begin{thm}\label{thm:Hausdorff_voronoi_Intro} Let $S \subset H$ be an admissible discrete set and $\innone{t}{\cdot\,, \cdot}$, $t \in \R_+$, a  tamely degenerating  family of scalar products for $\highinn{}{\cdot\, , \cdot}$.
For $\gamma \in S$, let $\Vor(\gamma)$ and $\Vor_t(\gamma)$, $t \in \R_+$, be the Voronoi cells for $\highinn{}{\cdot\, , \cdot}$ and $\innone{t}{\cdot\,, \cdot}$, respectively. 
	
Then, $\Vor_t(\gamma)$ converges compactly in the Hausdorff sense to $\Vor(\gamma)$ as $t \to \infty$. In particular,  if $\Vor (\gamma) \subset H$ is relatively compact, then $\Vor_t(\gamma)$ converges to $\Vor(\gamma)$ in the Hausdorff sense as $t \to \infty$. 
\end{thm} 
The Hausdorff convergence here is with respect to some (equivalently any) reference norm $\refnorm{\cdot} \colon H  \to [0, \infty)$. Compact Hausdorff convergence is a suitable generalization in the context of non-compact sets. We refer to Section~\ref{sec:Hausdorff_convergence} for the definition and to Theorem~\ref{thm:CompactHausdorffConvergence} proved in Section~\ref{ss:CompactHausdorffConvergence} for an alternative characterization. Relative compactness of subsets $V\subset H$ refers to the Euclidean topology on $H$.

\smallskip

The proof of the above result is based on a subtle finiteness lemma established in Section~\ref{sec:FinitenessLemma}.

\subsubsection{Metric degenerations of tori}
In the following, we summarize our results on metric degenerations of tori.
\smallskip

Let $\highinn{}{\cdot \,, \cdot} \colon H \times H \to \R^r$ be an inner product with induced filtration $\filter^\bullet$. For an admissible fulll rank lattice $\L \subset H$, the lattices $\L^j =  \L \cap \filter^j$, $j \in [r]$, are of full rank in $\filter^j$. Setting $\T^j \coloneqq \rquot{\filter^j}{\L^j}$, $j\in [r]$, we obtain a non-increasing filtration of tori
\[\T= \T^1\supseteq \T^2 \supseteq \dots \supseteq \T^r\supseteq \T^{r+1}=(0).\]

For $j \in [r]$, the lattice $\grm{}{j}\L = \rquot{\L^j}{\L^{j+1}}$ is of full rank in $\grm{}{j} H =\rquot{ \filter^j }{\filter^{j+1}}$. We equip the torus $\Theta_j \coloneqq \rquot{\grm{}{j}H}{\grm{}{j}\L}$ with the distance function $\dist_j$ induced by the scalar product $\highinn{j}{\cdot \, ,\cdot}$ on $\grm{}{j} H$.

The following theorem demonstrates that different parts of the torus exhibit different scales and, under proper normalization,  collapse to the tori $\Theta_j  = \rquot{\grm{}{j}H}{\grm{}{j}\L}$ (see Theorem~\ref{thm:GHConvergenceAbstract}).

\begin{thm} \label{thm:GHConvergenceIntro}
Consider an inner product $\highinn{}{\cdot \,, \cdot} \colon H \times H \to \R^r$ and an admissible lattice $\L \subset H$ of full rank. Let $\innone{t}{\cdot\,, \cdot}$, $t \in \R_+$, be an  $\omega$-tamely degenerating family of scalar products with parameters $\underline L_t =(L_{t,1}, \dots, L_{t,r})$. Then, 
\[
\Big (\T, \frac{1}{\sqrt{L_{t,1}}} \, \dist_t \Big ) \xrightarrow{} \Big ( \Theta_1 , \, \dist_{1} \Big )
\]
in the Gromov--Hausdorff sense as $t \to \infty$. More generally, for all $j=1, \dots, r$,
\[
\Big ( \T^j,\frac{1}{\sqrt{L_{t,j}}} \,\dist_t \rest{\T^j\times \T^j} \Big ) \xrightarrow{} \Big (\Theta_j , \,  \dist_{j} \Big )
\]
in the Gromov--Hausdorff sense as $t \to \infty$.
\end{thm}

Moreover, in Theorem~\ref{thm:volume_degeneration}, we describe the asymptotics of the volume of $\T$ for $t \to \infty$ in terms of volumes of the tori $\Theta_j$ and the degeneration parameters $L_{t,j}$, $j \in [r]$.

\smallskip

We would like to emphasize here that the above results are surprisingly subtle. For instance, the simplest application concerns the torus $\T= \rquot{\R^2}{\Z^2}$ equipped with the distance function induced by the scalar product $\innone{t}{x,y} = x_1 y_1 + t^{-1} x_2 y_2$, $x,y \in \R^2$. In this special case, the above results are easy to prove and in particular, the torus converges to the one-dimensional circle $\rquot{\R}{\Z}$ for $t \to\infty$. On the other hand, in Section~\ref{ss:PathologyTorus}, we present a rather similar example which shows a drastically different behavior (the torus collapses to a point). From this perspective, higher rank inner products and tame degenerations seem to provide a framework in which metric degenerations happen in a rather controlled way.

\smallskip

For pullback families, we obtain the following refined asymptotics.
Define the function $\dtor = (\dtor_1, \dots, \dtor_r) \colon \T\times \T \to \R^r$ by
\begin{align*}
	\dtor(p,q) \coloneqq \min_{\substack{x \in H\\ x = p-q \textrm{ in } \T}} \highinn{}{x,x},
\end{align*}
where the minimum is taken in the lexicographic order $\lexeq$ on $\R^r$. In Theorem~\ref{thm:metric-asymptotics-pullback},  we prove that for every fixed pair $(p,q) \in \T \times \T$,
\[
\dist_t(p,q)  = \sqrt{L_{t,1} \dtor_1(p,q) + \dots + L_{t,r} \dtor_r(p,q)}
\]
for all large $t \in \R_+$. Using the higher rank Voronoi cell $\Vor(0)$, one can formulate a uniform approximation result as well (see Theorem~\ref{thm:metric-asymptotics-pullback}).

\subsubsection{Geometric applications}
In Section~\ref{ss:degeneration-metric-graphs}, we discuss applications of our results to degenerations of metric graphs and Riemann surfaces.

Recall that a compact metric graph $\mgr$ is obtained from a finite combinatorial graph $G=(V,E)$ by identifying its edges $e \in E$ with intervals of certain lengths $\ell(e)$, $e \in E$. The first homology $H = H_1(\mgr, \R)$ carries a natural scalar product, which is called the \emph{polarization} of $\mgr$ \cite{BLN97, KS00, BR07, MZ08, CV10}. The homology $\L= H_1(\mgr, \Z)$ is a full rank lattice in $H_1(\mgr, \R)$ and the torus $\T= \rquot{H}{\L}$ is called the Jacobian of $\mgr$.
In Section~\ref{ss:degeneration-metric-graphs}, we describe the asymptotics of the Jacobian and the associated Voronoi decomposition, when the metric graph $\mgr$ degenerates, meaning that some of the edges become infinitely long. Here we use a notion of tropical curves introduced in \cite{AN,AN2}. Briefly speaking, these are multi-scale versions of metric graphs which allow to handle multi-scale effects appearing for degenerating metric graphs.

In Section~\ref{ss:DegenerationsRiemannSurfaces}, we describe informally the application of our results to Jacobians of degenerating Riemann surfaces, which will appear in our forthcoming work \cite{AN-AG-hybrid}. Loosely speaking, in the simple case discussed here, the Jacobians of a degenerating family of Riemann surfaces, upon proper rescaling, converge to the Jacobian of a metric graph. Moreover, certain subtori of the Jacobians converge to Jacobians of Riemann surfaces of smaller genera. (No previous knowledge of complex geometry is needed in this paper.)

\subsection{Related work} Starting from the pioneering work by Gromov~\cite{Gro78}, Ruh~\cite{Ruh82}, Cheeger--Gromov~\cite{CG86, CG90}, Fukaya~\cite{Fuk87, Fuk90}, Cheeger--Fukaya--Gromov~\cite{CFG92}, there has been a large amount of literature on metric degenerations of Riemannian manifolds, captured in the fundamental concept of Gromov--Hausdorff convergence. We refer to the survey paper by Fukaya~\cite{Fuk90}, and the books \cite{Burago} and~\cite{G07}.  In this context, the flat tori $\T_\varepsilon = \rquot{\R^2}{\Z \oplus\varepsilon \Z}$, $\varepsilon >0$, provide a basic example of a family of Riemannian manifolds exhibiting  a dimension collapse in the limit. Moreover, tori and their generalizations appear as fibers in the collapse of almost flat manifolds by the Gromov--Ruh theorem~\cite{Gro78, Ruh82}.

More recent work in connection to complex geometry concerns Gromov-Hausdorff  convergence of Calabi-Yau manifolds endowed with their canonical metric in a maximally degenerate family, which is the subject of a conjecture by Kontsevich and Soibelman~\cite{KS}, in connection to the SYZ conjecture in mirror symmetry. In this regard,  the result we prove in~\cite{AN-AG-hybrid} can be viewed as a multi-scale variant of the Kontsevich--Soibelman conjecture in the specific case of (non-necessarily maximally degenerate) families of Abelian varieties.

That said, the type of questions investigated in this paper and the multi-scale geometric phenomena developed here appear to be new and not previously considered in the literature.

\smallskip

Although higher rank valuations do not appear in this paper, the framework developed here is also related conceptually to recent developments around higher rank non-archimedean geometry, see e.g. \cite{Aro10, AI22, Ban15, FR16, HIL20, Iri22, JS23}. Note that the space $(\R^r, \lexeq)$ and its topology plays a key role in \cite{HIL20, AI22,Iri22}.  Moreover, limiting behavior of certain convex bodies arising in combinatorics and algebraic geometry has been studied recently in \cite{ABGJ18, ABGJ22, AI22}.

\smallskip

Taking into account the results of this paper, a natural question is whether higher rank inner products and tame degenerations can be used to define a meaningful (at least partial) compactification of the space of scalar products on a fixed vector space $H$. We discuss this in Section~\ref{sec:tame-compactifications}.

\smallskip

Voronoi decompositions play a central role in applications of computational geometry to other fields. We refer to the survey article~\cite{Aur91} and the books~\cite{BY98, BSCO09} for the computational geometry literature and a sample of these applications.

\subsection{Structure of the paper}
Section~\ref{sec:HigherRankInnerProducts} contains the definition of inner products with values in $(\R^r, \lexeq)$, the relationship to pullbacks and basic algebraic constructions. In Section~\ref{sec:HigherRankVoronoiDecompositions}, we introduce higher rank Voronoi cells, admissible discrete subsets, and prove that this leads to a Voronoi decomposition. In Section~\ref{sec:StructureVoronoiCells}, we study the structure of the higher rank Voronoi cells. Section~\ref{sec:tameness} contains the definition of tamely degenerating families of scalar products and their basic properties. Section~\ref{sec:FinitenessLemma} provides an auxiliary lemma for subsequent considerations. In Section~\ref{sec:Hausdorff_convergence}, we prove the convergence of Voronoi cells to higher rank Voronoi cells. In Section~\ref{sec:AdmissibleLattices}, we specialize the preceding results to the case of lattices. Section~\ref{sec:metric_degeneration_general_tori} contains our results on degenerations of tori. Section~\ref{sec:tropical_curves} concerns applications to metric graphs and Riemann surfaces. Finally, Section~\ref{sec:Discussion} contains further results, discussions and questions.

\subsection{Basic notations} \label{ss:BasicNotations}  For a point $a$ in the vector space $\R^r$, we denote by $a_1, \dots, a_r$ the coordinates of $a$. We set $\R_+ \coloneqq (0, +\infty)$ and denote by $\R_+^r =(0,+\infty)^r$ the set of $r$-vectors with positive coordinates.

Let $H$ be a finite dimensional real vector space. In this paper, a polyhedron is a subset $P \subset H$ which can be written as an intersection of finitely many half-spaces in $H$. A polytope is a bounded polyhedron. A generalized polyhedron is a subset $P \subset H$ such that $P \cap Q$ is a polytope for every polytope $Q \subset H$. The terminology is borrowed from~\cite[Page 250]{Gruber}.

A (generalized) polyhedral tiling of a vector space $H$ in this paper means a (possibly infinite) collection $\Sigma$ of full dimensional (generalized) polyhedra $\sigma \subseteq H$ which verifies the following properties: 
\begin{itemize}
	\item The union of $\sigma \in \Sigma$ is $H$.
	\item The interiors of $\sigma \in \Sigma$ are disjoint open sets in $H$.
\end{itemize}

Note that this definition allows two top-dimensional tiles to intersect only on part of their proper faces, see Figure~\ref{fig:quasi-tiling} for an example.

We use the notation $\highinn{}{\cdot\,,\cdot}$ for inner products with values in $\R^r$, $r\in \N$. For scalar products, we usually use $\innone{}{\cdot\,,\cdot}$.

\subsection*{Acknowledgments} O.\thinspace A. acknowledges support from Math+, the Berlin Mathematics Research
Center, and thanks the hospitality of the mathematics institutes at TU and HU Berlin where part of this research was carried out.  N.\thinspace N. acknowledges financial support by the Austrian Science Fund (FWF) under Grant No. J 4497.


\section{Higher rank inner products} \label{sec:HigherRankInnerProducts}

\subsection{Totally ordered vector spaces}
Let $(\Lambda, \preceq)$ be a totally ordered real vector space of finite dimension. It is known that $\Lambda$ is isomorphic to $\R^r$ with its lexicographic order $\preceq=\lexeq$ where $r$ is the dimension of $\Lambda$~\cite{Bir-Lattice, HW52}. This means there is no loss of generality in assuming that $(\Lambda, \preceq) = (\R^r, \lexeq)$. For $a =(a_1, \dots, a_r)$ and $b=(b_1, \dots, b_r)$ in $\R^r$, we have $a \lexeq b$ if either $a=b$, or there exists $j\in [r]$ such that $a_1=b_1, \dots, a_{j-1}=b_{j-1}$ and $a_j<b_j$.

The vector space $\Lambda = \R^r$ with its lexicographic order $\lexeq$ has a natural decreasing filtration 
\begin{equation} \label{eq:FiltrationRr}
\Lambda^\bullet  \colon  \qquad \Lambda^1 \coloneq \Lambda\supset \Lambda^2 \supset \dots \supset \Lambda^{r} \supset \Lambda^{r+1} \coloneq (0)
\end{equation}
defined by
\[
\Lambda^j \coloneqq \left\{a = (a_1, \dots, a_r)  \in \R^r \st a_1 = \dots =a_{j-1}=0\right\}.
\]
Note that $\Lambda^j$ is naturally isomorphic to $\R^{r+1-j}$ endowed with its lexicographic order.

In this paper, we introduce and study inner products with values in $(\R^r, \lexeq)$. The notion can be generalized to handle inner products with values in general totally ordered vector spaces $(\Lambda, \preceq)$. However, upon fixing an isomorphism $(\Lambda, \preceq) \cong (\R^r, \lexeq)$, we can assume that $(\Lambda, \preceq) = (\R^r, \lexeq)$, which has the advantage of simplifying the notation. For this reason, we treat inner products with values in $ (\R^r,   \lexeq)$, and discuss the general case in Section~\ref{ss:InnerProductsGeneralValueGroup}.

\subsection{Inner product spaces of higher rank}
Let $H$ be a real vector space of finite dimension and $\Lambda =\R^r$ with its lexicographical order. 

Consider a symmetric bilinear form $\highinn{}{\cdot \,, \cdot} \colon H \times H \to \R^r$. 
We denote by $\highinn{j}{\cdot\,, \cdot}$ the $j$-th coordinate of $\highinn{}{\cdot\,,\cdot}$, so that $\highinn{}{\cdot\,,\cdot}$ has the form
\[\highinn{}{\cdot\,,\cdot} = \bigl(\highinn{1}{\cdot\,, \cdot}, \dots, \highinn{r}{\cdot\,, \cdot}\bigr).\] 
In case that $r=1$, we call $\highinn{}{\cdot \,, \cdot}$ a {\em scalar form}. 

We say that $\highinn{}{\cdot \,, \cdot}$ is {\em non-negative} if $\highinn{}{x\,, x} \succeq  0$ for all $x \in H$, and {\em positive} if $\highinn{}{x\,, x} \succ 0$ for all $x \in H \setminus \{0\}$.  In the case $r=1$, a positive scalar form is called a \emph{scalar product}. 

A symmetric bilinear form $\highinn{}{\cdot \,, \cdot} \colon H \times H \to \R^r$ defines a filtration on $H$ as follows. Consider the decreasing filtration $\Lambda^\bullet$ on $\Lambda = \R^r$ given in \eqref{eq:FiltrationRr}.
Setting
\begin{align}
\filter^j  &\coloneqq \left\{ x \in H \st \,\highinn{}{x \,, y} \in \Lambda^j \text{ for all $y \in H$}\right\}\\
&= \left\{ x \in H \st \,\highinn{}{x \,, y}_i = 0 \text{ for all $i <j$ and all $y \in H$}\right\}
\end{align}
we obtain a non-increasing filtration 
\[
\filter^\bullet\colon \qquad \filter^1 = H \supseteq \filter^2 \supseteq \filter^3 \supseteq \dots \supseteq \filter^r \supset \filter^{r+1} = (0)
\]
of $H$. We call $\filter^\bullet$ the {\em filtration induced by} $\highinn{}{\cdot \,, \cdot}$.

We have the following basic examples of bilinear forms and their induced filtrations.
\begin{example}[Euclidean products of higher rank] \label{ex:HigherRankEuclidean} Let $H = \R^r$. Then \[
\highinn{}{x\,, y} \coloneq \left(x_1y_1, x_2 y_2, \dots, x_r y_r\right), \qquad x, y \in H,
\] 
defines a positive bilinear form with values in $\Lambda = \R^r$. The filtration $\filter^\bullet$ on $H$ is given by 
\[
\filter^j = \left\{ x \in H \st \, x_i = 0 \text{ for all $i <j$}\right\}, \qquad j=1, \dots, r+1.\qedhere
\]
\end{example}

 \begin{example}[Euclidean products of higher rank associated to ordered partitions] \label{ex:HigherRankProjection} More generally, let $\pi=(\pi_1,\pi_2, \dots, \pi_r)$ be an ordered partition of a finite set $E$. (For the terminology of ordered partition, see~\cite[Section 2.2]{AN}.)
 Consider a real vector subspace $H \subseteq \R^E$. We define a positive bilinear form $\highinn{}{\cdot\,, \cdot} \colon H \times H \to \R^r$ as follows.

For each $j$, let $\projeuc_{j} \colon \R^E \to \R^{\pi_j}$ be the projection map. For any pair of vectors $x, y \in H$, set
 \[
 \highinn{}{x, y} \coloneqq \bigl( \innone{\pi_1}{\projeuc_{1}(x),\projeuc_{1}(x)}, \dots, \innone{\pi_r}{\projeuc_{r}(x),\projeuc_{r}(x)} \big)
 \]
 where $\innone{\pi_j}{\cdot \, , \cdot}$ is the Euclidean scalar product in $\R^{\pi_j}$.
 The induced filtration $\filter^\bullet$ on $H$ is given by
 \[
 \filter^j = \{ x \in H | \, \projeuc_i(x) = 0 \text{ for all $i<j$} \}, \qquad j = 1, \dots, r+1.
 \]
  Note that the previous example corresponds to setting $E = [r]$ and the ordered partition $\pi=(\pi_1, \dots, \pi_r)$ of $E = [r]$ with $\pi_j\coloneqq\{j\}$. 
\end{example}

\begin{example} \label{ex:PositiveButNotInner}
Let $H = \R^3$. Consider the bilinear form $\highinn{}{\cdot \,, \cdot} \colon H \times H \to \R^3$ given by
\[
\highinn{}{x , y} = \big ( x_1 y_1, x_2 y_2 + (x_1y_3 + x_3 y_1), x_3 y_3 \big ), \qquad x,y \in H.
\]
Then $\highinn{}{\cdot \, , \cdot}$ is positive. The filtration $\filter^\bullet \colon \filter^1 \supseteq \filter^2 \supseteq \filter^3  \supseteq \filter^{4}$ is given by
\begin{align*}
\filter^1 = H = \R^3, \qquad \filter^2 =  \{ \gamma \in H | \, \gamma_1 = 0 \} = \{0\} \times \R^2, \qquad \textrm{and} \qquad \filter^3 = \filter^4 = (0).\qedhere
\end{align*}
 
\end{example}
The following properties are immediate from the definition of the filtration $\filter^\bullet$.

\begin{prop} \label{prop:GradedForms} Let $\highinn{}{\cdot \,, \cdot} \colon H \times H \to \R^r$ be a symmetric bilinear form. For $j \in [r]$, the $j$-th component $\highinn{j}{\cdot \,, \cdot}$ of $\highinn{}{\cdot \,, \cdot}$ induces a symmetric bilinear form
\[
\highinn{j}{\cdot\,, \cdot} \colon \rquot{\filter^j}{\filter^{j+1}} \times \rquot{\filter^j}{\filter^{j+1}} \to \R
\]  
on the quotient $\rquot{\filter^j}{\filter^{j+1}}$.

If $\highinn{}{\cdot \,, \cdot}$ is non-negative on $H$, then $\highinn{j}{\cdot \,, \cdot}$ is non-negative on $\rquot{\filter^j}{\filter^{j+1}}$ for all $j$. If $\highinn{j}{\cdot \,, \cdot}$ is positive on $\rquot{\filter^j}{\filter^{j+1}}$ for all $j$, then $\highinn{}{\cdot \,, \cdot}$ is positive on $H$.
\end{prop}

\begin{defi}[Inner product of higher rank] \label{defi:InnerProduct}
 $-$ An \emph{inner product on $H$ with values in $\Lambda = \R^r$} is a symmetric bilinear form 
\[
\highinn{}{\cdot\,, \cdot} \colon H \times H \to \R^r
\]
such that for all $j \in [r]$, the induced form $\highinn{j}{\cdot\,, \cdot}$ on $\filter^j / \filter^{j+1}$ is positive, equivalently, if $\highinn{j}{\gamma, \gamma} >0$ for all $\gamma \in \filter^j \setminus \filter^{j+1}$. The pair $(H, \highinn{}{\cdot\,, \cdot})$ is called {\em an inner product space}.

\smallskip

$-$ For $j \in [r]$, the  $j$-th {\em graded piece} $\grm{}{j}H$ of $(H, \highinn{}{\cdot\,, \cdot})$ is the quotient
\[
\grm{}{j}H = \rquot{\filter^j}{\filter^{j+1}}
\]
endowed with the scalar product
\[\highinn{j}{\cdot\,, \cdot} \colon \grm{}{j}H\times \grm{}{j}H \to \R.
\]
We denote by $\proj_j \colon \filter^j \to \grm{}{j}H$ the projection map.

\smallskip

$-$ By a slight abuse of notation, we refer to inner products with values in $\R^r$, for arbitrary positive integer $r$, as \emph{inner products of higher rank} or simply as \emph{inner products} (although, strictly speaking, $r=1$ is allowed and, for larger $r$, their image is allowed to be a one-dimensional subspace of $\R^r$).\qedhere
\end{defi}

\begin{remark}
Note that with this definition, an inner product of rank one is isomorphic to a scalar product (up to fixing an isomorphism between the space of values and $\R$).
\end{remark}

\begin{remark}
We stress that every inner product $\highinn{}{\cdot \,, \cdot} \colon H \times H \to \R^r$ is necessarily a positive bilinear form by Proposition~\ref{prop:GradedForms}. 
On the other hand, as Example~\ref{ex:PositiveButNotInner} shows, not every positive bilinear form $\highinn{}{\cdot \,, \cdot} \colon H \times H \to \R^r$ is an inner product. In this example, the quotient $(\filter^2/\filter^3, \highinn{2}{\cdot \,, \cdot})$ is isomorphic to $\R^2$ with the semi-definite form $(x,y) = x_1y_1$.
\end{remark}

We have the following basic examples of higher rank inner products.
\begin{example}The bilinear forms in Examples~\ref{ex:HigherRankEuclidean} and~\ref{ex:HigherRankProjection} are inner products. In Example~\ref{ex:HigherRankProjection}, the graded piece $\grm{}{j} H$ is isomorphic to the image $\projeuc_{j}(\filter^j) \subseteq \R^{\pi_j}$ equipped with the induced Euclidean product from $\R^{\pi_j}$.
\end{example}

\begin{defi}[Quadratic form associated to an inner product] For an inner product $\highinn{}{\cdot\,, \cdot} \colon H \times H \to \R^r$, we denote by $\qf$ its {\em quadratic form}
\begin{align*}
\begin{array}{cccc}\qf = (\qf_1 , \dots, \qf_r ) \colon &H &\rightarrow &\R^r \\  & \gamma &\mapsto &\highinn{}{\gamma,\gamma}.\qedhere
\end{array}
\end{align*}
\end{defi}
\begin{prop} Let $\highinn{}{\cdot\,, \cdot} \colon H \times H \to \R^r$ be an inner product. Then,
\[
\filter^j = (\qf)^{-1}(\Lambda^j) = \left\{x \in H \st \qf_1(x)=\dots=\qf_{j-1}(x)=0\right\}.
\]
\end{prop}
\begin{proof} The inclusion $\filter^j\subseteq (\qf)^{-1}(\Lambda^j)$ follows from the definition of the filtration $\filter^\bullet$. The other inclusion follows from the positivity of $\qf_i$ on $\rquot{\filter^i}{\filter^{i+1}}$: if $x \in \filter^i \setminus \filter^{i+1}$ for $i<j$, then $\qf_i(x)>0$ which shows $x \notin (\qf)^{-1}(\Lambda^j)$.
\end{proof}

\begin{remark}\label{rem:filtration-realizability} Note that any filtration $\filter^1=H \supseteq \filter^2 \supseteq \dots \supseteq \filter^r \supseteq \filter^{r+1} =(0)$ of a finite dimensional real vector space can arise as the filtration associated to an inner product $\highinn{}{\cdot\,,\cdot} \colon H\times H \to \R^r$.
\end{remark}

\subsection{Higher rank inner products as limits of scalar products} \label{ss:PullbackFamilies} Consider a vector space $H$ and a symmetric bilinear form 
\begin{align*}
\highinn{}{\gamma, \eta}  = \bigl({\highinn{1}{\gamma, \eta}}, \dots, \highinn{r}{\gamma, \eta}\bigr), \qquad \gamma, \eta \in H, %
\end{align*}
on $H$ with values in $\Lambda = \R^r$. Let $\R_+^r=(0,+\infty)^r \subset \R^r$.

\begin{defi}[Pullback of a bilinear form] \label{defi:PullbackBilinearForm} Let $\underline L=(L_1, \dots, L_r)$ be a vector in $\R_+^r$. Then, the bilinear form $\innone{{\underline L}}{\cdot\,, \cdot} \colon H \times H \to \R$ defined by
\begin{align*}
\innone{{\underline L}}{\gamma, \eta} \coloneq L_1 {\highinn{1}{\gamma, \eta}} + \dots +L_r{\highinn{r}{\gamma, \eta}}, \qquad \forall\, \gamma, \eta\in H,
\end{align*}
is called the {\em pullback} by $\underline L$ of $\highinn{}{\cdot\,, \cdot}$. 
\end{defi}

We have the following fundamental result.

\begin{thm}\label{thm:pullback} Let $\highinn{}{\cdot\,, \cdot} \colon H \to \R^r$ be a symmetric bilinear form. The following statements are equivalent:

(i) $\highinn{}{\cdot\,, \cdot} \colon H \to \R^r$ is an inner product.

(ii) There exists a positive constant $C$ such that for any vector $\underline L = (L_1, \dots, L_r)$ of positive reals $L_1, \dots, L_r$ with
\[
L_j/L_{j+1}> C, \qquad j=1, \dots, r-1,
\]
the pullback $\innone{\underline L}{\cdot\,, \cdot}\colon H\times H \to \R$ is a scalar product. 
\end{thm}

\begin{proof} $-$ We first prove that (i) implies (ii). We proceed by induction on $r$. For $r=1$, there is nothing to prove. 

Suppose that the result holds for $r-1$. For the sake of a contradiction, assume that $\highinn{}{\cdot\,, \cdot} \colon H \to \R^r$ is an inner product which does not verify the claim.  Then, for any  $n \in \N$, there exists $\underline{L}_n=(L_{n,1}, \dots, L_{n,r})$ in $\R_+^r$ with $L_{n,j}/L_{n, j+1} > n$ for all $j$ and a point $\gamma_n\in H \setminus \{0\}$ such that
\[
\innone{\underline L_n}{\gamma_n, \gamma_n}\le0.
\]
Without loss of generality, we can assume that $\gamma_n$ lives on the unit sphere for a reference norm $\refnorm{\cdot}$ on $H$. Passing to a subsequence, we infer that the sequence $(\gamma_n)_n$ has a limit $\gamma$. There exists a $k \in [r]$ such that $\gamma$ belongs to $\filter^k \setminus \filter^{k+1}$. The $k$-th coordinate $\highinn{k}{\gamma, \gamma}$ is strictly positive, and, by continuity, there exists a small neighborhood $U$ of $\gamma$ in the unit sphere and a positive real $\delta>0$ such that 
\[\highinn{k}{\eta, \eta} > 2 \delta \sum_{j>k}^r \abs{\highinn{j}{\eta, \eta}}\]
 for all $\eta$ in $U$. For $n$ large, we then obtain that 
\begin{align*}
\sum_{j \ge k} L_{n,j} \, \highinn{j}{\gamma_n, \gamma_n} &\ge \frac 12 L_{n,k} \highinn{k}{\gamma_n, \gamma_n} + \sum_{j > k} L_{n,j} \Big( \delta { L_{n,k}}/{ L_{n,j}} - 1\Big) \abs{\highinn{j}{\gamma_n, \gamma_n}} \\
&\ge  \frac 12 L_{n,k} \highinn{k}{\gamma_n, \gamma_n}  > 0.
\end{align*}
As is easily verified, the first $k-1$ components of $\highinn{}{\cdot \,, \cdot}$ define an inner product  
\begin{align*}
[\cdot \,, \cdot]&\colon \rquot{H}{\filter^k}\times \rquot{H}{\filter^k} \to \R^{k-1}\\
[x \,, y] &\coloneqq ( \highinn{1}{x, y}, \dots, \highinn{k-1}{x, y}) \qquad \forall\, x, y \in\rquot{H}{\filter^k}
\end{align*} 
on $\rquot{H}{\filter^k}$ with values in $\R^{k-1}$. Applying the induction hypothesis to $[\cdot \,, \cdot]$, it follows that
\[
\sum_{j < k} L_{n,j} \, \highinn{j}{\gamma_n, \gamma_n} \ge 0
\]
for large $n$. This means that $\innone{\underline L_n}{\gamma_n, \gamma_n}>0$ for large $n$, which is a contradiction.

\smallskip

$-$ It remains to prove that (ii) implies (i). For the sake of a contradiction, assume that $\highinn{}{\cdot\,, \cdot}$ is not an  inner product. Fix a sequence of vectors $\underline{L}_n=(L_{n,1}, \dots, L_{n,r})$, $n \in \N$, satisfying that $L_{n,j}/L_{n, j+1} > n$ and, additionally, $L_{n,j}^2 / (L_{n, j+1} L_{n, 1})> n$ for all $j$. We prove that, for large $n$, the pullback form $\innone{\underline L_n}{\gamma_n, \gamma_n}$ is not an inner product on $H$.

By assumption, there exists a $j \in [r]$ and an element $x \in \filter^j \setminus \filter^{j+1}$ with $\highinn{j}{x,y} \le 0$. Since $x$ lies in $\filter^j\setminus \filter^{j+1}$, we can find an element $y \in H$ with $\highinn{j}{x,y} = -1$.
Consider now elements of the form $\lambda_n x + y$ for some $\lambda_n >0$. We prove that $\innone{\underline L_n}{\lambda_n x + y, \lambda_n x + y} < 0$  for a suitable choice of $\lambda_n$.

Observe first that for $n \in \N$ and any $\lambda_n > 0$, we have
\begin{align*}
\innone{\underline L_n}{\lambda_n x + y, \lambda_n x + y} = \sum_{k \le j} L_{n,k}\highinn{k}{y, y} + L_{n,j}(\lambda_n^2 \highinn{j}{x, x} - 2\lambda_n) + \sum_{k > j} L_{n,k} \highinn{k}{\lambda_n x + y, \lambda_n x + y}.
\end{align*}
Moreover, we can find a positive $\varepsilon >0$ such that 
\[\sum_{k > j} \abs{\highinn{k}{y, y}} + \abs{\highinn{k}{x, x}} + 2\abs{\highinn{k}{x, y}} < \varepsilon \qquad \textrm{and} \qquad \sum_{k \le j} \abs{\highinn{k}{y, y}}< \varepsilon.
\]
Choosing
\[
\lambda_n = \varepsilon L_{n,1}/L_{n,j},
\]
and recalling that $\highinn{j}{x, x} \le 0$, and $L_{n,j}^2 / (L_{n, j+1} L_{n, 1})> n$ by assumption, it follows that
\[
\innone{\underline L_n}{\lambda_n x + y, \lambda_n x + y} \le \varepsilon^3 \frac{L_{n,1}^2 L_{n, j+1}}{L_{n, j}^2}    - \varepsilon  L_{n,1} < 0
\]
for all large $n \in \N$. This completes the proof.
\end{proof}

We note the following immediate consequence, which will be important in the sequel.
\begin{cor} \label{cor:OullbackFamiliesInnerProducts} Let $\highinn{}{\cdot\,, \cdot} \colon H \to \R^r$ be an inner product. Suppose $\underline L_t = (L_{t,1}, \dots, L_{t,r})$, $t \in \R_+$, are vectors of positive reals with
\begin{equation} \label{eq:MultiScaleLimit}
\lim_{t \to \infty} \frac{L_{t,j}}{L_{t, j+1}} =   \infty, \qquad j=1, \dots, r-1.
\end{equation}
The pullbacks $\innone{\underline L_t}{\cdot\,, \cdot}\colon H\times H \to \R$, $t\in \R_+$, are scalar products for any large enough $t \in \R_+$.
\end{cor}
Given a family of vectors $\underline L_t \in \R_+^r$, $t \in \R_+$, satisfying \eqref{eq:MultiScaleLimit}, we call the resulting family of scalar products $\innone{\underline L_t}{\cdot\,, \cdot}$, $t \in \R_+$, a {\em pullback family} for the inner product $\highinn{}{\cdot\,, \cdot}$. The vectors $\underline L_t$, $t \in \R_+$, are called the \emph{parameters of the pullback family}.

Here and below, we slightly abuse the notation and ignore the fact that the pullbacks $\innone{\underline L_t}{\cdot\,, \cdot}$ are scalar products only for large $t \in \R_+$. Since we will be always interested in the asymptotic properties for large values of $t \in \R_+$ when considering pullback families, this will not be an issue.

\begin{remark}
We stress that Theorem~\ref{thm:pullback} states that a symmetric bilinear form $\highinn{}{\cdot\,, \cdot} \colon H \times H\to \R^r$ is an inner product exactly when the pullbacks $\innone{\underline L_t}{\cdot\,, \cdot}$ are eventually scalar products for {\em all families} $\underline L_t = (L_{t,1}, \dots, L_{t,r})$, $t \in \R_+$, satisfying \eqref{eq:MultiScaleLimit} (see however the next remark). 

If $\highinn{}{\cdot\,, \cdot}$ is not an inner product, it can still happen that the pullbacks $\innone{\underline L_t}{\cdot\,, \cdot}$ are scalar products for certain sequences $\underline L_t$, $t \in \R_+$, satisfying \eqref{eq:MultiScaleLimit} (see Section~\ref{ss:PathologyTorus} for an interesting example). On the other hand, if such a sequence $\underline L_t$, $t \in \R_+$, exists, then the bilinear form $\highinn{}{\cdot\,, \cdot}$ is necessarily positive. 
\end{remark}

\begin{remark} The proof of Theorem~\ref{thm:pullback} shows that if the pullbacks $\innone{\underline L_t}{\cdot\,,\cdot}$ are scalar products for some family $\underline L_t$, $t \in \R_+$, which satisfies the stronger condition that
\begin{equation} \label{eq:StrongSeparation}
\lim_{t \to \infty} \frac{L_{t,j}}{L_{t, j+1}} \cdot \frac{L_{t,j}}{L_{t,1}} = + \infty, \qquad \text{for all $j=1, \dots, r-1$},
\end{equation}
then $\highinn{}{\cdot\,, \cdot}$ is necessarily an inner product. This condition implies \eqref{eq:MultiScaleLimit} and merely means that the different scales $L_{t,1}, \dots ,L_{t,r}$ are ``very far" from each other.
\end{remark}


\subsection{The almost orthogonal decomposition} \label{ss:AlmostOrthogonalDecomposition} Let $\highinn{}{\cdot\,,\cdot} \colon H \times H \to \R^r$ be an inner product and $\filter^\bullet$ the induced filtration on $H$.

Consider two elements $\gamma \in \filter^{j} \setminus \filter^{j+1}$ and $\gamma' \in \filter^{j'} \setminus \filter^{j'+1}$. Then $\highinn{}{\gamma\,,\gamma'}$ lies in $\Lambda{\max\{j, j'\}}$, that is, it has the form
\[
\highinn{}{\gamma\,,\gamma'} = \Big (0, \dots, 0, \highinn{\max\{j, j'\}}{\gamma\,,\gamma'}, \dots, \highinn{r}{\gamma\,,\gamma'} \Big ).
\]
This motivates the following definition.

\begin{defi}[Almost orthogonality]
Two elements $\gamma, \gamma' \in H$ are called {\em almost orthogonal} and we write $\gamma \aperp\gamma'$ if $\highinn{\max\{j,j'\}}{\gamma, \gamma'} = 0$ where $j, j' \in [r]$ are the indices with $\gamma \in \filter^j \setminus \filter^{j+1}$ and $\gamma' \in \filter^{j'} \setminus \filter^{j'+1}$.

Two subsets $A, B \subset E$ are called {\em almost orthogonal} and we write $A \aperp B$, if any pair of elements $\gamma \in A$ and $\eta \in B$ are almost orthogonal. We use the notation $A \aplus B$ for $A+B$, indicating that $A$ and $B$ are almost orthogonal.  
\end{defi}
The terminology is justified as follows.
\begin{lem}\label{lem:orthogonality_vanishing} If $\gamma \in H$ is almost orthogonal to itself, then $\gamma = 0$. In particular, if $U, V \subset H$ are almost orthogonal subspaces, then $U \aplus V$ is a direct sum. That is, if $u+v = 0$ for $u\in U$ and $v \in V$, then $u=v = 0$.  
\end{lem}
\begin{proof} If $\gamma \neq 0$, then there exists $k \in [r]$ with $\gamma \in \filter^k \smallsetminus \filter^{k+1}$. Moreover,  $\highinn{k}{\gamma, \gamma} >0$ since $\highinn{}{\cdot, \cdot}$ is an inner product. On the other hand, if $\gamma$ is orthogonal to itself, then $\highinn{k}{\gamma, \gamma} =0$. We have arrived at a contradiction.
\end{proof}

As we discuss next, the notion of almost orthogonality allows to canonically decompose $H$ in terms of its graded pieces $\grm{}{j}H$, $j \in [r]$.

\begin{lem}[Lifting lemma] \label{lem:LiftingLemma}Let $\highinn{}{\cdot\,,\cdot} \colon H \times H \to \R^r$ be an inner product. Then for $j\in [r]$, there exists a unique linear map $\proj_j^* \colon \grm{}{j}H \to \filter^j$ with the following properties:
\begin{enumerate}
\item $\proj_j \circ \proj_j^* = \id$ and, in particular, $\proj_j^* \colon \grm{}{j}H \hookrightarrow \filter^j$ is an embedding.
\item The image of $\proj_j^*$ is almost orthogonal to $\filter^{j+1}$.
\end{enumerate}
The inner product space $(H, \highinn{}{\cdot, \cdot})$ has the almost orthogonal decomposition
\[H = \proj_1^*\left(\grm{}{1}H\right) \aplus \proj_2^*\left(\grm{}{2}H\right) \aplus \dots \aplus \proj_r^*\left(\grm{}{r}H\right).\]
That is, $H$ is the direct sum of the subspaces $\proj_j^*(\grm{}{j}(H))$, $j \in [r]$, and they are pairwise almost orthogonal.
\end{lem}
\begin{defi} We refer to the embedding $\proj_j^*\colon \grm{}{j} H \hookrightarrow H$ as the \emph{canonical lifting} of $\grm{}{j}H$ to $H$, and denote its image by $H_j \coloneqq \proj_j^*( \grm{}{j} H)$. The resulting almost orthogonal decomposition is denoted by
\[H = H_1\aplus H_2 \aplus \dots \aplus H_r. \qedhere\]
\end{defi}

\begin{proof}[Proof of Lemma~\ref{lem:LiftingLemma}]
Fix $j \in [r]$ and let $\gamma\in \grm{}{j}H$. We show that there exists a unique element $\theta \in \filter^j$ such that $\proj_j(\theta) =\gamma$, and $\theta$ is almost orthogonal to $\filter^{j+1}$. This defines a unique linear map $\proj_j^*$ with the above properties. 

To prove the uniqueness of $\theta$, let $\theta$ and $\theta'$ be two elements of $\filter^j$ such that $\proj_j(\theta) =\proj_j(\theta') = \gamma$, and both are almost orthogonal to $\filter^{j+1}$. Then the difference $\theta-\theta'$ belongs to $\filter^{j+1}$ and is moreover almost orthogonal to $\filter^{j+1}$. Lemma~\ref{lem:orthogonality_vanishing} then ensures that $\theta =\theta'$.

To prove the existence, we proceed inductively, and show that for each $k=0, \dots, r-j$, we can find an element $\theta_{j+k}\in \filter^j$ such that 
\begin{itemize}
\item[($i$)] $\proj_j(\theta_{j+k}) = \gamma$, and 

\item[($ii$)] for each $h=1, \dots, k$, we have $\highinn{j+h}{\theta_{j+k}, \eta} = 0 \textrm{ for all  $\eta\in \filter^{j+h}$}.$
\end{itemize}
Setting $\theta = \theta_{r}$, we then obtain the desired element. 

For $k=0$, choose any $\theta_j \in \filter^j$ with $\proj_j(\theta_j) =\gamma$. Then $\theta_j$ verifies ($i$) and ($ii$) is empty.

Assume we have found $\theta_{j+k}$ for some $k< r-j$. Define the linear form 
\begin{align*}  f_{j+k}&\colon \filter^{j+k+1} \to \R\\
\eta &\longmapsto \highinn{j+k+1}{\theta_{j+k}, \eta}.
\end{align*}
Since $\highinn{j+k+1}{\cdot \,, \cdot}$ is an inner product on $\grm{}{j+k+1}H$, there exists an element $\alpha \in \filter^{j+k+1}$ with
\[f_{j+k}(\eta) = \highinn{j+k+1}{\alpha, \eta}, \qquad \forall \eta\in \filter^{j+k+1}.  \]
Let $\theta_{j+k+1} \coloneqq \theta_{j+k} - \alpha$ and note that $\proj_j(\theta_{j+k+1}) = \proj_j(\theta_{j+k}) =\gamma$. Since $\alpha \in \filter^{j+k+1}$, we also have $\highinn{j+h}{\alpha, \eta}=0$ for any $\eta \in \filter^{j+h}$ for all $h=1, \dots, k$. Finally,
\[\highinn{j+k+1}{\theta_{j+k+1}, \eta} = \highinn{j+k+1}{\theta_{j+k}, \eta} - \highinn{j+k+1}{\alpha, \eta}= 0, \qquad \forall \eta\in \filter^{j+k+1},\]
by construction. Thus, $\theta_{j+k+1}$ satisfies the properties ($i$) and ($ii$). This finishes the proof of existence and uniqueness of the embedding $\proj_j^* \colon \grm{}{j}H \hookrightarrow H$.

The properties of $\proj_j^*$ clearly imply that the subspaces $\proj_i^*(\grm{}{j}H)$, $j \in [r]$, are pairwise almost orthogonal. It remains to prove that $H$ is the direct sum of these subspaces. Since $\dim(H) = \sum_{j} \dim(\grm{}{j}H)$, the desired claim follows from Lemma~\ref{lem:orthogonality_vanishing}.
\end{proof}

\section{Higher rank Voronoi decompositions} \label{sec:HigherRankVoronoiDecompositions}
Throughout this section, let $H$ be a finite dimensional vector space and $\highinn{}{\cdot\, ,\cdot}\colon H\times H \to \R^r$ an inner product with values in $\Lambda = \R^r$.

Let $S \subset H$ be a discrete subset of $H$, that is, $S$ has no accumulation point in $H$ endowed with its Euclidean topology. We define the \emph{Voronoi cell} of $\gamma \in S$ as
\begin{equation} \label{eq:DefinitionVoronoiCell}
\Vor_S(\gamma) \coloneqq \left\{x\in H \, \st\, \highinn{}{x-\gamma, x- \gamma} \preceq \highinn{}{x-\eta, x-\eta} \,\textrm{ for all } \eta \in S\right\}.
\end{equation}
If $S$ is understood from the context, we simply write $\Vor(\gamma)$.

The analogy with the case $r = 1$ and the relation to degenerating scalar products (see Theorem~\ref{thm:pullback} ) suggest to interpret $\highinn{}{x, x}$ as the (vector-valued) length of an element $x \in H$. The Voronoi cell of $\gamma$ is precisely the set of all points $x \in H$ whose closest point in $S$ is $\gamma$.

\subsection{Basic properties and pathologies} \label{ss:Pathologies} In the classical setting of a scalar product, that is,  $r=1$, the Voronoi cells $\Vor_S(\gamma)$, $\gamma \in S$, are convex, closed and have non-empty, mutually disjoint interiors. If $S \subset H$ is finite or $S= \L$ is a lattice in $H$,  then each Voronoi cell is a polyhedron, that is, an intersection of finitely many half-spaces in $H$. For a general discrete subset $S \subset H$, the Voronoi cells are generalized polyhedra (see Section~\ref{ss:BasicNotations} for the terminology). Moreover, they provide a decomposition $H = \bigcup_{\gamma \in S} \Vor(\gamma)$ of $H$ into (generalized) polyhedra, which is called the \emph{Voronoi decomposition}, see e.g.~\cite[Page 465]{Gruber}.

These properties can drastically fail for inner products of higher rank, as shown by the following examples.

\begin{example}\label{ex:pathology1} Consider $H =\R^2$, $\Lambda=\R^2$ with the inner product from Example~\ref{ex:HigherRankEuclidean},
\begin{align*}
\highinn{}{\cdot\,,\cdot} \colon H\times H &\to \R^2 \\
\Bigl((x_1, x_2), (y_1, y_2)\Bigr) &\mapsto (x_1\cdot y_1 , x_2\cdot y_2).
\end{align*}
Let  $S = \{0\}\sqcup \bigl\{(a_n,b_n) \st n\in \Z_{\ge 1} \bigr\} \sqcup \bigl\{(c_n, d_n) \st n \in \Z_{\le- 1}\bigr\}$ be a sequence of points in $\R^2$ with coordinates verifying
\begin{align*}
a_1 > a_2 > a_3 > \dots >0,& \qquad  0<b_1<b_2< b_3 < \dots, \quad \textrm{and}\\
 c_{-1} < c_{-2} < c_{-3}< \dots <0, & \qquad  0>d_{-1}>d_{-2}> d_{-3} > \dots
\end{align*}
such that
\[\lim_{n\to +\infty}a_n = \lim_{n \to -\infty} c_n = 0, \qquad \lim_{n\to +\infty}b_n = +\infty, \qquad \qquad \lim_{n \to - \infty} d_n = -\infty.\]

 The Voronoi decomposition of the sequence is depicted in Figure~\ref{fig:voronoi-path1}.
 \begin{figure}[!t]
\centering
   \scalebox{.32}{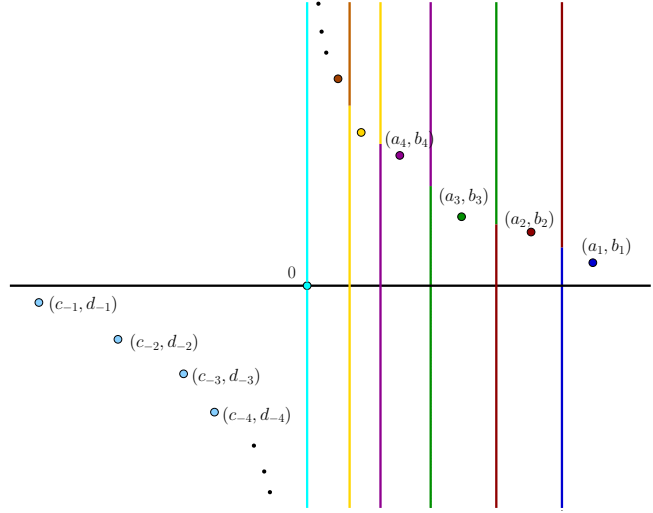}
\caption{A set $S$ and its Voronoi decomposition as described in Example~\ref{ex:pathology1}. Only the Voronoi cells of the origin and points $(a_n,b_n)$, $n\in\N$, are shown. The colored parts of the straight lines belong to the Voronoi cells of the points with the same color.}
\label{fig:voronoi-path1}
\end{figure}
In particular, we have $\Vor(0) = \{0\}\times \R$, and for $n \geq 1$, we have
\[\Vor((a_n,b_n)) =  \big\{ x \in \R^2 \st A_n < x_1 < A_{n-1}\big\} \sqcup  \{  A_n \} \times (-\infty, B_n]  \sqcup \{ A_{n-1} \} \times [B_{n-1}, + \infty)\]
with 
\[A_0 =\infty, \qquad A_n = \frac{a_n+a_{n+1}}2, \qquad  \textrm{and } \qquad B_n= \frac{b_n + b_{n+1}}2\qquad  \forall n\geq 1.\]
The cell $\Vor((c_n,d_n))$ has a similar description.

In particular, the Voronoi cell of $\Vor(0)$ has empty interior in $H$. Moreover, the Voronoi cells $\Vor( (a_n,b_n ))$, $n \ge 1$, and  $\Vor( (c_n,d_n ))$, $n \le -1$, are not closed. 
\qedhere
\end{example}

We will be particularly interested in the case where $S=\L$ is a lattice of full dimension in $H$. In this case, we have
\[
\Vor(\gamma) = \Vor(0)+\gamma, \qquad  \gamma \in \L.
\]
In the classical setting, when $r=1$, the Voronoi cell $\Vor(0)$ is a polytope of full dimension. Moreover, the cells $\Vor(\gamma)$, $\gamma \in S$, provide an $\L$-periodic tiling of $H$ by polytopes.

However, even in the case of lattices, Voronoi cells of higher rank inner products behave different from their classical counterparts.

\begin{example}\label{ex:pathology2} Consider $H =\R^2$ with the $(\R^2$)-valued inner product from Example~\ref{ex:HigherRankEuclidean}. Let $\L$ be the lattice of rank two generated by the vectors $(\pi, 1)$ and $(1,0)$, so that 
\[\L= \left\{(a\pi+b, a)\st a,b\in \Z\right\}.\]
 Then, we have 
\[\Vor(0) = \{0\}\times \R\]
and more generally 
\[\Vor\left((a\pi+b, a)\right) = \{(a\pi+b)\} \times \R.\]
Thus, the Voronoi cells $\Vor(\gamma)$, $\gamma \in \L$, are unbounded and have empty interiors. Note that they do not provide a decomposition of $H$. For instance, a point $x \in H$ of the form $x = (x_1, x_2)$ with $x_1 \in \Q \setminus \Z$ does not belong to any of the Voronoi cells.
\end{example}

\begin{example}\label{ex:pathology3} Figure~\ref{fig:voronoi} shows an example with $H$ of dimension two, the space of values $\Lambda =\R^2$, and $S=\L$ a lattice of full rank in $H$. The Voronoi cell $\Vor(0)$ is bounded and has non-empty interior, but is not closed. The closures  $\overline \Vor(\gamma)$, $\gamma \in \L$, provide a polyhedral tiling of $H$.
\end{example}

The following proposition summarizes the basic properties of higher rank Voronoi cells.

\begin{prop}\label{prop:general_voronoi} Let $\highinn{}{\cdot\, ,\cdot}\colon H\times H \to \R^r$ be an inner product and $S$ a discrete subset of $H$. Then for any $\gamma\in S$, the Voronoi cell $\Vor(\gamma)$ is a non-empty convex set. Moreover, the interiors of the Voronoi cells are mutually disjoint.

In general, $\Vor(\gamma)$ is not closed and can have empty interior. The union of the Voronoi cells can be a proper subset of $H$. Moreover, even in case that $S = \L$ is a lattice of full rank in $H$, the Voronoi cell $\Vor(\gamma)$ can be unbounded.
\end{prop}

Suppose that $0 \in S$. A point $x$ belongs to the Voronoi cell $\Vor(0)$ of $0 \in S$ exactly when
\[
 \highinn{}{x,x}  \preceq  \highinn{}{x - \eta, x-\eta} \qquad \text{ for all $\eta \in S$}.
\]
For $\eta \in S$, the above inequality can be written as
\[
2 \highinn{}{x,\eta} \preceq \highinn{}{\eta,\eta}.
\]
This formulation involves the $\R^r$-valued linear form $2 \highinn{}{\cdot,\eta}$  on $H$ and the constant $\highinn{}{\eta,\eta} \in \R^r$. The set of points $x \in H$ satisfying this inequality may be viewed as a higher rank version of a half-space in $H$.
 
\begin{proof}[Proof of Proposition~\ref{prop:general_voronoi}]
Obviously, $\Vor(\gamma)$ is non-empty as it contains $\gamma$. In proving that $\Vor(\gamma)$ is convex and $\interior {\Vor(\gamma)} \cap \interior{ \Vor(\eta)} = \emptyset$ for $\gamma \neq \eta$, without loss of generality, we can assume that $\gamma=0$. 
For $x, y \in \Vor(0)$, $t\in [0,1]$, and $\eta\in S$, we have 
\[2 \highinn{}{tx+(1-t)y,\eta} = 2t\highinn{}{x,\eta} + 2(1-t) \highinn{}{y,\eta} \preceq t \highinn{}{\eta,\eta} +(1-t)\highinn{}{\eta,\eta}=\highinn{}{\eta,\eta}\] 
and hence $\Vor(0)$ is convex.

Next, we prove that the interiors of $\Vor(0)$ and $\Vor(\eta)$, $\eta \neq 0$ are disjoint. Suppose that $x$ belongs to both $\Vor(0)$ and $\Vor(\eta)$. Then, $\highinn{}{x,x} = \highinn{}{x-\eta,x-\eta}$ and hence $2 \highinn{}{x,\eta} = \highinn{}{\eta,\eta}$. On the other hand, if $x$ lies in the interior of $\Vor(0)$, then $x + \varepsilon \eta$ belongs to $\Vor(0)$ for some small $\varepsilon >0$. In particular, $2\highinn{}{x+\varepsilon\eta, \eta} \preceq \highinn{}{\eta,\eta}$. This however implies that $\highinn{}{\eta,\eta} = 2 \highinn{}{x,\eta}  \preceq (1-2 \varepsilon) \highinn{}{\eta,\eta}$, a contradiction.

Examples~\ref{ex:pathology1} and \ref{ex:pathology2} demonstrate that $\Vor(\gamma)$ can be unbounded, non-closed and have empty interior. In Example~\ref{ex:pathology2}, the Voronoi cells do not cover $H$ and are unbounded.
\end{proof}

\subsection{Admissible discrete sets and their Voronoi decomposition} \label{ss:AdmissibleDiscreteSetsAndVoronoiDecomposition} The above examples suggest that discrete sets in the Euclidean topology on $H$ are not well compatible with higher rank inner products. In this section, we formulate a notion of discreteness in the higher rank setting, and prove that this leads to well-behaving Voronoi decompositions.

\subsubsection{Admissible discrete sets} Let $\highinn{}{\cdot \, , \cdot} \colon H \times H \to \R^r$ be an inner product with values in $\R^r$. Let $H = H_1 \aplus \dots \aplus H_r$ be the almost orthogonal decomposition. We write elements $x \in H$ as $x=x_1+\dots+x_r$ with $x_j \in H_j$, $j= 1, \dots, r$.
 
 Consider a subset $S\subset H$. For each $j \in [r]$ and a tuple $z = (z_1, \dots, z_{j-1})$, with $z_i \in H_i$ for $i=1,\dots, j-1$,  define the subset  $S_{/z} \subseteq S$ as
\begin{equation} \label{eq:SOver}
 S_{/z} \coloneqq \left\{ \gamma=\gamma_1 + \dots + \gamma_r \in S \st \gamma_i = z_i \text{  for $i=1, \dots, j-1$} \right\}.
 \end{equation}
Moreover, we define the translated set 
 \begin{equation} \label{eq:Translation}
 TS(z) \coloneqq S_{/z}  - z_1 - \dots -z_{j-1}.
 \end{equation}
We have $TS(z) \subset \filter^j$. In the special case $j=1$, we have the empty vector $z = \varnothing$ and recover
 \[
 S_{/\varnothing} = TS (\varnothing) = S \subset \filter^1.
\]
\begin{defi}[Admissible discrete sets]
A subset $S \subset H$ is called {\em admissible discrete} for the inner product $\highinn{}{\cdot\,,\cdot}$, if for all $j \in [r]$ and all tuples $z= (z_1, \dots, z_{j-1})$, with $z_i \in H_i$ for $i=1,\dots, j-1$, the projection $\proj_j(TS(z))$ is a discrete set in $\grm{}{j} H$.
\end{defi}
Example~\ref{ex:pathology2} shows that not all discrete sets are admissible discrete. However, the converse is always true.

\begin{prop}\label{prop:discrete-discrete} Every admissible discrete set $S$ in an inner product space $(H,  \highinn{}{\cdot \, , \cdot})$ is discrete with respect to the Euclidean topology on $H$. 
\end{prop}
\begin{proof}  Assume that this is not the case and let $x$ be an accumulation point of $S$. Then, there exists a sequence $(x^n)_{n \in \N}$ with $ x^n\in S \setminus \{x\}$ for all $n$, which converges to $x$. Proceeding by induction on $j\in [r]$, we show that the equality $x^n_i =x_i$ hold for all $i=1, \dots, j$, and $n \in \N$ large enough. Setting $j=r$, this shows that $x^n = x$ for all large $n \in \N$, contradicting the assumption.

 For $j=1$, we use the discreteness of the set $\proj_1(S)$ to infer that the projections $\proj_1(x^n)$ are eventually equal to $\proj_1(x)$. This implies that $x^n_1 = \proj_1^*(\proj_1(x^n)) = \proj_1^*(\proj_1(x)) =x_1$, as required. Assume that the statement holds for $j-1$ and let $z = (x_1, \dots, x_{j-1})$. We use the discreteness of the set $\proj_j(TS(z))$ to infer that for large enough $n$, we have 
 \[
 \proj_j(x^n_j) = \proj_j(x^n - x^n_1-\dots -x^n_{j-1}) = \proj_j(x_j).
 \]
Applying $\proj_j^*$, we infer that $x^n_j = x_j$ for large enough $n$, as required. \end{proof}

\begin{remark}
In Section~\ref{ss:HybridTopology}, we introduce a new topology on $H$ whose discrete sets are precisely the admissible discrete sets for the inner product $\highinn{}{\cdot \, , \cdot}$. Since this topology is coarser than the Euclidean topology on $H$, this gives an alternative proof of Proposition~\ref{prop:discrete-discrete}.
\end{remark}

\begin{remark}[Admissible discrete sets for a filtration] \label{rem:AdmissibleDiscreteFiltration} Admissible discrete sets can equivalently be defined by using the filtration $\filter^\bullet$ induced by the inner product $\highinn{}{\cdot \, \cdot}$.

More generally, if $H$ is a vector space of finite dimension and $\filter^\bullet : \filter^1 = H \supseteq \dots \supseteq \filter^r \supseteq \filter^{r+1}=(0)$ a non-increasing filtration, then we say that a subset $S\subset H$ \emph{is admissible discrete for $\filter^\bullet$} if recursively the following properties (1) and (2) hold...
\begin{itemize}
\item[(1)] The projection $S_1\coloneqq \proj_1(S)$ of $S$ into $\grm{}{1}H =\rquot{\filter^1}{\filter^2}$ is a discrete set.  
\item[(2)] For each $\gamma \in S$, the subset $TS(\gamma) \coloneqq S_{/\proj_1(\gamma)} -\gamma$ of $\filter^2$ is admissible discrete in $\filter^2$ for the induced filtration $\filter^2 \supseteq \filter^3\supseteq \dots \supseteq \filter^r \supseteq \filter^{r+1} =(0)$.
\end{itemize}
Here, for $a \in S_1$, we set 
\begin{align}\label{eq:SOver-filter}
S_{/a} \coloneqq \proj_1^{-1}(a) = \left\{\theta \in S \st \proj_1(\theta) = a\right\} \subset S.
\end{align}
 One can show that a subset $S \subset H$ is admissible discrete for an inner product  $\highinn{}{\cdot \, , \cdot}$ if and only if it is admissible discrete for the filtration $\filter^\bullet$  induced by $\highinn{}{\cdot \, , \cdot}$. In particular, if two inner products induce the same filtration $\filter^\bullet$ on $H$, then they have the same admissible discrete sets.
\end{remark}

\subsubsection{Voronoi decomposition induced by an admissible discrete set} We have the following Voronoi decomposition theorem for admissible discrete sets.
\begin{thm}\label{thm:voronoi-discrete} Let $\highinn{}{\cdot \,, \cdot} \colon H \times H \to \R^r$ be an inner product and $S \subset H$ an admissible discrete set. Then, the sets $\Vor_S(\gamma)$, $\gamma\in S$, cover the full space, that is, $H = \bigcup_{\gamma \in S} \Vor_S(\gamma)$.
\end{thm}
We call the above decomposition $H = \bigcup_{\gamma \in S} \Vor_S(\gamma)$ the {\em Voronoi decomposition} associated to the inner product $\highinn{}{\cdot \,, \cdot} $ and the admissible discrete set $S$.

\smallskip

In the proof, we use the almost orthogonal decomposition $H= \aplus_{j=1}^r H_j$ and decompose each element $x \in H$ as $x = x_1 + \dots + x_r$ with $x_j \in H_j$.

\begin{proof} We have to prove that for every $x \in H$, there exists an element $\gamma \in S$ such that $\qf(x-\gamma) \lexeq \qf(x-\eta)$ for all $\eta \in S$. Since $\lexeq$ is a total order on $\R^r$, it suffices to construct a finite set $\mathbb{S} \subset S$ with the property that for every $\eta \in S \setminus \S$, there exists an element $\gamma \in \S$ with $\qf(x-\gamma) \lexst \qf(x-\eta)$. 

In order to prove this, we proceed by induction on the dimension $r$ of $\R^r$. For $r = 1$, we are in the case of a scalar product and the claim is trivial.

Let $H'\coloneqq \rquot{H}{\filter^r}$ and denote by $\proj \colon H \to H'$ the projection map. The first $r-1$ components of $\highinn{}{\cdot \,, \cdot}$ define an inner product with values in $\R^{r-1}$
\begin{align*}
[\cdot \,, \cdot]&\colon H' \times H' \to \R^{r-1}\\
[x \,, y] &\coloneqq ( \highinn{1}{x, y}, \dots, \highinn{r-1}{x, y}) \qquad \forall\, x, y \in H'.
\end{align*}  
The almost orthogonal decomposition $H'= \aplus_{j=1}^{r-1} H_j'$ of $H'$ is given by $H_j' = \proj(H_j)$. Moreover, the projection $S ' \coloneqq \proj(S)$ of $S$ in $H'$ is an admissible discrete set in $H'$.

Using the induction hypothesis for the inner product $[\cdot \,, \cdot]$, the point $\proj(x)$, and the set $S'$, we infer the existence of a finite, non-empty subset $\S' \subset S'$ with the property that for any $\eta' \in S' \setminus \S'$, there exists an element $\gamma' \in \S'$ with $[\proj(x)-\gamma',  \proj(x)-\gamma'] \lexst [\proj(x)-\eta',  \proj(x)-\eta'] $.

This implies that for $\gamma \in \proj^{-1}(\S')$ and $\eta \in S \setminus \proj^{-1}(\S')$, we have
\begin{align*}
\qf(x - \gamma) = \Big([\proj(x - \gamma), \proj(x - \gamma)] , \qf_r(x - \gamma)\Big) \lexst \Big([\proj(x - \eta), \proj(x-\eta)] , \qf_r(x - \eta) \Big) 
= \qf(x - \eta).
\end{align*}
Moreover, setting $S_u = \left\{\gamma \in S \st \proj(\gamma) =u\right\}$, we can write $\proj^{-1}(\S')$ as a finite union
\[
\proj^{-1}(\S') = \bigcup_{u \in \S'} S_u.
\]
It thus suffices to construct for $u \in \S'$ a finite subset $\S_u \subset S_u$ with the following property: 
\begin{align} \label{eq:prop-finiteness}
\textrm{\emph{For every $\eta \in S_u \setminus \S_u$, there exists a $\gamma \in \S_u$ with $\qf(x - \gamma) \lexst \qf(x - \eta)$.}}
\end{align}
 Indeed, setting $\mathbb{S} = \bigcup_{u \in \S'} \S_u$, we obtain a set with the desired property. 

\smallskip

Given $u \in \S'$, we write it as $u = u_1 + \dots u_{r-1}$ with $u_j \in H'_j$. For $j=1, \dots, r-1$,  $\proj$ is injective on $H_j$ and $\proj(H_j) = H'_j$. Hence, there exists a unique element $z_j \in H_j$ with $\proj(z_j) = u_j$. Considering the tuple $(z_1, \dots, z_{r-1})$, we can identify
\begin{align*}
S_u &= \left\{\gamma \in S \st \proj(\gamma)_i = u_i \, \textrm{ for  } i=1, \dots, r-1 \right\}  \\
&= \left\{\gamma \in S \st \gamma_i = z_i \,\textrm{ for  } i=1, \dots, r-1\right\} \\
&=S_{/(z_1, \dots z_{r-1})}.
\end{align*}
Since $S$ is admissible discrete in $H$, we infer that the set $A = \left \{\gamma_r \st \gamma \in S_{/(z_1, \dots z_{r-1})}\right\}$ is discrete in $H_r$. Thus, there exists a finite, non-empty subset $I_u \subset A$ with
\[
\forall\,   b \in A \setminus I_u \,, \,\exists\, a \in I_u \qquad \qf_r(x_r - a) < \qf_r(x_r - b).
\]
Then $\S_u \coloneqq \left \{ \gamma \in S_{/(z_1, \dots z_{r-1})} \st \gamma_r \in I_u\right\}$ is a finite, non-empty subset of $S_u$. We show that this set verifies the property \eqref{eq:prop-finiteness}. 

Let  $\eta\in S_u \setminus \S_u$.  Then, we choose $\gamma \in \S_u$ with the property that $\qf_r(x_r-\gamma_r) < \qf(x_r-\eta_r)$. Using the almost orthogonality,
\begin{align*}
\qf_r(x-\gamma) &= \qf_r \big ( x_r-\gamma_r + \sum_{j=1}^{r-1} x_j - z_j  \big ) =  \qf_r(x_r-\gamma_r) + \qf_r  \big (  \sum_{j=1}^{r-1} x_j - z_j  \big )  \\
& <   \qf_r(x_r-\eta_r  ) + \qf_r  \big (  \sum_{j=1}^{r-1} x_j - z_j  \big) = \qf_r(x-\eta).
\end{align*}
Since $\qf_j(x-\gamma) = \qf_j(x-\eta)$ for $j=1, \dots, r-1$, it follows that $\qf(x-\gamma) \lexst \qf(x-\eta)$. Thus, $\S_u$ has the desired properties and the proof is complete.
\end{proof}


\section{Structure of Voronoi cells} \label{sec:StructureVoronoiCells}
Let $\highinn{}{\cdot \, , \cdot} \colon H \times H \to \R^r$ be an inner product with values in $\Lambda = \R^r$ and $S \subset H$ an admissible discrete subset. Consider the associated Voronoi decomposition $\bigcup_{\gamma\in S} \Vor_S(\gamma)=H$ (see Theorem~\ref{thm:voronoi-discrete}). In this section, we provide a description of the closed Voronoi cells $\overline{\Vor}_{S}(\gamma)$. 

\smallskip

Consider the filtration $\filter^\bullet = (\filter^j)_{j}$ on $H$ induced by $\highinn{}{\cdot\, ,\cdot}$. Recall that the $j$th graded piece $\grm{}{j}H = \rquot{\filter^j}{\filter^{j+1}}$, $j \in [r]$, is endowed with the scalar product $\highinn{j}{\cdot\, ,\cdot} \colon \grm{}{j}H \times \grm{}{j}H \to \R$. We have the almost orthogonal decomposition $H = \aplus_{j=1}^r H_j$.  In the following, we decompose each $x \in H$ as $x = x_1 + \dots +x_r$ with $x_j \in H_j$, $j=1, \dots, r$.  Note that $\highinn{j}{\cdot\, ,\cdot}$ is a scalar product on $H_j$. Moreover, we have $H_j \subseteq \filter^j$ and $H_j \simeq \grm{}{j}H$, where the isomorphism is given by $\proj_j^*$ in the Lifting Lemma~\ref{lem:LiftingLemma}.

Consider an element $\gamma =\gamma_1 + \dots + \gamma_r$ in $S$. For $j \in [r]$, define the subset $S_{\gamma, j} \subset H_j$ by
\[
S_{\gamma, j} = \left\{ \theta_j \st \theta \in S \text{ and } \theta_i = \gamma_i \text{ for } i \in\{1, \dots, j-1\} \right\}. 
\]
Note that $S_{\gamma, 1} = \left\{\theta_1 \st \theta \in S\right\}$ does not depend on $\gamma$.

Using the tuple $(\gamma_1, \dots, \gamma_{j-1})$ and the notation from \eqref{eq:SOver} and \eqref{eq:Translation}, we have
\[
S_{\gamma, j} = \left\{ \theta_j \st \theta \in S_{/(\gamma_1, \dots, \gamma_{j-1})} \right\} = \proj_j^\ast(\proj_j(TS(\gamma_1, \dots, \gamma_{j-1}))).
\]
Since $S$ is admissible discrete, $S_{\gamma,j}$ is a discrete subset of $H_j$ with $\gamma_j \in S_{\gamma,j}$. Let 
\[
V_{\gamma, j} \coloneqq \Vor_{S_{\gamma,j}}(\gamma_j) \subseteq H_j
\]
 be the Voronoi cell of $\gamma_j$ in $H_j$ associated to $S_{\gamma,j}$ and the scalar product $\highinn{j}{\cdot\,,\cdot}$.

The main result of this section reads as follows.

\begin{thm}[Almost orthogonal decomposition theorem for Voronoi cells]\label{thm:VoronoiCellDecomposition-discrete}  Notations as above, let $S \subset H$ be an admissible discrete set and $\gamma \in S$. Then, the closure of $\Vor_{S}(\gamma)$ has the almost orthogonal decomposition
\[
\overline \Vor_{S}(\gamma) = V_{\gamma,1} \aplus \dots \aplus V_{\gamma,r}.
\]
That is, $\overline \Vor_{S}(\gamma)$ coincides with the Minkowski sum of the $V_{\gamma,j}$, $j \in [r]$, and these sets are pairwise almost orthogonal.

It follows that $\overline \Vor_{S}(\gamma)$ is a (generalized) polyhedron with a non-empty interior.  
\end{thm}

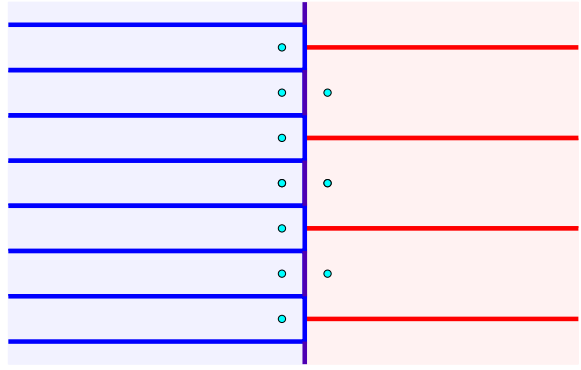
\begin{figure}[!t]
\centering
\begin{tikzpicture}[scale=.6]

\filldraw[red!5!white] (6.5, 4) -- (.5, 4) -- (.5, 3) -- (6.5, 3) -- (6.5, 4);
\filldraw[red!5!white] (6.5, -4) -- (.5, -4) -- (.5, -3) -- (6.5, -3) -- (6.5, -4);

\filldraw[red!5!white] (6.5, 1) -- (.5, 1) -- (.5, -1) -- (6.5, - 1) -- (6.5, 1);
\draw[line width=0.6mm, red] (6.5, 1) -- (.5, 1) -- (.5, -1) -- (6.5, - 1);

\filldraw[red!5!white] (6.5, 3) -- (.5, 3) -- (.5, 1) -- (6.5, 1) -- (6.5, 3);
\draw[line width=0.6mm, red] (6.5, 3) -- (.5, 3) -- (.5, 1) -- (6.5, 1);
\filldraw[red!5!white] (6.5, -3) -- (.5, -3) -- (.5, -1) -- (6.5, -1) -- (6.5, -3);
\draw[line width=0.6mm, red] (6.5, -3) -- (.5, -3) -- (.5, -1) -- (6.5, -1);

\filldraw[blue!5!white] (-6, {4}) -- (.5, {4}) -- (.5, {3.5}) -- (-6, {3.5}) -- (-6, {4});

\filldraw[blue!5!white] (-6, {-4}) -- (.5, {-4}) -- (.5, {-3.5}) -- (-6, {-3.5}) -- (-6, {-4});

\foreach \x in {0, ..., 3}
{
\filldraw[blue!5!white] (-6, {\x+.5}) -- (.5, {\x+.5}) -- (.5, {\x-.5}) -- (-6, {\x-.5}) -- (-6, {\x+.5});
\draw[line width=0.6mm, blue] (-6, {\x+.5}) -- (.5, {\x+.5}) -- (.5, {\x-.5}) -- (-6, {\x-.5});
}

\foreach \x in {1, ..., 3}
{
\filldraw[blue!5!white] (-6, {-\x+.5}) -- (.5, {-\x+.5}) -- (.5, {-\x-.5}) -- (-6, {-\x-.5}) -- (-6, {-\x+.5});
\draw[line width=0.6mm, blue] (-6, {-\x+.5}) -- (.5, {-\x+.5}) -- (.5, {-\x-.5}) -- (-6, {-\x-.5});
}

\draw[line width=0.6mm, red!30!blue] (.5, -4) -- (.5, -3.5);

\draw[line width=0.6mm, red!30!blue] (.5, -2.5) -- (.5, -1.5);

\draw[line width=0.6mm, red!30!blue] (.5, .5) -- (.5, -.5);

\draw[line width=0.6mm, red!30!blue] (.5, 2.5) -- (.5, 1.5);

\draw[line width=0.6mm, red!30!blue] (.5, 4) -- (.5, 3.5);

\foreach \x in {0,...,3}
{
\filldraw[aqua] (0,{\x}) circle (0.7mm);
\draw[line width=0.1mm] (0,{\x}) circle (0.8mm);
\filldraw[aqua] (0,{-\x}) circle (0.7mm);
\draw[line width=0.1mm] (0,{-\x}) circle (0.8mm);
}
\foreach \x in {0,1}
{
\filldraw[aqua] (1,{2*\x}) circle (0.7mm);
\draw[line width=0.1mm] (1,{2*\x}) circle (0.8mm);
\filldraw[aqua] (1,{-2*\x}) circle (0.7mm);
\draw[line width=0.1mm] (1,{-2*\x}) circle (0.8mm);
}

\end{tikzpicture}
\caption{The Voronoi decomposition of Example~\ref{ex:quasi-tiling1}. The common parts of red and blue Voronoi cells are shown in purple. Note that the Voronoi cell of the point $(1,2n)$, $n\in \Z$, drawn in red, is not closed. The blue Voronoi cells are closed, and their vertices are purple. }
\label{fig:quasi-tiling}
\end{figure}

\begin{figure}[!t]
\centering
\begin{tikzpicture}
\draw[line width=.50mm, byzant] (0,2) -- (2,2);
\draw[line width=.50mm, byzant] (0,2) arc (150:30:1.15cm);
\draw[line width=.50mm, pomme] (0,2) arc (-150:-30:1.15cm);
\filldraw[aqua] (0.03,2) circle (0.7mm);
\draw[line width=0.1mm] (0.03,2) circle (0.8mm);
\filldraw[aqua] (1.97,2) circle (0.7mm);
\draw[line width=0.1mm] (1.97,2) circle (0.8mm);
\draw[-stealth, line width=.1mm, dashed] (.5, 3.4) arc (90:0:0.8cm);
\draw[-stealth, line width=.1mm, dashed] (.5, 3.4) arc (150:235:1cm);
\node at (.5, 3.6) {$\pi_1$};
\node at (.3, 1) {$\pi_2$};
\draw[-stealth, line width=.1mm, dashed] (.55, 1) arc (270:330:.8cm);
\end{tikzpicture}
\hspace{2cm}
\begin{tikzpicture}
\fill[blue!5!white] (0,0) -- (2,0) -- (2,2) -- (0,2) -- (0,0);
\draw[line width=0.40mm, lblue, dashed] (0.15,0) -- (.95,0);
\draw[line width=0.70mm, lblue] (1,0) -- (2.035,0);
\draw[line width=0.70mm, lblue] (2,0) -- (2,1);
\draw[line width=0.40mm, lblue, dashed] (2,1.1) -- (2,1.85);
\draw[line width=0.40mm, lblue, dashed] (1.85,2) -- (1.1,2);
\draw[line width=0.70mm, lblue] (1,2) -- (0,2);
\draw[line width=0.70mm, lblue] (0,2.035) -- (0,1);
\draw[line width=0.40mm, lblue, dashed] (0,0.9) -- (0,0.10);
\filldraw[lblue] (0,0) circle (0.4mm);
\filldraw[lblue] (2,2) circle (0.4mm);
\end{tikzpicture}
\caption{A tropical curve $\curve$ of rank two and its Voronoi cell. The vertices are all inside the Voronoi cell, the dashed lines are in the exterior. The closure of the Voronoi cell is isomorphic to the Voronoi cell of the first graded minor.  The second graded minor has trivial homology.}
\label{fig:voronoi}
\end{figure}
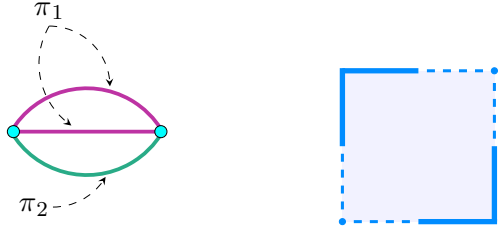

\begin{remark}
We stress that taking the closure of the Voronoi cell in Theorem~\ref{thm:VoronoiCellDecomposition-discrete} is necessary,  since the Voronoi cells are not always closed (see Figure~\ref{fig:quasi-tiling}, as well as Figure~\ref{fig:voronoi} and Section~\ref{sec:tropical_curves} for examples).
\end{remark}

\begin{cor}\label{cor:covering-discrete}
Let $S \subset H$ be an admissible discrete set. Then, the closures of Voronoi cells provide a tiling of $H$ with (generalized) polyhedra 
\[
H= \bigcup_{\gamma \in S} \overline\Vor_{S}(\gamma).
\]
Moreover, the polyhedra $\overline\Vor_{S}(\gamma)$, $\gamma \in S$, have mutually disjoint, non-empty interiors. 
\end{cor}
\begin{proof} By Theorem~\ref {thm:voronoi-discrete} and Theorem~\ref{thm:VoronoiCellDecomposition-discrete}, the closed Voronoi cells $\overline\Vor_{S}(\gamma)$ are (generalized) polyhedra and cover $H$. The remaining claim follows from Proposition~\ref{prop:general_voronoi}. Note that, due to convexity, $\overline\Vor_{S}(\gamma)$ and $\Vor_{S}(\gamma)$ have the same interior. 
\end{proof}

\begin{proof}[Proof of Theorem~\ref{thm:VoronoiCellDecomposition-discrete}] We first prove the claim in the case $0\in S$ and  $\gamma=0$, and then explain how this implies the general case. In this case, we have
\[S_j \coloneqq S_{\gamma,j}= \left\{\theta_j \st  \text{$\theta \in S$ and $\theta_1=\dots = \theta_{j-1}=0$} \right\} = \left\{\theta_j \st \theta \in S \cap \filter^j\right\}
\]
and
\[V_j \coloneqq V_{\gamma,j} = \Vor_{S_j}(0) \qquad j=1,\dots, r.\]
Since $V_{j} \subset H_j$, we have $V_i \aperp V_j$ for $i\neq j$. We will show that $\overline \Vor_{S}(0) =  V_1 \aplus \dots \aplus V_r$.

We first prove the inclusion $\overline \Vor_{S}(0) \subseteq V_1 \aplus \dots \aplus V_r$. It will suffice to show that
\[
\Vor_{S}(0) \subseteq V_1 \aplus \dots \aplus V_r,
\]
 as the right hand side is clearly closed in $H$. Let $x\in \Vor_{S}(0)$. Then we have the inequalities
\[
2 \highinn{}{x, \gamma} \lexeq \highinn{}{\gamma, \gamma} \qquad \forall\, \gamma \in S.
\]
Looking at the first coordinate in both sides, we obtain the inequalities
\[
2 \highinn{1}{x, \gamma} \leq \highinn{1}{\gamma,\gamma}\,\qquad \forall \gamma \in S.
\]
Using the almost orthogonal decomposition, these inequalities can be rewritten as 
\[2\highinn{1}{x_1, \gamma_1} \leq \highinn{1}{\gamma_1,\gamma_1}\,\qquad \forall \gamma \in S.\]
This shows that $x_1$ belongs to the Voronoi cell $V_1=\Vor_{S_1}(0)$ in $H_1$. 

The difference $x-x_1$ belongs to $\filter^2$, and for all $\gamma \in \filter^2 \cap S$, we have 
\[
2\highinn{2}{x-x_1,\gamma} =  2\highinn{2}{x,\gamma}  \lexeq \highinn{2}{\gamma,\gamma},
\]
which is equivalent to the inequalities 
\[
2\highinn{2}{x_2,\gamma_2}   \lexeq \highinn{2}{\gamma_2,\gamma_2} \qquad \forall \gamma \in \filter^2 \cap S.
\]
This means that $x_2$ belongs to the Voronoi cell $V_2 = \Vor_{S_2}(0)$ in $H_2$. Proceeding inductively in $j$, we find that $x_1 \in V_1, x_2\in V_2, \dots, x_j\in V_j$, and
\[
2\highinn{j+1}{x_{j+1},\gamma_{j+1}}   \lexeq \highinn{j+1}{\gamma_{j+1},\gamma_{j+1}} \qquad \forall \gamma \in \filter^{j+1} \cap S.
\] 
For $j=r$, we obtain $x \in V_1\aplus  \dots \aplus V_r$, as required.

 We now show the reverse inclusion $V_1 \aplus \dots \aplus V_r \subseteq \overline \Vor_{S}(0)$. Since the interior $\interior{V_j}$ is dense in $V_j$ (viewed as a subset of $H_j$), it will be enough to prove the inclusion 
 \[ \ring V_1 \aplus \dots  \aplus \ring V_r \subseteq  \Vor_{S}(0)  \]
 for the sets
 \[
 \ring V_j \coloneqq  \interior{V_j} = \interior{\Vor_{S_j}(0)}.
 \]
More precisely, here we take the interior of $V_j$ as a subset of $H_j$, and view it as a subset of $H$. Let $x =x_1+ \dots +x_r$ be a point in the left hand side, with $x_j \in \ring V_j$. We prove that $x \in \Vor_{S}(0)$ by verifying the inequalities
 \[
 2\highinn{}{x, \gamma} \lexeq \highinn{}{\gamma, \gamma} \qquad \forall \gamma \in S.
 \]
We will in fact prove the strict inequalities. 

Take an element $\gamma \in S \setminus \{0\}$, and let $j\in [r]$ be the integer with $\gamma \in \filter^j \smallsetminus \filter^{j+1}$.  Then,
\[
\highinn{1}{\gamma,\gamma}=\dots = \highinn{j-1}{\gamma,\gamma}=0, \qquad \highinn{j}{\gamma, \gamma}>0.
\]
The Lifting Lemma~\ref{lem:LiftingLemma} implies that
\[
\highinn{j}{x_i , \gamma} =0 \qquad \textrm{ for } i \neq j.
\]
For $i= j$, we have
\[
\highinn{j}{x_j , \gamma} = \highinn{j}{x_j , \gamma_j}.
\]
Since $x_j \in \ring V_j$, and $\gamma \in S\cap \filter^{j} \smallsetminus \filter^{j+1}$, it follows that $\gamma_j \neq 0$ and
\[
2\highinn{j}{x_j , \gamma_j} < \highinn{j}{\gamma_j , \gamma_j} = \highinn{j}{\gamma, \gamma}.
\]
Combining the above, we infer that, as required,
\[
2\highinn{}{x, \gamma} = (0, \dots, 0, 2\highinn{j}{x_j , \gamma} , \ast, \dots, \ast) \lexst (0, \dots, 0,  \highinn{j}{\gamma, \gamma}, \ast, \dots, \ast) = \highinn{}{\gamma, \gamma}.
\]

 The remaining claims on $\overline \Vor_{S}(0)$ follow from the decomposition just established and the properties of the (generalized) polyhedra $V_j$, $j=1, \dots, r$.

 \smallskip
 
The statement for $\gamma \in S$ follows by applying the decomposition to the translated set $S' = S-\gamma$, which is again admissible discrete, and noticing that 
\begin{align*}
\Vor_S(\gamma) &= \Vor_{S'}(0)+\gamma =(\Vor_{S'_1}(0) +\gamma_1) \aplus \dots \aplus (\Vor_{S'_r}(0)+\gamma_r) = V_{\gamma, 1}\aplus \dots \aplus V_{\gamma, r}. \qedhere
\end{align*}
  \end{proof}

\begin{example} \label{ex:quasi-tiling1} Consider $H = \R^2$ with the rank two Euclidean product $\highinn{}{\cdot\,,\cdot} \colon H\times H\to \R^2$ from Example~\ref{ex:HigherRankEuclidean}. Then, $S = \left \{(0, n) \st n \in \Z \right\}  \cup \left\{(1, 2n)\st n \in \Z\right \}$ is an admissible discrete set in $H$. The induced Voronoi decomposition is depicted in Figure~\ref{fig:quasi-tiling}. Taking the closures of the Voronoi cells, we obtain a polyhedral tiling of $\R^2$. Note that the intersection of two polyhedral tiles is not necessarily a common face. 
\end{example}

\section{Tame degenerations} \label{sec:tameness}
Throughout this section, let $H$ be a real vector space of finite dimension endowed with an inner product
\[
\highinn{}{x,y} = \bigl(\highinn{1}{x, y},\highinn{2}{x, y}, \dots, \highinn{r}{x, y} \bigr), \qquad x,y \in H,
\]
with values in $\Lambda = \R^r$ with the lexicographic order. Let 
\begin{equation} \label{eq:AlmostOrthoDecompositionState}
	H = H_1 \aplus H_2 \aplus  \dots \aplus H_r
\end{equation}
be the almost orthogonal decomposition with $H_j \coloneqq \proj_j^*(\grm{}{j}H)$ (see Lemma~\ref{lem:LiftingLemma}).
As before, each element $x\in H$ is written as $x=x_1+ \dots + x_r$ with $x_j \in H_j$, $j \in [r]$.

\begin{convention}[Reference norm] \label{con:ref-norm}
 In the following, we will need to endow $H$ with a {\em reference norm} denoted by $\refnorm{\cdot}$. By equivalence of norms on a given finite dimensional vector space, the discussion in the rest of the paper is independent of the choice of this reference norm. In order to streamline the presentation, we however make the natural choice given by 
\begin{equation} \label{eq:ReferenceNorm}
		\refnorm{x} = \Big ( \sum_{j=1}^r \highinn{j}{x_j\, ,x_j}  \Big)^{\frac{1}{2}}, \qquad x \in H. \qedhere
	\end{equation}
	\end{convention}
	The fact that the $j$th component $\highinn{j}{\cdot \,, \cdot}$ of the inner product $\highinn{}{\cdot \,, \cdot}$ is a scalar product on $H_j$ ensures that $\refnorm{\cdot}$ is indeed a norm on $H$.


\subsection{Properties of pullback families} \label{ss:PullbackProperties}
Consider a sequence $\underline L_t = (L_{t,1}, \dots, L_{t,r}) \in \R_+^r$, $t\in \R_+$, such that
\[
\lim_{t \to \infty} \frac{L_{t, j}}{L_{t,j+1}} = + \infty, \qquad 1 \le j \le r-1.
\] 
Consider the pullback family $\innone{t}{\cdot\,,\cdot}$, $t\in \R_+$,  given by
\begin{align}\label{eq:basic_tame_family}
	\innone{t}{x,y} = \innone{\underline L_t}{x, y} =  L_{t,1}\highinn{1}{x,y}+ \dots +L_{t,r}\highinn{r}{x,y}, \qquad x , y\in H,
\end{align}
Corollary~\ref{cor:OullbackFamiliesInnerProducts} ensures that $\innone{t}{\cdot\,,\cdot}$ is a scalar product on $H$ for $t$ large enough.

\begin{prop}\label{prop:tame_properties_toy} Let $\innone{t}{\cdot\,,\cdot}$, $t \in \R_+$, be a pullback family with parameters $\underline L_t = (L_{t,1}, \dots, L_{t,r})$. Then the following properties hold:
	\begin{itemize}
		\item There exists a constant $C_1>0$ such that for any distinct integers $i,j\in [r]$ and all large enough $t$, we have
		\begin{align} \label{eq:toy_tame1}
			\abs{\innone{t}{x,y}} \le C_1 L_{t,\max\{i,j\}+1}\, \refnorm{x} \cdot \refnorm{y} \qquad \forall x\in H_i, \, y \in H_j.
		\end{align}
		\item There is a constant $C_2>0$ such that for any $j\in [r]$ and all large enough $t$, we have 
		\begin{align} \label{eq:toy_tame2}
			\abs{ \innone{t}{x, y} - L_{t,j} \highinn{j}{x, y}} \le C_2 L_{t,j+1}  \refnorm{x} \cdot \refnorm{y} \qquad \forall x,y\in H_j.
		\end{align}
	\end{itemize}
	Here, $L_{t, r+1}\coloneqq0$ by convention.
\end{prop}

\begin{proof}
	For each $k \in [r]$, consider the bilinear form $(x,y) \mapsto \highinn{k}{x,y}$ on $H \times H$. Since $H$ has finite dimension, there exists a constant $c_k >0$ with
	\[
	\abs{\highinn{k}{x,y}} \le c_k \refnorm{x} \cdot \refnorm{y} \qquad \forall x,y \in H.
	\]
	By the properties of the subspaces $H_j$, $j \in [r]$, we have $\highinn{k}{x,y} = 0$ whenever $x \in H_i$ and $y \in H_j$ for indices $i,j$ with $k < \max\{i,j\}$ or $k = \max\{i,j\}$ and $i< j$. Thus,
	\[
	\innone{t}{x, y} = \sum_{k = \max \{i,j\}}^r L_{t,k} \highinn{k}{x,y}, \qquad x \in H_i, y \in H_j.
	\]
	Taking into account the preceding discussion, the above claims hold for
	\[
	C_1 \coloneqq \max_{i \neq j} \sum_{k = \max \{i,j\}+1}^r c_k, \qquad C_2 \coloneqq \max_j \sum_{k=j+1}^r c_k. \qedhere
	\]
\end{proof}

\subsection{Tameness: axiomatic of higher rank degenerations of scalar products}\label{ss:tameness}
In this section, we introduce \emph{tamely degenerating families of scalar products}. Tamely degenerating families are a slight generalization of pullback families that naturally appear when studying asymptotic geometry of complex algebraic varieties. Their properties simply allow to carry through the same proofs as for pullback families in the subsequent sections.  We will thus state our main results in this generalized setting. 

This generalization is necessary for geometric applications, where the appearing scalar products are usually not pullbacks (see for example our results in Section~\ref{sec:tropical_curves}). However, a reader who is mainly interested in polyhedral geometric results can skip this section and assume that all subsequently appearing families are pullback families. The polyhedral geometric results remain equally interesting for pullback families.

\smallskip

Let $\innone{t}{\cdot\,, \cdot} \colon H \times H \to \R$, $t\in \R_+$, be scalar products on $H$. Consider a family of vectors $\underline L_t = (L_{t,1}, \dots, L_{t,r}) \in \R_+^r$, $t\in \R_+$. 

\begin{defi}[Tame degeneration] \label{def:TameDegenerations} We say that a family $\innone{t}{\cdot\,, \cdot}$, $t\in \R_+$, of scalar products on $H$ {\em tamely degenerates} to the inner product $\highinn{}{\cdot\,, \cdot}$ with {\em parameters} $\underline L_t = (L_{t,1}, \dots, L_{t,r})\in \R_+^r$ if the following conditions hold:
	
	\begin{enumerate}[label=(\alph*)]
		\item \label{item:tame-a} For every $j =1, \dots, r-1$, we have
		\begin{equation} \label{eq:abstract_tame0}
			\lim_{t \to \infty} \frac{L_{t,j}}{L_{t,j+1}} = \infty.
		\end{equation}
		\item \label{item:tame-b}

		For any two distinct indices $i , j \in [r]$, and for every given $\varepsilon>0$, we have
		\begin{align} \label{eq:abstract_tame1}
			\abs{\innone{t}{x,y}} \le \varepsilon L_{t,\max\{i,j\}}\, \refnorm{x}\cdot\refnorm{y} \qquad \forall\,\, x \in H_i,  \, y \in H_j,
		\end{align}
		provided that $t\in \R_+$ is large enough.
		\item \label{item:tame-c} For every $j \in [r]$ and for every given $\delta >0$, we have
		\begin{align} \label{eq:abstract_tame2}
			\abs{L_{t,j}^{-1} \innone{t}{x, y} -  \highinn{j}{x, y}} \le \delta \refnorm{x}\cdot\refnorm{y}, \qquad \forall x,y \in H_j,
		\end{align}
		provided that $t\in \R_+$ is large enough. \qedhere
	\end{enumerate}
\end{defi}

In the geometric applications we are interested in, see for example Section~\ref{sec:tropical_curves}, it is sometimes necessary to replace \ref{item:tame-b} by the following weaker property:
\begin{enumerate}[label=(\alph*$'$)]
	\setcounter{enumi}{1}
	\item \label{item:weaktame-b}  There exists a constant $C>0$ such that for all $i, j \in [r]$ and all large $t$, we have
	\begin{equation} \label{eq:WeakTameEstimate}
		\abs{\innone{t}{x, y} } \le C L_{t, \max\{i,j\}} \refnorm{x} \refnorm{y}, \qquad \forall\,\, x \in \filter^i,\, y \in \filter^j.
	\end{equation}
\end{enumerate}
For this reason, we introduce the following variant of tame degenerations.
\begin{defi}[$\omega$-tame degeneration] Notations as above, 
	let $\innone{t}{\cdot, \cdot}$, $t\in \R_+$, be a family of scalar products on $H$ and $\underline L_t = (L_{t,1}, \dots, L_{t,r}) \in \R_+^r$, $t\in \R_+$, a family of vectors. 
	
	We say that the family $\innone{t}{\cdot, \cdot}$, $t\in \R_+$, {\em degenerates tamely in the weak sense}, or {\em degenerates $\omega$-tamely}, to $\highinn{}{\cdot\,, \cdot}$ with {\em parameters} $\underline L_t = (L_{t,1}, \dots, L_{t,r})$ if the properties \ref{item:tame-a}, \ref{item:weaktame-b}  and \ref{item:tame-c} are satisfied.
\end{defi}

The next proposition clarifies the relationship between pullback families, tamely degenerating families, and $\omega$-tamely degenerating families. 

\begin{prop}\label{prop:tameness_basic_family}
	Let $\highinn{}{\cdot \, , \cdot} \colon H \times H \to \R^r$ be an inner product on $H$. Then:
	\begin{itemize}
		\item [$(i)$] Let  $\innone{t}{\cdot\,,\cdot}$, $t \in \R_+$, be a pullback family of $\highinn{}{\cdot \, , \cdot}$ with parameters $\underline L_t,$ $t \in \R_+$. Then $\innone{t}{\cdot\,,\cdot}$, $t \in \R_+$, tamely degenerates to $\highinn{}{\cdot \, , \cdot}$ with parameters $\underline L_t,$ $t \in \R_+$.
		\item [(ii)] Let $\innone{t}{\cdot\,,\cdot}$, $t \in \R_+$, be a family of scalar products which tamely degenerates to $\highinn{}{\cdot \, , \cdot}$ with parameters $\underline L_t,$ $t \in \R_+$. Then, property \emph{\ref{item:weaktame-b}}  holds. In particular, $\innone{t}{\cdot\,,\cdot}$, $t \in \R_+$, $\omega$-tamely degenerates to $\highinn{}{\cdot \, , \cdot}$ with parameters $\underline L_t,$ $t \in \R_+$.
	\end{itemize}
\end{prop}
\begin{proof}
	In proving $(i)$, note that property \ref{item:tame-a} is part of the definition of pullback families. Properties \ref{item:tame-b} and \ref{item:tame-c} follow from \eqref{eq:toy_tame1} and \eqref{eq:toy_tame2}, respectively, using that $L_{t, j+1}/L_{t,j}$ tends to zero as $t$ goes to infinity.
	
	It remains to prove $(ii)$. Fix $i, j \in [r]$. We decompose $x \in \filter^i, y\in \filter^j$ as $x = x_i + \dots +x_r$, and $y \in \filter^j$ as $y = y_j + \dots +y_r$ with $x_k, y_k \in H_k$ for $k\in [r]$. Applying properties \ref{item:tame-b} and \ref{item:tame-c}  and taking into account \ref{item:tame-a}, we obtain
	\[
	\abs{\innone{t}{x, y} } \le \sum_{m = i}^r \sum_{n = j}^r \abs{\innone{t}{x_i, y_j} } \le \sum_{m = i}^r \sum_{n = j}^r  L_{\max\{i, j\}}\refnorm{x_i} \refnorm{y_j} \le L_{\max\{i, j\}} \big( \sum_{m = i}^r \refnorm{x_i} \big) \big( \sum_{n = j}^r \refnorm{y_j} \big) 
	\]
	for all $x\in \filter^i$ and $y\in \filter^j$. Since any two norms on $H$ are equivalent, there exists a constant $c>0$ such that $\sum_{k=1}^r \refnorm{z_k} \le c \refnorm{z}$ for all $z \in H$ decomposed as $z=z_1+\dots+z_r$. It follows that property \ref{item:weaktame-b}  holds for $C =c^2$.
\end{proof}
\begin{remark}
	One might wish to interpret ($\omega$-)tame degeneration as a notion of \emph{``convergence of scalar products to inner products of higher rank"}. However, this point of view should be taken with some precaution. For instance, \emph{``\,limits"} in this context are not unique, that is, a family of scalar products $\innone{t}{\cdot\,,\cdot}$, $t \in \R_+$, can ($\omega$-)tamely degenerate to two different inner products $\highinn{}{\cdot\,,\cdot} \neq \highinn{}{\cdot\,,\cdot}'$. For a simple example, see Proposition~\ref{prop:orthogonalization}\ref{prop:orthogonalization3}. Moreover, the parameters $\underline L_t \in \R_+^r$ for a ($\omega$)-tamely degenerating family of scalar products $\innone{t}{\cdot\,,\cdot}$, $t \in \R_+$, are not unique. In our geometric applications, the ($\omega$-)tame property of degenerations of scalar products allow to deduce interesting features in the limit. For a more detailed discussion, we refer to Section~\ref{ss:TameEquivalence}. 
\end{remark}

\subsection{Uniform equivalence to orthogonalization in tame degenerations}
Consider the almost orthogonal decomposition $\aplus_{j=1}^r H_j$ induced by the inner product $\highinn{}{\cdot\, ,\cdot}\colon H \times H \to \R^r$. As before, each element $x\in H$ is written as $x=x_1+ \dots + x_r$ with $x_j \in H_j$, $j = 1, \dots, r$.

The \emph{orthogonalization} of the inner product $\highinn{}{\cdot\, ,\cdot}$ is the new inner product
\begin{align*}
	\highinnoneorth{}{\cdot\,,\cdot} &\colon H \times H \to \R^r\\
	\highinnoneorth{}{x,y} &\coloneqq \bigl(\highinn{1}{x_1, y_1},\highinn{2}{x_2, y_2}, \dots, \highinn{r}{x_r, y_r} \bigr).
\end{align*} 
The following property holds by construction.

\begin{prop} \label{prop:orthogonalization} Let $\highinnoneorth{}{\cdot\,,\cdot} \colon H \times H \to \R^r$ be the orthogonalization of $\highinn{}{\cdot\, ,\cdot}$. Then:
	\begin{enumerate}[label=(\roman*)]
		\item The almost orthogonal decomposition induced by $\highinn{}{\cdot\,,\cdot}$ is orthogonal with respect to $\highinnoneorth{}{\cdot\,,\cdot}$. That is, if $x \in H_i$ and $y \in H_j$ with $i \neq j$, then $\highinnoneorth{}{x\,,y} = 0$.
		\item The reference norm $\refnorm{\cdot}$ coincides with the square root of the pullback $\innoneorth{\underline L}{\cdot\,,\cdot}$ of $\highinnoneorth{}{\cdot\,,\cdot}$ by the vector $\underline L = (1, 1, \dots, 1)$. More precisely,
		\[
		\refnorm{x} = \sqrt{\innoneorth{\underline L}{x\,,x}}, \qquad x \in H.
		\]
		\item \label{prop:orthogonalization3} The inner products $\highinn{}{\cdot\, ,\cdot}$ and $\highinnoneorth{}{\cdot\,,\cdot}$ have the same tamely degenerating families and the same $\omega$-tamely degenerating families.
	\end{enumerate}
\end{prop}
\begin{proof}
	Properties (i) and (ii) are clear from the definition of $\highinnoneorth{}{\cdot\,,\cdot}$ and $\refnorm{\cdot}$. Property (iii) follows from the fact that $\highinn{}{\cdot\,,\cdot}$ and $\highinnoneorth{}{\cdot\,,\cdot}$ induce the same almost orthogonal decomposition of $H$ and $\highinn{j}{x,y} = \highinnoneorth{j}{x,y}$ for all $x,y \in H_j$, see also Theorem~\ref{thm:TameEquivalence}. 
\end{proof}

Given a family of vectors $\underline L_t =\left(L_{1,t}, \dots, L_{r,t}\right) \in \R_+^r$, $t \in \R_+$, consider the pullback family $\innoneorth{t}{\cdot\,,\cdot} \coloneqq \innoneorth{\underline L_t}{\cdot\,,\cdot}$, $t \in \R_+$ given by
\[
\innoneorth{t}{x,y} = L_{t,1}\highinn{1}{x_1, y_1} + L_{t,2}\highinn{2}{x_2, y_2}+ \dots+L_{r,t}\highinn{r}{x_r, y_r}.
\]
Then, $\innoneorth{t}{\cdot\,, \cdot}$ is a scalar product on $H$ for all $t \in \R_+$. We denote the associated norm on $H$ by
\[
\normorth{x}_t \coloneqq \sqrt{\innoneorth{t}{x, x}} = \Big( \sum_{j=1}^r L_{t,j} \qf_j(x_j) \Big )^{1/2}, \qquad x \in H,
\]
where $\qf_j$ is the $j$-th component of the quadratic form $\qf$ associated to $\highinn{}{\cdot\,,\cdot}$, and we have
\[
\qf_j(z) = \highinn{j}{z, z}, \qquad z \in H_j.
\]
\begin{prop}[Uniform equivalence to orthogonalization pullbacks] \label{prop:apriori_estimate}
	Consider an inner product $\highinn{}{\cdot\,,  \cdot} \colon H \times H \to \R^r$.  Let $\innone{t}{\cdot \,, \cdot} \colon H \times H \to \R$, $t \in \R_+$, be a family of scalar products which $\omega$-tamely degenerates to $\highinn{}{\cdot\,,  \cdot}$ with parameters $\underline L_t \in \R_+^r$. Denote by $\norm{x}_t \coloneqq \sqrt{\innone{t}{x,x}}$, $t \in \R_+$, the associated norms on $H$.
	
	Then, the two families of norms $\norm{\cdot}_t $ and $\normorth{\cdot}_t$ are uniformly equivalent for large $t$. That is, there exists a constant $D >0$ such that
	\[
	D^{-1/2} \normorth{\cdot}_t \le \norm{\cdot}_t \le  D^{1/2} \normorth{\cdot}_t,
	\]
	provided that $t \in \R_+$ is large enough. Equivalently,
	\begin{equation} \label{eq:apriori_estimate}
		D^{-1} \sum_{j=1}^r L_{t,j} \normsq{x_j}_j \le \normsq{x}_t \le D \sum_{j=1}^r L_{t,j} \normsq{x_j}_j \qquad \forall x\in H,
	\end{equation}
	provided that $t \in \R_+$ is large enough.
\end{prop}

\begin{proof}
	Writing $x \in H$ as $x = x_1+ \dots + x_r$, with $x_j\in H_j$ for $j=1, \dots, r$, we have
	\[
	\normsq{x}_t = \sum_{j=1}^r \normsq{x_j}_t + 2 \sum_{i <j} \innone{t}{x_i, x_j}.
	\]
	It follows from \eqref{eq:abstract_tame2} that there exists a uniform constant $D' >0$ such that
	\[
	\frac 1{D'} \sum_{j=1}^r L_{t,j} \normsq{x_j}_j \ge \sum_{j=1}^r \normsq{x_j}_t  \ge D' \sum_{j=1}^r L_{t,j} \normsq{x_j}_j \qquad \forall \,\, x\in H \textrm{ and }  \forall \, \, t\in \R_+ \textrm{ large enough}.
	\]
	Fix now a pair of indices $i,j\in [r]$ with $i<j$. Combining the inequality
	\begin{equation} \label{eq:trick}
		2 ab \le\frac{1}{\delta} a^2 + \delta b^2, \qquad a,b \ge 0, \, \delta >0,
	\end{equation}
	with assumption \eqref{eq:WeakTameEstimate}, we conclude that for every $\delta >0$, there exists a constant $C_\delta >0$ with
	\[
	\abs{\innone{t}{x_i, x_j}} \le C L_{t,j} \norm{x_i}_i \cdot \norm{x_j}_j \le L_{t,j}\left(C_\delta \normsq{x_i}_i +  \delta \normsq{x_j}_j\right) \qquad \forall \,\, x\in H \textrm{ and } t \in \R_+ \textrm{ large enough}.
	\]
	Taking into account that $L_{t,i} /L_{t,j} \to \infty$ for $t \to \infty$, it follows that
	\[
	\abs{\innone{t}{x_i, x_j}} \le \delta L_{t,i} \normsq{x_i}_i +  \delta L_{t,j} \normsq{x_j}_j 
	\]
	for $t$ large. Choosing $\delta >0$ small and combining the above estimates, we conclude.
\end{proof}


\section{A finiteness lemma} \label{sec:FinitenessLemma}
In this section we state and prove the following finiteness lemma. It is a key ingredient in the arguments of the next sections. 

Consider an inner product $\highinn{}{\cdot\,,\cdot} \colon H\times H \to \R^r$. Let $\innone{t}{\cdot\, ,\cdot} \colon H\times H \to \R$, $t\in \R_+$, be scalar products which $\omega$-tamely degenerate to $\highinn{}{\cdot\,,\cdot}$ with parameters $\underline L_t \in \R_+^r$. Denote by $\norm{x}_t^2 = \innone{t}{x,x}$, $t \in \R_+$, the associated norms on $H$.

\begin{lem}[Finiteness Lemma] \label{lem:finiteness-discrete}
	Let $S \subset H$ be an admissible discrete set for $\highinn{}{\cdot\,,\cdot}$. Then, the following holds. For any compact subset $B \subseteq H$, there exists a finite subset $\S_B \subset S$ such that for any large enough $t$, we have
	\[
	\min_{\gamma \in S} \norm{x -\gamma}_t = \min_{\gamma \in \S_B} \norm{x -\gamma}_t \qquad \forall \,\, x \in B.
	\]
\end{lem}
For scalar products (i.e., when $r=1$), this lemma trivially follows from the compactness of $B$. On the contrary, the proof for higher values of $r$ is quite involved. 

As in the previous section, we use the almost orthogonal decomposition $H = \aplus_{j=1}^r H_j$ and decompose each $x \in H$ as $x = x_1 + \dots +x_r$ with $x_j \in H_j$, $j\in[r]$. Moreover by our Convention~\ref{con:ref-norm}, we endow $H$ with the reference norm 
\[
\refnorm{x} = \Big ( \sum_{j=1}^r \highinn{j}{x_j\, ,x_j}  \Big)^{\frac{1}{2}}, \qquad x \in H.
\]
\begin{proof} 
	The construction of the set $\S_B$ will be by induction. We will successively choose appropriate radii $R_1, \dots, R_r \in \R_+$ such that the sets 
	\[S^0 \coloneqq S\supseteq S^1 \supseteq S^2 \supseteq S^3 \supseteq \dots \supseteq S^{r}\]
	defined by
	\begin{equation} \label{eq:define-successive-sets}
		S^j = S^j(R_1, \dots, R_j) \coloneqq \Bigl\{ \gamma = \gamma_1 + \dots +\gamma_r \in S\,\, \st \, \, \refnorm{\gamma_k} \le R_k \text{ for all } k \le j \Bigr\}
	\end{equation}
	satisfy the property
	\begin{equation} \label{eq:finiteness-wish}
		\min_{\gamma \in S} \normsq{x -\gamma}_t = \min_{\gamma \in S^j} \normsq{x -\gamma}_t \qquad \forall \,\, x\in B,
	\end{equation}
	provided that $t$ is large enough. Moreover, we define $Z^j=Z^j(R_1, \dots, R_j)$ to be the set of all tuples $(\gamma_1, \dots, \gamma_j)$ appearing in the decomposition $\gamma =\gamma_1 + \dots +\gamma_r$ for $\gamma \in S^j$, and show that
	\begin{equation} \label{eq:finiteness-wish2}
		Z^j \textrm{ is a finite set.}
	\end{equation} 
	The last set $S^r$ coincides with $Z^r$, which is finite. The claim then follows by setting $\S_B \coloneqq S^r$.

	\smallskip
	
	The properties \eqref{eq:finiteness-wish} and \eqref{eq:finiteness-wish2} are trivially satisfied for $S^0 = S$ and $Z^0 =\emptyset$. 
	
	Suppose we have already constructed radii $R_1, \dots, R_j \in \R_+$ such that \eqref{eq:finiteness-wish} holds for $S^j$ defined by \eqref{eq:define-successive-sets}, and \eqref{eq:finiteness-wish2} is also satisfied for $Z^j$.

	For $\gamma \in S^j$ and $x \in H$, we then write
	\begin{align*}
		\normsq{x - \gamma}_t =  \normsq{x-(\gamma_1 + \dots +\gamma_j)}_t + \normsq{\gamma_{j+1} + \dots + \gamma_r}_t  - 2\sum_{k = j+1}^r \innone{t}{x- (\gamma_1 + \dots + \gamma_j), \gamma_k}
	\end{align*}
	for the decomposition $\gamma = \gamma_1 + \dots + \gamma_r$, $\gamma_j\in H_j$. Combing \eqref{eq:WeakTameEstimate} and the inequality \eqref{eq:trick}, we conclude that for every every $\delta >0$, there exists a constant $C_\delta >0$ such that
	\[
	\sum_{k = j+1}^r \abs{\innone{t}{x- (\gamma_1 + \dots + \gamma_j), \gamma_k}} \le \sum_{k =j+1}^r L_{t,k} \left(C_\delta \refnormsq{x- (\gamma_1 + \dots + \gamma_j)} + \delta \refnormsq{\gamma_k}\right)
	\]
	holds true for all $x\in B$ and $\gamma \in S^j$. Note that, by the definition of $S^j = S^j(R_1, \dots, R_j)$, and the compactness of $B$, the elements $x - (\gamma_1+ \dots+\gamma_j)$, for $\gamma \in S^j$ and $x \in B$, are all contained in a compact subset $\~ B$ of $H$. In particular, using \eqref{eq:abstract_tame0}, there exists a constant $D_\delta >0$ such that
	\[
	2 \sum_{k = j+1}^r \abs{\innone{t}{x- (\gamma_1 + \dots + \gamma_j), \gamma_k}} \le D_\delta L_{t,j+1} + \sum_{k =j+1}^r L_{t,k} \delta \refnormsq{\gamma_k}\qquad \forall \,\, \gamma \in S^j \textrm{ and } \forall \,\, x\in B.
	\]
	Applying Proposition~\ref{prop:apriori_estimate} and choosing $\delta >0$ small, we get the estimate
	\begin{align}\label{eq:estimate-induction1}
		\nonumber \normsq{x - \gamma}_t -  \normsq{x-(\gamma_1 + \dots +\gamma_j)}_t  &\ge L_{t,j+1} \left(D \refnormsq{\gamma_{j+1}} - \delta \refnormsq{\gamma_{j+1}} - D_\delta \right) + \sum_{k=j+2}^r L_{t,k} (D - \delta) \refnormsq{\gamma_k}  \\
		& \ge L_{t,j+1} \left(\frac{D}{2} \refnormsq{\gamma_{j+1}} - D_\delta\right),
	\end{align}
	holding uniformly for all $x \in B$,  $\gamma \in S^j$, and $t$ large.
	
	To proceed, we choose a large positive real $D' \in\R_+$,  to be determined below,  and choose $R_{j+1} \in \R_+$ large enough so that $\frac{D}{2} R_{j+1}^2 - D_\delta >3D'$. We then define the set $S^{j+1} = S^j(R_1, \dots, R_{j+1}) \subseteq S$ by \eqref{eq:define-successive-sets}.

	We claim that, for large $t$, and all $\gamma \in S^j \setminus S^{j+1}$ and $x \in B$, there exists $\tilde \gamma \in S^{j+1}$ such that we have the estimate
	\begin{equation} \label{eq:FinitenessFinalStep}
		\normsq{x - \gamma}_t > \normsq{x - \tilde \gamma}_t.
	\end{equation}
	This will prove the desired claim \eqref{eq:finiteness-wish}.

	Indeed for all $\gamma \in S^j \setminus S^{j+1}$, we have $\refnormsq{\gamma_{j+1}} \ge R_{j+1}^2$. The estimate~\eqref{eq:estimate-induction1} and the choice of $R_{j+1}$ gives
	\begin{align} \label{eq:estimate-induction2}
		\normsq{x - \gamma}_t & >  \normsq{x-(\gamma_1 + \dots +\gamma_j)}_t  + 3D' L_{t, j+1}
	\end{align}
	for large enough $t$.

	For all $z = (\gamma_1, \dots, \gamma_j) \in Z^j$, set 
	\[S_{/z} \coloneqq \Bigl\{ \theta =\theta_1+ \dots+\theta_r \in S\,\, \st\,\, \theta_1=\gamma_1, \dots, \theta_j = \gamma_j\Bigr\}.\]
	
	Since $B$ is compact and $Z^j$ is finite, the union $\bigcup_{z\in Z^j} \left(B -z\right)$ is compact. Using the properties \eqref{eq:abstract_tame0}-\eqref{eq:abstract_tame1}-\eqref{eq:abstract_tame2}-\eqref{eq:WeakTameEstimate}, we choose $D'$ large enough so that for all $z =(\gamma_1, \dots, \gamma_j)\in Z^j$, we can find an element $\tilde \gamma = \tilde\gamma_z \in S_{/z}$ such that the following inequalities hold...
	\begin{align}\label{eq:estimate-induction3}
		\normsq{\tilde \gamma_{j+1} + \dots  + \tilde \gamma_r}_t &\leq D' L_{t, j+1} \qquad \textrm{and} \\ 
		\label{eq:estimate-induction4} \abs{\innone{t}{x-\gamma_1-\dots-\gamma_j,  \tilde \gamma_{j+1} + \dots  + \tilde \gamma_r }} &\leq D' L_{t,j+1} \qquad  \textrm{for all $x\in B$},
	\end{align}
	provided that $t$ is large enough. Choosing $R_{j+1}$ large enough, we can ensure that $\tilde \gamma = \tilde\gamma_z$ belongs to $S^{j+1}$ for all $z \in Z^j$.
	
	To conclude, we combine \eqref{eq:estimate-induction2}, \eqref{eq:estimate-induction3}, and \eqref{eq:estimate-induction4} to get for all $\gamma \in S^j \setminus S^{j+1}$ with $\gamma \in S_{/z}$,
	\begin{align*}
		\normsq{x - \gamma}_t & >  \normsq{x-(\gamma_1 + \dots +\gamma_j)}_t  + 3D' L_{t, j+1} \\
		&\ge \normsq{x-(\gamma_1 + \dots +\gamma_j)}_t + \normsq{\tilde \gamma_{j+1} + \dots  + \tilde \gamma_r}_t + 2 D' L_{t,j+1} \\
		& =  \normsq{x- \tilde \gamma}_t - 2 \innone{t}{x-(\gamma_1 + \dots +\gamma_j),  \tilde \gamma_{j+1} + \dots  + \tilde \gamma_r } + 2 D' L_{t,j+1} \\
		&\ge \normsq{x- \tilde \gamma}_t .
	\end{align*}
	
	This establishes \eqref{eq:FinitenessFinalStep}.
	
	\smallskip
	
	It remains to show that the set $Z^{j+1} =Z^{j+1}(R_1, \dots, R_{j+1})$ is finite. This will be a direct consequence of the finiteness of $Z^j$ and the fact that $S \subset H$ is admissible discrete. Fix $z =(z_1, \dots, z_j)$ in $Z^j$ and consider the subset $S_{/z}  \subset S$. Since $S \subset H$ is admissible discrete, the subset $\bigl\{\gamma_{j+1} \,\st\, \gamma = \gamma_1 + \dots + \gamma_r \in S_{/z}\bigr\} \subset H_{j+1}$ is a discrete subset of $H_{j+1} \cong \grm{}{j+1}H$. Restricting to $S^{j+1} \cap S_{/z}$, we obtain that $\{\gamma_{j+1} \st \gamma = \gamma_1 + \dots + \gamma_r \in S_{/z} \cap S^{j+1}\}$ is finite. By finiteness of $Z^j$, we  conclude that $Z^{j+1}$ is finite.
	\end{proof}

\section{Hausdorff convergence of Voronoi cells}\label{sec:Hausdorff_convergence}
Throughout this section, let $\highinn{}{\cdot\, , \cdot} \colon H \times H \to \R^r$ be an inner product on $H$ and $S \subset H$ an admissible discrete set. In this section, we prove that the closed Voronoi cells $\overline \Vor_S(\gamma)$, $\gamma\in S$, appear as limits of Voronoi cells associated to any family of scalar products tamely degenerating to $\highinn{}{\cdot\, , \cdot}$.

\smallskip

We continue to denote by $\refnorm{\cdot} \colon H \to \R$ the reference norm on $H$ chosen according to our Convention~\ref{con:ref-norm}. The distance between a point $a \in H$ and a non-empty subset $X \subset H$ is denoted by $\dist(a, X) = \inf_{x \in X} \refnorm{a-x}$. Recall that the {\em Hausdorff distance} $\hdist(X, Y)$ between two non-empty subsets $X, Y \subset H$ is defined by 
\[
\hdist(X,Y) \coloneqq \max\Bigl\{ \sup_{x \in X} \dist(x,Y),  \, \sup_{y \in Y} \dist(y,X) \Bigr\}.
\]

A sequence $(X_t)_{t \in \R_+}$ of non-empty subsets $X_t \subset H$ {\em converges in the Hausdorff metric} to a non-empty subset $X \subset H$ if $\hdist(X_t,X) \to 0$ as $t \to \infty$. Convergence in the Hausdorff metric is equivalent to requiring that the distance functions $\dist(\cdot, X_t)$, $t\in \R_+$, converge to the distance function $\dist(\cdot, X)$ uniformly on $H$. That is, $\sup_{a \in H} | \dist(a, X_t) - \dist(a,X)| \to 0$ for $t \to\infty$.

\smallskip
Note that for unbounded sets $X,Y \subset H$ in general we have $\hdist(X,Y) = + \infty$.
We now deal with a variant of Hausdorff convergence for non-necessarily bounded subsets. We say that a sequence $(X_t)_{t \in \R_+}$ of non-empty subsets $X_t \subset H$ {\em converges compactly in the Hausdorff metric} to a non-empty subset $X \subset H$ if the distance functions $\dist(\cdot, X_t)$, $t\in \R_+$, converge to the distance function $\dist(\cdot, X)$ uniformly on compact subsets $P \subset H$. By Theorem~\ref{thm:CompactHausdorffConvergence}, proved later in Section~\ref{thm:CompactHausdorffConvergence}, this is equivalent to requiring that $(X_t\cap B_k)_{t\in \R_+}$ converges in the Hausdorff distance to $X \cap B_k$ for all sets $B_k$ in a covering $H = \bigcup_{k\in \N} B_k$ by open, bounded sets $B_k \subset H$.

(Note that since any two norms on $H$ are equivalent, replacing the reference norm $\refnorm{\cdot}$ by another norm on $H$ leads to the same notion of Hausdorff (compact) convergence.)

\smallskip

The next theorem is the main result of this section.

\begin{thm}[Hausdorff convergence of Voronoi cells]\label{thm:Hausdorff_voronoi} Let $\highinn{}{\cdot\, , \cdot} \colon H \times H \to \R^r$ be an inner product on $H$ and $S \subset H$ an admissible discrete set. Let $\innone{t}{\cdot\, , \cdot}$, $t\in \R_+$, be a family of scalar products which tamely degenerates to $\highinn{}{\cdot\,,\cdot}$ with parameters $\underline L_t \in \R_+^r$.  For $\gamma \in S$, let $W_t \coloneqq  \Vor_{S, \, \innone{t}{\cdot\,, \cdot}}(\gamma)$, $t \in \R_+$, be the Voronoi cell of $\gamma$ for $S$ and the scalar product $\innone{t}{\cdot\,, \cdot}$. 

Then, the Voronoi cells $(W_t)_{t \in \R_+}$ converge compactly in the Hausdorff metric to $\overline \Vor_{S}(\gamma)$.
In case that $\overline \Vor_S(\gamma)$ is compact, the Voronoi cells $(W_t)_{t \in \R_+}$ converge in the Hausdorff metric to $\overline \Vor_{S}(\gamma)$. 
\end{thm} 

As a special case, we get the following statement.

\begin{cor} Let $\innone{t}{\cdot\, , \cdot}$, $t\in \R_+$, be a pullback family for $\highinn{}{\cdot\,,\cdot}$ and $S\subset H$ an admissible discrete set. Then, we have...
\begin{enumerate}
\item  $\innone{t}{\cdot\, , \cdot} \colon H \times H \to \R$ is a scalar product for all large $t$, and
\item for $\gamma \in S$, the associated Voronoi cells $W_t = \Vor_{S, \, \innone{t}{\cdot\,, \cdot}}(\gamma)$, $t \in \R_+$, converge  compactly in the Hausdorff metric to $\overline \Vor_{S}(\gamma)$ as $t$ tends to infinity.  If $\overline \Vor_S(0)$ is compact, then the Voronoi cells $W_t$, $t \in \R_+$, converge to $\overline \Vor_{S}(\gamma)$ in the Hausdorff metric.
\end{enumerate} 
\end{cor}
\begin{proof} The first statement follows from Corollary~\ref{cor:OullbackFamiliesInnerProducts}. The second is a consequence of Theorem~\ref{thm:Hausdorff_voronoi} and Proposition~\ref{prop:tameness_basic_family}.
\end{proof}

The rest of this section is devoted to the proof of Theorem~\ref{thm:Hausdorff_voronoi}. First, we observe that, replacing $S$ by $S-\gamma$ if necessary, we can assume $0\in S$ and prove the result for $\gamma=0$.

We need some preparation. 

 \smallskip

Let $P$ be a polytope in $H$ containing the origin in its interior. Since such polytopes cover $H$, by Theorem~\ref{thm:CompactHausdorffConvergence}, it will be enough to prove the convergence of $W_t\cap P$ to $\Vor_S(0)\cap P$ in the Hausdorff metric.
 
 \smallskip
 
For $\varepsilon \in (0,1)$, define the set 
\[
U_\varepsilon \coloneqq (1-\varepsilon)\left(\overline \Vor_{S}(0)\cap P\right).
\]

\begin{lem} 
\label{lem:inclusion_voronoi} Let $\varepsilon \in (0,1)$. Then for all large $t \in \R_+$, we have the inclusion
\[
U_\varepsilon \subset W_t.
\] 
\end{lem}

\begin{proof} Since $\overline\Vor_{S}(0)$ is a (generalized) polyhedron and $P$ is a polytope, the intersection $\overline\Vor_{S}(0) \cap P$ is a polytope. Since $W_t$ is convex, it will be enough to show that for each vertex $x$ of $\overline\Vor_{S}(0) \cap P$, we have $(1-\varepsilon) x \in W_t$ for large enough $t$.
This amounts in proving the inequalities
\[
\innone{t}{ (1-\epsilon)x ,  (1-\epsilon) x } \leq \innone{t}{ (1-\epsilon)x -\gamma, (1-\epsilon)x -\gamma }\]
or, equivalently,  
\begin{equation} \label{eq:EpsilonInclusionInequality}
2(1-\epsilon) \innone{t}{x, \gamma} \leq \innone{t}{\gamma, \gamma}
\end{equation}
for all $\gamma \in S$ and $t$ large enough.  

Applying the Finiteness Lemma~\ref{lem:finiteness-discrete} to $B =\{(1-\epsilon)x\}$, consisting of a single point, we can find a finite subset $\S_B \subset S$ such that if the above inequalities hold for all $\gamma \in \S_B$, then they also hold for all $\eta \in S$.  Since $\S_B$ is finite, it thus suffices to prove the inequality for a single $\gamma$ and large $t$.

Let $i, j \in [r]$ such that $x \in \filter^i \smallsetminus \filter^{i+1}$ and $\gamma\in \filter^j\smallsetminus \filter^{j+1}$. Using the almost orthogonal decomposition $H =  H_1\aplus \dots \aplus H_r$,  we write $x = x_i + \dots + x_r$ and $\gamma = \gamma_j + \dots + \gamma_r$ with $x_k \in H_k$, $k=i, \dots, r$, and $\gamma_k\in H_k$, $k=j, \dots, r$.  Note that, by \ref{item:tame-a}, \ref{item:tame-c},  and \ref{item:weaktame-b}, the right hand side in \eqref{eq:EpsilonInclusionInequality} can be estimated by
\begin{equation} \label{eq:EpsilonInclusionAPriori}
\innone{t}{\gamma, \gamma} \ge (1-\varepsilon^2) L_{t,j} \highinn{j}{\gamma, \gamma} = (1-\varepsilon^2) L_{t,j} \highinn{j}{\gamma_j, \gamma_j}
\end{equation}
for all large $t$. (This holds for any constant $0<c<1$ instead of $(1-\varepsilon^2)$.)

\smallskip

The rest of the proof is by a case analysis depending on whether $i=j$, $i<j$, or $i>j$.

\smallskip

\noindent $\bullet$ Consider first the case $i = j$. Then, we decompose
\[
2(1-\varepsilon) \highinn{t}{x, \gamma} = 2(1-\varepsilon)\Big( \innone{t}{x_j, \gamma_j} + \innone{t}{x - x_j, \gamma_j}  +\innone{t}{x, \gamma - \gamma_j} \Big).
\]
By the Decomposition Theorem~\ref{thm:VoronoiCellDecomposition-discrete} for Voronoi cells, we have $x_j \in \Vor_{S_j}(0)$, for $S_j =\proj^*\left(\proj_j(\filter^j \cap S)\right) \subset H_j$.  Since $\gamma \in \filter^j\cap S$, we infer that $2\highinn{j}{x_j,\gamma_j}  \leq \highinn{j}{\gamma_j,\gamma_j}$.  Combining this with \ref{item:tame-a} and \ref{item:weaktame-b}, it follows that for every fixed $\delta >0$, we have
\begin{align}\label{eq:EpsilonInclusionAPriori2}
2(1-\varepsilon) \innone{t}{x, \gamma} \le (1-\varepsilon) L_{t,j} \highinn{j}{\gamma_j,\gamma_j} +\delta {L_{t,j}}
\end{align}
for large enough $t$. Choosing $\delta$ small, we obtain the desired estimate \eqref{eq:EpsilonInclusionInequality} from \eqref{eq:EpsilonInclusionAPriori}.

\smallskip

\noindent $\bullet$ Consider now the case $i<j$. Then, we can again decompose
\[
2(1-\varepsilon) \innone{t}{x, \gamma} = 2(1-\varepsilon)\Big( \innone{t}{\sum_{k < j} x_k, \gamma} + \innone{t}{\sum_{k \ge j} x_k, \gamma} \Big) 
\]
Applying the steps from the previous case $i =j$ to $\tilde x =\sum_{k\ge j} x_k$, we obtain the same estimate as \eqref{eq:EpsilonInclusionAPriori2}:
\[
2  (1-\varepsilon) \innone{t}{\sum_{k \ge j} x_k, \gamma} \le (1-\varepsilon) L_{t,j} \highinn{j}{\gamma_j,\gamma_j} +\delta {L_{t,j}}
\]
for every fixed $\delta >0$ and $t$ large enough. On the other hand, \ref{item:tame-a} and
\ref{item:tame-b} imply that
\[
2(1-\varepsilon)  \innone{t}{\sum_{k < j} x_k, \gamma}  \le \delta L_{j,t}
\]
for every fixed $\delta >0$ and $t$ large enough. Choosing again $\delta$ small enough, the desired estimate \eqref{eq:EpsilonInclusionInequality} follows again from \eqref{eq:EpsilonInclusionAPriori}.

\smallskip

\noindent $\bullet$ Finally, consider the case $i>j$. By \ref{item:weaktame-b} and \ref{item:tame-a}, we then have
\begin{align*}
 2 (1-\varepsilon) \innone{t}{x, \gamma} \leq 2 (1-\varepsilon)C   L_{t,i} \refnorm{x} \refnorm{\gamma} \leq  (1-\varepsilon^2) L_{t,j} \highinn{j}{\gamma_j, \gamma_j}
 \end{align*}
for all large $t$. Taking into account \eqref{eq:EpsilonInclusionAPriori}, the estimate \eqref{eq:EpsilonInclusionInequality} follows.
\end{proof}

\begin{lem} \label{lem:DegeneratingCellsConstained} Suppose that $\varepsilon \in (0,1)$. Then, all accumulation points of the sequence of sets $(1-\varepsilon) W_t$, $t \in \R_+$, are in $\Vor_{S}(0)$. Moreover, if $P \subset H$ is a polytope with $0 \in \interior{P}$, then 
\[
(1-\varepsilon)\left(W_t\cap P\right) \subset \Vor_{S}(0)
\]
for all large $t$. Finally, if $\overline \Vor_S(0)$ is compact, then
\[
(1-\varepsilon) W_t \subset \Vor_{S}(0)
\]
for all large $t$.
\end{lem}
\begin{proof} Let $x_t \in (1-\varepsilon)W_t$ be a sequence of points, which converges to a point $x \in H$ as $t$ tends to infinity. We need to prove the inequalities 
\[2\highinn{}{x,\gamma} \lexeq \highinn{}{\gamma,\gamma} \qquad \forall \gamma \in S.\]
Since $x_t \in (1-\varepsilon)W_t$, we have the inequalities
\[
2 \innone{t}{x_t,\gamma} \leq (1-\varepsilon)\innone{t}{\gamma, \gamma} \qquad \forall \gamma \in S
\]
for all $t$. Fix $\gamma \in S$ and let $j\in [r]$ be so that $\gamma \in \filter^j \smallsetminus \filter^{j+1}$. Decompose $x_t = x_{t,1}+ \dots+x_{t,r}$ and $x = x_1+\dots+x_r$ according to the almost orthogonal decomposition $H =\aplus_j H_j$. Then, $x_{t,j}$ converges to $x_j \in H_j$ for all $j$.

Observe first that, since $\gamma \in \filter^j$, we have $\highinn{i}{x,\gamma} =\highinn{i}{\gamma,\gamma}=0$ for all $i<j$. We now prove that
 $2\highinn{j}{x,\gamma} \leq (1- \varepsilon^2)\highinn{j}{\gamma, \gamma}$. This implies that $2\highinn{}{x,\gamma} \lexeq \highinn{}{\gamma,\gamma}$, as required.  

We first decompose
\[
\innone{t}{x_t,\gamma} = \innone{t}{x_{t,j},\gamma_j} + \innone{t}{x_{t,j},\gamma - \gamma_j} +  \sum_{i>j} \innone{t}{x_{t,i},\gamma} + \sum_{i<j} \innone{t}{x_{t,i},\gamma_j}  + \sum_{i<j} \innone{t}{x_{t,i},\gamma - \gamma_j}. 
\]
Fix $\delta >0$ and let $C>0$ be the constant in \ref{item:weaktame-b}. Applying \ref{item:tame-b}, \ref{item:weaktame-b}, and \ref{item:tame-c}, we conclude that
\begin{align*}
\innone{t}{x_t,\gamma} &\ge L_{t,j} \innone{t}{x_{t,j},\gamma_j} - \delta L_{t,j} \refnorm{x_{t,j}} \refnorm{\gamma_j} - C L_{t,  j+1} \refnorm{x_{t,j}} \refnorm{\gamma-\gamma_j} \\
&\quad - C L_{t,j+1} \sum_{i >j} \refnorm{ x_{t,i}} \refnorm{\gamma} - \delta L_{t,j} \sum_{i<j} \refnorm{x_{t,i}} \refnorm{\gamma_j} - CL_{t, j+1}\sum_{i <j} \refnorm{ x_{t,i}} \refnorm{\gamma -\gamma_j}
\end{align*}
for $t$ large. Choosing $\delta>0$ small and applying \ref{item:tame-a}, it follows that
\[
2 L_{t,j} \highinn{j}{x_{t,j},\gamma} = 2 L_{t,j} \highinn{j}{x_{t,j},\gamma_j} \leq  (1- \varepsilon^2)\innone{t}{\gamma, \gamma}
\]
for large $t$. Here, we have used that $\innone{t}{\gamma, \gamma} \ge c L_{t,j}$ for some $c>0$, which follows from \ref{item:tame-a}, \ref{item:weaktame-b}, and \ref{item:tame-c}. We then get the desired inequality by dividing by $L_{t,j}$ and taking the limit for $t$ tending to $\infty$,
\[
2\highinn{j}{x,\gamma}   = \lim_{t \to \infty} 2 \highinn{j}{x_{t,j},\gamma} \le \lim_{t \to \infty} (1-\varepsilon^2)L_{t,j}^{-1}\innone{t}{\gamma, \gamma} = (1-\varepsilon^2) \highinn{j}{\gamma, \gamma}. 
\]
(Here, to compute the limit on the right hand side, we have used again \ref{item:tame-a}, \ref{item:weaktame-b}, and \ref{item:tame-c}.)

\smallskip

We now prove that $(1-\varepsilon) \left(W_t\cap P \right) \subset \Vor_{S}(0)$ for all large $t$. For the sake of a contradiction, suppose that this is not the case. Then, we can find a subsequence $x_t \in W_t \cap P$ such that $(1-\varepsilon) x_t \notin \Vor_{S}(0)$ for all $t$. 
 Up to choosing a subsequence, we can assume that $x_t$ converges to a point $x \in H$. In particular, the sequences $(1-\varepsilon^2) x_t$ and $(1-\varepsilon) x_t$ converge to $y = (1-\varepsilon^2) x$ and $z =(1- \varepsilon)x$,  respectively. By the above result, we have $y \in \Vor_S(0)$. Since $\Vor_S(0)$ is convex and contains $0$ in its interior, it follows that $z = \frac1{1+\varepsilon} \cdot y$ belong to the interior of $\Vor_S(0)$. However, then by the convergence, the points $(1-\varepsilon) x_t$ lie eventually in $\Vor_S(0)$. This contradiction completes the proof.

Finally, assume that $\overline{\Vor}_S(0)$ is compact. Let $B \subset H$ be an open, relatively compact set with $2\cdot \overline{\Vor}_S(0) \subset B$. Then $W_t \subset B$ for all $t$ large, which proves the remaining claim upon choosing a large polytope $P$ with $B \subset P$. Namely, suppose there exists an infinite sequence of points $x_t$ with  $x_t \in W_t \setminus B$. Using that $W_t$ is convex and contains $0 \in H$ and passing to a subsequence, we can assume that $x_t$ remains bounded and converges to some point $x \in  H$. By the above result, we have that $x = 2 \cdot\frac{1}{2} x$ lies in $2 \overline{\Vor}_S(0) \subset B$. Thus, $x_t$ belongs to $B$ for all large $t$, a contradiction.
\end{proof}

\begin{proof}[Proof of Theorem~\ref{thm:Hausdorff_voronoi}] By the preceding lemmas, the inclusions $(1-\varepsilon)\left(\Vor_{S}(0)\cap P\right) \subset W_t$ and $(1-\varepsilon)\left(W_t\cap P\right) \subset \overline \Vor_{S}(0)$ hold for $t$ large enough. In particular, for a given $\varepsilon >0$, we have
\[
\sup_{x \in W_t\cap P} \inf_{y \in \overline \Vor_{S}(0)\cap P} \refnorm{x - y} \le  \epsilon \sup_{x \in W_t \cap P} \refnorm{x} + \sup_{x \in W_t \cap P} \inf_{y \in \overline \Vor_{S}(0)\cap P} \refnorm{(1-\varepsilon) x  - y} = \epsilon \sup_{x \in W_t\cap P} \refnorm{x}
\]
for all $t$ large. Arguing similar for $\sup_{y \in \overline \Vor_{S}(0)\cap P} \inf_{x \in W_t\cap P} \refnorm{x - y}$, we obtain the claim.

 If $\overline\Vor_S(0)$ is compact, then from Lemma~\ref{lem:DegeneratingCellsConstained}, we infer that $W_t$ remains bounded as $t$ goes to infinity. Choosing a polytope $P$ which contains $\overline \Vor_S(0)$ and all the Voronoi cells $W_t$ for $t$ large, we conclude. 
\end{proof}


\section{Admissible lattices} \label{sec:AdmissibleLattices}
Let $H$ be a real vector space endowed with an inner product $\highinn{}{\cdot\, ,\cdot}\colon H\times H \to \R^r$. Consider the filtration
\[
\filter^\bullet \colon \quad \filter^1=H \supseteq \filter^2 \supseteq \dots \supseteq \filter^r \supset \filter^{r+1} =(0)
\]
on $H$ induced by $\highinn{}{\cdot\, ,\cdot}$.

Let $\L$ be a lattice in  $H$. The filtration $\filter^\bullet$ induces a non-increasing filtration 
\[
\L^\bullet \colon \quad \L^1=\L \supseteq \L^2 \supseteq \dots \supseteq \L^r \supseteq \L^{r+1} =(0)
\]
on $\L$ by setting $\L^j \coloneqq \filter^j\cap \L $ for all $j \in [r]$. We call a lattice $\L$ in $H$ \emph{admissible} if it is an admissible discrete subset of $H$. 

The following theorem gives a characterization of admissible lattices. For a lattice $\L$ in $H$, we define the $j$th {\em graded piece} of $\L$ as the quotient
\[
\grm{}{j}\L\coloneqq \rquot{\L^j}{\L^{j+1}}, \qquad j \in [r].
\]
Note that the inclusion $\L^j \subset \filter^j$ defines a natural embedding $\grm{}{j}\L  \hookrightarrow \grm{}{j}H$, so that we can identify $\grm{}{j}\L$ with the corresponding subset of $\grm{}{j}H$. 
\begin{thm} \label{thm:AdmissibleDiscreteLattices} The following are equivalent for a lattice $\L$ in $H$:
\begin{enumerate}
\item $\L$ is admissible.
\item Each graded piece $\grm{}{j}\L = \rquot{\L^j}{\L^{j+1}}$, $j \in [r]$, is a lattice in $\grm{}{j}H = \rquot{\filter^j}{\filter^{j+1}}$. 
\end{enumerate}
Moreover, if $\L$ is of full rank, then $\L$ is admissible if and only if each $\L^j$, $j \in [r]$, is a lattice of full rank in $\filter^j$.
\end{thm}

\begin{proof} 
We prove the equivalence of (1) and (2). It will be more convenient to use the characterization of admissible discreteness given in Remark~\ref{rem:AdmissibleDiscreteFiltration} in terms of the filtration. Let $\L$ be a lattice in $H$. We proceed by induction on the dimension $r$ of $\R^r$. For $r=1$, the claim is trivial. Suppose the claim holds for $r-1$ and consider an inner product $\highinn{}{\cdot\, ,\cdot}\colon H\times H \to \R^r$. The last $r-1$ components of $\highinn{}{\cdot\, ,\cdot}$ define an inner product $[\cdot \, , \cdot] =(\highinn{2}{\cdot\, ,\cdot}, \dots, \highinn{r}{\cdot\, ,\cdot})$ on $\filter^{2}$.

 We first prove that (2) $\Rightarrow$ (1). Assume that $\rquot{\L^j}{\L^{j+1}}$ is a lattice in $\rquot{\filter^j}{\filter^{j+1}}$, for $j\in[r]$. The projection of $\L$ into $\grm{}{1}H$ is the graded piece $\grm{}{1}\L$, which is a lattice. Therefore, it is discrete. By Remark~\ref{rem:AdmissibleDiscreteFiltration}, using the notation of \eqref{eq:SOver-filter}, it remains to show that $\L_{/\proj_1(\gamma)}-\gamma$ is admissible discrete in $\filter^2$ for all $\gamma\in \L$. However, we have that $\L_{/\proj_1(\gamma)}-\gamma = \L^2$. Applying the induction hypothesis to $[\cdot \, , \cdot] $ and $\L^2$, we conclude.

We now prove the reverse implication (1) $\Rightarrow$ (2). Suppose $\L$ is admissible. Proceeding by induction, we prove that $\grm{}{j}\L$ is a lattice in $\grm{}{j}H$. 

Since $\L$ is admissible, $\grm{}{1}\L$ is a discrete subgroup in $\grm{}{1}H$, and hence a lattice. By Remark~\ref{rem:AdmissibleDiscreteFiltration}, using~\eqref{eq:SOver-filter}, admissibility implies that $\L^2 =\L_{/0_1}$ is admissible in $\filter^2$, where $0_1 \in \grm{}{1}H$ is the origin in $\grm{}{1}H$. Applying once again the induction hypothesis to $[\cdot \, , \cdot] $ and $\L^2$, the proof is complete.

Having established this, from the equalities
\[
\sum_{j=1}^r \rk \left(\grm{}{j}\L\right)= \rk(\L), \qquad  \sum_{j=1}^r \dim\left(\grm{}{j}H\right) = \dim(H),\]
we infer that, if $\L$ has full rank, then $\L$ is admissible discrete exactly when $\grm{}{j}\L$ is a lattice of full rank in $\grm{}{j}H$ for all $j$.
\end{proof}

Recall that the $j$th graded piece $\grm{}{j}H = \rquot{\filter^j}{\filter^{j+1}}$ is endowed with the scalar product $\highinn{j}{\cdot\, ,\cdot} \colon \grm{}{j}H \times \grm{}{j}H \to \R$. For an admissible lattice $\L$, we denote by $\Vor_{\grm{}{j}\L}(\gamma)$, $\gamma \in \grm{}{j}\L$, the Voronoi cells of the lattice $\grm{}{j}\L$ in $\grm{}{j}H$ and consider the corresponding Voronoi decomposition
\[
 \grm{}{j}H = \bigcup_{\gamma \in \grm{}{j}\L} \Vor_{\grm{}{j}\L}(\gamma).
\]

The following theorem allows to decompose the Voronoi cells of admissible lattices in terms of the Voronoi cells of their graded pieces. 

\begin{thm}[Almost orthogonal decomposition for Voronoi cells of admissible lattices]\label{thm:VoronoiCellDecomposition-lattices} Let $\highinn{}{\cdot\, ,\cdot}\colon H\times H \to \R^r$ be an inner product and $\L \subset H$ be an admissible lattice. For each $j \in [r]$, let 
\[
V_j \coloneqq \proj_j^*\left(\Vor_{\grm{}{j}\L}(0)\right) \subset H_j \subseteq H.
\]
Then, the closure of $\Vor_{\L}(0)$ has the almost orthogonal decomposition
\[
\overline \Vor_{\L}(0) = V_1 \aplus \dots \aplus V_r.
\]
In words, the closure of $\Vor_{\L}(0)$ coincides with the Minkowski sum of the canonical liftings of the Voronoi cells $\Vor_{\grm{}{j}\L}(0) \subset \grm{}{j}H$ to $H_j \subseteq H$. Moreover, the sets $V_j$, $j \in [r]$, in the decomposition are pairwise almost orthogonal.

For any $\gamma \in \L$, we have $\Vor_{\L}(\gamma) = \gamma+ \Vor_{\L}(0)$. The closure of each Voronoi cell $\Vor_{\L}(\gamma)$, $\gamma \in \L$, is a polyhedron. If $\L$ is of full rank, then $\overline \Vor_{\L}(\gamma)$, $\gamma \in \L$, is a polytope with a non-empty interior, and $\Vor_{\L}(\gamma)$ is a bounded convex set with a non-empty interior.
\end{thm}
\begin{proof} The theorem follows directly from Theorem~\ref{thm:VoronoiCellDecomposition-discrete}. Note that the Voronoi cell $V_j$ appearing in the decomposition of $\overline \Vor_{\L}(0)$ is congruent to the Voronoi cell of the origin in $\grm{}{j}H$ with respect to the lattice $\grm{}{j}\L$. This is always a polyhedron.
\end{proof}
We restate Corollary~\ref{cor:covering-discrete} in the case of admissible lattices.
\begin{cor}\label{cor:covering-lattices}
Let $\L$ be an admissible lattice in an inner product space $(H, \highinn{}{\cdot\, ,\cdot})$. Then, the Voronoi cells cover the space $H$, that is,
\[
H= \bigcup_{\gamma \in \L} \Vor_{\L}(\gamma),
\]
and they have mutually disjoint, non-empty, interiors. Moreover, the closures $\overline \Vor_{\L}(\gamma)$, $\gamma \in \L$, provide a covering of $H$ by polyhedra. If $\L$ is of full rank, the closure of each Voronoi cell is a polytope. 
\end{cor}

The following example shows that in general, the intersection of two Voronoi cells is not necessarily a common face of both of them. Moreover, the Voronoi cells are not always closed. Even more surprisingly, they may not even contain any of their vertices.

\begin{figure}[!t]
\centering
\begin{tikzpicture}[scale=.7]
\fill[blue!5!white] (0,-1.03) -- (2,-0.03) -- (2,2) -- (0,1) -- (0,0);
\fill[red!5!white] (0,0) -- (2,1) -- (2,3) -- (0,2) -- (0,0);
\fill[blue!5!white] (0,2) -- (2,3) -- (2,5) -- (0,4) -- (0,2);
\fill[red!5!white] (0,4) -- (2,5) -- (2,7) -- (0,6) -- (0,4);

\fill[brown!10!white] (2,0) -- (4,1) -- (4,3) -- (2,2) -- (2,0);
\fill[green!15!white] (2,2) -- (4,3) -- (4,5) -- (2,4) -- (2,2);
\fill[brown!10!white] (2,4) -- (4,5) -- (4,7) -- (2,6) -- (2,4);
\fill[green!15!white] (2,6) -- (4,7) -- (4,8.03) -- (2,7.03) -- (2,6);

\fill[blue!5!white] (4,.97) -- (6,1.97) -- (6,3) -- (4,2) -- (4,.97);
\fill[red!5!white] (4,2) -- (6,3) -- (6,5) -- (4,4) -- (4,2);
\fill[blue!5!white] (4,4) -- (6,5) -- (6,7) -- (4,6) -- (4,4);
\fill[red!5!white] (4,6) -- (6,7) -- (6,9) -- (4,8) -- (4,6);

\draw[line width=0.50mm, byzant] (0,0) -- (2,1);
\draw[line width=0.50mm, byzant] (0,2) -- (2,3);
\draw[line width=0.50mm, byzant] (0,4) -- (2,5);
\draw[line width=0.50mm, byzant] (0,6) -- (2,7);

\draw[line width=0.50mm, ggray] (2,0) -- (4,1);
\draw[line width=0.50mm, ggray] (2,2) -- (4,3);
\draw[line width=0.50mm, ggray] (2,4) -- (4,5);
\draw[line width=0.50mm, ggray] (2,6) -- (4,7);

\draw[line width=0.50mm, byzant] (4,2) -- (6,3);
\draw[line width=0.50mm, byzant] (4,4) -- (6,5);
\draw[line width=0.50mm, byzant] (4,6) -- (6,7);
\draw[line width=0.50mm, byzant] (4,8) -- (6,9);

\draw[line width=0.55mm, lblue] (0,-1.01) -- (0,-.5);
\draw[line width=0.55mm, dred] (0,.5) -- (0,1.5);
\draw[line width=0.55mm, lblue] (0,2.5) -- (0,3.5);
\draw[line width=0.55mm, dred] (0,4.5) -- (0,5.5);

\draw[line width=0.55mm, lblue] (2,-0.02) -- (2,.5);
\draw[line width=0.55mm, dred] (2,1.5) -- (2,2.5);
\draw[line width=0.55mm, lblue] (2,3.5) -- (2,4.5);
\draw[line width=0.55mm, dred] (2,5.5) -- (2,6.5);

\draw[line width=0.55mm, lblue] (4,1) -- (4,1.5);
\draw[line width=0.55mm, dred] (4,2.5) -- (4,3.5);
\draw[line width=0.55mm, lblue] (4,4.5) -- (4,5.5);
\draw[line width=0.55mm, dred] (4,6.5) -- (4,7.5);

\draw[line width=0.55mm, lblue] (6,2) -- (6,2.5);
\draw[line width=0.55mm, dred] (6,3.5) -- (6,4.5);
\draw[line width=0.55mm, lblue] (6,5.5) -- (6,6.5);
\draw[line width=0.55mm, dred] (6,7.5) -- (6,8.5);

\draw[line width=0.55mm, brown] (0,-.5) -- (0,0.5);
\draw[line width=0.55mm, pgreen] (0,1.5) -- (0,2.5);
\draw[line width=0.55mm, brown] (0,3.5) -- (0,4.5);
\draw[line width=0.55mm, pgreen] (0,5.5) -- (0,6);

\draw[line width=0.55mm, brown] (2,.5) -- (2,1.5);
\draw[line width=0.55mm, pgreen] (2,2.5) -- (2,3.5);
\draw[line width=0.55mm, brown] (2,4.5) -- (2,5.5);
\draw[line width=0.55mm, pgreen] (2,6.5) -- (2,7);

\draw[line width=0.55mm, brown] (4,1.5) -- (4,2.5);
\draw[line width=0.55mm, pgreen] (4,3.5) -- (4,4.5);
\draw[line width=0.55mm, brown] (4,5.5) -- (4,6.5);
\draw[line width=0.55mm, pgreen] (4,7.5) -- (4,8);

\draw[line width=0.55mm, brown] (6,2.5) -- (6,3.5);
\draw[line width=0.55mm, pgreen] (6,4.5) -- (6,5.5);
\draw[line width=0.55mm, brown] (6,6.5) -- (6,7.5);
\draw[line width=0.55mm, pgreen] (6,8.5) -- (6,9);

\filldraw[lblue] (1,-.5) circle (0.7mm);
\draw[line width=0.1mm] (1,-.5) circle (0.8mm);
\filldraw[dred] (1,1.5) circle (0.7mm);
\draw[line width=0.1mm] (1,1.5) circle (0.8mm);
\filldraw[lblue] (1,3.5) circle (0.7mm);
\draw[line width=0.1mm] (1,3.5) circle (0.8mm);
\filldraw[dred] (1,5.5) circle (0.7mm);
\draw[line width=0.1mm] (1,5.5) circle (0.8mm);

\filldraw[brown] (3,1.5) circle (0.7mm);
\draw[line width=0.1mm] (3,1.5) circle (0.8mm);
\filldraw[pgreen] (3,3.5) circle (0.7mm);
\draw[line width=0.1mm] (3,3.5) circle (0.8mm);
\filldraw[brown] (3,5.5) circle (0.7mm);
\draw[line width=0.1mm] (3,5.5) circle (0.8mm);
\filldraw[pgreen] (3,7.5) circle (0.7mm);
\draw[line width=0.1mm] (3,7.5) circle (0.8mm);

\filldraw[lblue] (5,1.5) circle (0.7mm);
\draw[line width=0.1mm] (5,1.5) circle (0.8mm);
\filldraw[dred] (5,3.5) circle (0.7mm);
\draw[line width=0.1mm] (5,3.5) circle (0.8mm);
\filldraw[lblue] (5,5.5) circle (0.7mm);
\draw[line width=0.1mm] (5,5.5) circle (0.8mm);
\filldraw[dred] (5,7.5) circle (0.7mm);
\draw[line width=0.1mm] (5,7.5) circle (0.8mm);

\foreach \x in {1,...,4}
{
\foreach \y in {1,...,7}
{
\filldraw[red!30!blue] ({2*\x-2.03},{\x+\y+.5-3.04}) -- ({2*\x-1.97},{\x+\y+.5-3.01}) -- ({2*\x-1.97},{\x+\y+.5-2.95}) -- ({2*\x-2.03},{\x+\y+.5-2.98}) -- ({2*\x-2.03},{\x+\y+.5-3.04}) -- ({2*\x-1.97},{\x+\y+.5-3.01});
}
}

\foreach \x in {1,...,4}
{
\filldraw[lblue] ({2*\x-2.029},{\x-2.04}) -- ({2*\x-1.97},{\x-2.01}) -- ({2*\x-1.97},{\x-1.95}) -- ({2*\x-2.029},{\x-1.98}) -- ({2*\x-2.029},{\x-2.04}) -- ({2*\x-1.97},{\x-2.01});
}

\foreach \x in {1,...,4}
{
\filldraw[pgreen] ({2*\x-2.03},{\x+8-3.04}) -- ({2*\x-1.97},{\x+8-3.01}) -- ({2*\x-1.97},{\x+8-2.95}) -- ({2*\x-2.03},{\x+8-2.98}) -- ({2*\x-2.03},{\x+8-3.04}) -- ({2*\x-1.97},{\x+8-3.01});
}

\end{tikzpicture}
\caption{ The Voronoi tiling of the lattice $\Z^2$ in $\R^2$ endowed with the inner product of rank two described in Example~\ref{ex:quasi-tiling2}. The boundary of each Voronoi cell contains two closed segments in the middle of the two vertical sides, as well as the interior of the two other parallel sides. The vertices are excluded from the Voronoi cell.}
\label{fig:quasi-tiling2}
\end{figure}
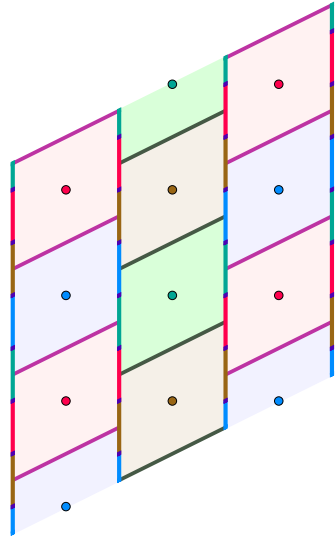

\begin{example} \label{ex:quasi-tiling2}
Consider $H=\R^2$ with the lattice $\L = \Z^2$ and rank two inner product $\highinn{}{\cdot\,,\cdot}$ given by 
\[\highinn{1}{x,y} = x_1  y_1 \textrm{ and } \highinn{2}{x,y} = (x_2 - \frac{1}{2}x_1)(y_2 - \frac{1}{2}y_1).\]
As shown in Figure~\ref{fig:quasi-tiling2}, the corresponding Voronoi decomposition is a tiling of $H$. The intersection of the closures of two Voronoi cells is not necessarily a common face. 
\end{example}

The following result is a consequence of Theorem~\ref{thm:Hausdorff_voronoi}.

\begin{thm}[Hausdorff convergence of Voronoi cells of admissible lattices of full rank]\label{thm:Hausdorff_voronoi-lattices} Let $\highinn{}{\cdot\, ,\cdot}\colon H\times H \to \R^r$ be an inner product and $\L \subset H$ be an admissible lattice of full rank. Let $\innone{t}{\cdot\, , \cdot}$, $t\in \R_+$, be a family of scalar products on $H$ which tamely degenerate to $\highinn{}{\cdot\,,\cdot}$ with parameters $\underline L_t \in \R_+^r$. 

Then, for each $\gamma \in \L$, the Voronoi cells $\Vor_{\L, \, \innone{t}{\cdot\,, \cdot}}(\gamma)$, $t \in \R_+$, converge to $\overline \Vor_{\L}(\gamma)$ in the Hausdorff metric as $t$ goes to infinity.
\end{thm} 
\begin{proof} By Theorem~\ref{thm:VoronoiCellDecomposition-lattices}, $\overline\Vor_\L(\gamma)$ is compact. Therefore, Theorem~\ref{thm:Hausdorff_voronoi} applies.
\end{proof}


\section{Metric degenerations of tori}\label{sec:metric_degeneration_general_tori}
Throughout this section, let $H$ be a real vector space and $\highinn{}{\cdot\,, \cdot} \colon H\times H \to \R^r$ an inner product. Moreover, let $\L$ be an admissible lattice of full rank in $H$. We consider the torus
\[
\T\coloneqq \rquot{H}{\L}.
\]

\smallskip

Let $\innone{t}{\cdot\,, \cdot}$, $t \in \R_+$, be scalar products on $H$, and denote by $\norm{x}_t \coloneqq \sqrt{\innone{t}{x, x}}$, $t \in \R_+$, the corresponding norms. The scalar product $\innone{t}{\cdot\,, \cdot}$ on $H$ induces a Riemannian  metric $\varphi_t$ on the torus $\T$. The associated distance function $\dist_t \colon \T \times \T \to [0, +\infty)$ on $\T$ has the following explicit form. For $u, v \in \T$, we have
\begin{equation} \label{eq:DistanceFunctionExplicit}
\dist_t(u,v) = \min_{\substack{z \in H \\ z=u-v \textrm{ in }\T}} \norm{z}_t = \min_{\gamma \in \L} \norm{x-y-\gamma}_t,
\end{equation}
where $x \in H$ and $y \in H$ are any representatives of $u$ and $v$, respectively. We are interested in the degeneration of the metric space $(\T, \dist_t)$, when the scalar products $\langle \cdot\,, \cdot \rangle_t$ behave asymptotically as in Section~\ref{ss:tameness}. We thus consider a family of vectors $\underline L_t = (L_{t,1}, \dots, L_{t,r}) \in \R_+^r$, $t\in \R_+$, and assume that the family $\innone{t}{\cdot\,, \cdot}$, $t \in \R_+$, $\omega$-tamely degenerates to $\highinn{}{\cdot\,, \cdot}$ with parameters $\underline L_t$, $t \in \R_+$, so that \ref{item:tame-a}, \ref{item:weaktame-b} and \ref{item:tame-c} are verified. In particular, by Proposition~\ref{prop:tameness_basic_family}, every pullback family for $\highinn{}{\cdot \, , \cdot}$ tamely degenerates to $\highinn{}{\cdot \, , \cdot}$.

\subsection{Middle tori} \label{ss:middle_tori}
Consider the non-increasing filtration
\[\filter^1=H \supseteq \filter^2 \supseteq \dots \supseteq \filter^r \supseteq \filter^{r+1}=(0)\]
of $H$ and the non-increasing filtration
\[\L^1=\L \supseteq \L^2 \supseteq \dots \supseteq \L^r \supseteq \L^{r+1}=(0)\]
of $\L$, where $\L^j =  \L \cap \filter^j$. From these two filtrations, we obtain a non-increasing sequence of tori
\[\T^1 = \T \supseteq \T^2 \supseteq \dots \supseteq \T^r\supseteq \T^{r+1}=(0),\]
 with $\T^j \coloneqq \rquot{\filter^j}{\L^j}$, $j\in [r]$. Each $\T^j$ is a closed subset of $\T$.
 
By Theorem~\ref{thm:AdmissibleDiscreteLattices}, the lattice $\grm{}{j}\L = \rquot{\L^j}{\L^{j+1}}$ is of full rank in $\grm{}{j} H =\rquot{ \filter^j }{\filter^{j+1}}$ for every $j \in [r]$. This defines a torus
\[
\Theta_j \coloneqq  \rquot{\grm{}{j}H }{ \grm{}{j}\L }.
\]
We call $\Theta_j$ the $j$th {\em graded piece} of the torus $\T$ with respect to the inner product $\highinn{}{\cdot \,, \cdot}$. The terminology is justified as follows. The projection maps $\proj_j \colon \L^j \to \grm{}{j}\L$ and $\proj_j \colon \filter^j \to \grm{}{j}H$ give rise to a projection map $\proj_j \colon \T^j \to \Theta_j$ with $\ker(\proj_j) =  \T^{j+1}$. Thus, we obtain an isomorphism
\[
\rquot{\T^j}{\T^{j+1}} \cong \Theta_j.
\]

The scalar product $\highinn{j}{\cdot\, ,\cdot}\colon \grm{}{j}H  \times\grm{}{j} H \to \R$ induces a Riemannian metric $\varphi_j$ on $\Theta_j$. We denote the induced distance function by $\dist_j \colon \Theta_j \times \Theta_j \to [0, +\infty)$.

\subsection{Volume formula and eigenvalue degeneration}
Consider the almost orthogonal decomposition $H= \aplus_{j=1}^r H_j$. Denote by $n \coloneq \dim_\R(H)$ and $n_j \coloneq \dim_\R(\grm{}{j} H)$, $ j \in [r]$. We first prove the following result. 
\begin{thm} \label{thm:volume_degeneration} 
Let $\innone{t}{\cdot\,, \cdot}$, $t \in \R_+$, be a family of scalar products which $\omega$-tamely degenerates to $\highinn{}{\cdot\,, \cdot}$ with parameters $\underline L_t =(L_{t,1}, \dots, L_{t,r})$, $t \in \R_+$. Then, the volume of $\T$ with respect to the Riemannian metric $\varphi_t$, properly rescaled, has the asymptotic behavior
\[
\lim_{t \to \infty} \prod_{j=1}^r L_{t,j}^{-n_j/2}   \cdot \vol(\T, \phi_t) =  \prod_{j=1}^r \vol (\Theta_{j}, \varphi_j).
\]
\end{thm}
Choose a basis $e_1, \dots, e_n$ of the lattice $\L$ in $H$. Let $M_t \in \R^{n \times n}$ be the Gram matrix of the scalar product $\innone{t}{\cdot\,,\cdot}$ in the basis $e_1, \dots, e_n$:
\begin{align*}
M_t = (M_t(i,j))_{1 \le i,j \le n} \qquad  \textrm{ with } \qquad M_t(i,j) \coloneqq \innone{t}{e_i, e_j}.
\end{align*} 
The volume of $\T$ can be expressed as
\begin{equation} \label{eq:torus_volume_matrix}
\vol(\T, \phi_t) = \det(M_t)^{1/2} = \prod_{i=1}^n \lambda_i(M_t)^{1/2},
\end{equation}
where $0 < \lambda_1(M_t) \le \lambda_2(M_t) \le \dots \le \lambda_n(M_t)$ are the eigenvalues of $M_t$.

Recall that $n_k \coloneqq \dim_\R(H_k) = \dim_\R(\grm{}{k} H)$, $k\in [r]$. Define the index sets
\begin{align*}
I_k &\coloneqq \{m_k +1, \dots, m_{k+1}\}, \quad k \in [r],\\
m_k &\coloneq \sum_{j =1}^{k-1} n_j, \quad k\in [r+1].
\end{align*}
We can choose a basis $e_1, \dots, e_n$ of $\L$ so that for each $k \in [r]$, the vectors $e_i$, $i \in I_k$, belong to $\filter^k$ and their  projections $\proj_k(e_i)$, $i\in I_k$, form a basis of $\grm{}{k}\L$. Let $N_k \in \R^{I_k \times I_k}$ be the $n_k\times n_k$ Gram matrix of the scalar product $\highinn{k}{\cdot\,,\cdot}$ on $\grm{}{k}H$ in the basis $\proj_k(e_i)$, $i\in I_k$. The volume of the graded piece $\Theta_k$, $k \in [r]$, is given by
\begin{equation} \label{eq:torus_volume_matrix2}
\vol(\Theta_k, \varphi_k) = \det(N_k)^{1/2} = \prod_{j=1}^{n_k} \lambda_j(N_k)^{1/2},
\end{equation}
where $0 < \lambda_1(N_{k}) \le \dots \le \lambda_{n_k}(N_{k})$ are the eigenvalues of $N_k$.

 In order to prove Theorem~\ref{thm:volume_degeneration}, we will relate the eigenvalues of the matrices $M_t$ and $N_{k}$. More precisely, we will establish the following theorem.
\begin{thm}[Degeneration of eigenvalues]\label{thm:limiting_eigenvalues} Assumptions as in Theorem~\ref{thm:volume_degeneration}, we have
\begin{equation} \label{eq:limiting_eigenvalues}
\lim_{t \to \infty} \frac{\lambda_{m_k + i} (M_t)}{L_{t,k}} =\lambda_i(N_{k})
\end{equation}
for all $k \in [r]$ and any $1 \le i \le n_k$.
\end{thm}
\begin{proof} Fix $k \in [r]$. For $t >0$, consider the polynomial
\[
f_t(\lambda) \coloneqq \prod_{k' < k} L_{t,k'}^{-n_k} \cdot L_{t,k}^{-(n - m_k)} \cdot \det \Big ( M_t - L_{t,k} \lambda \id \Big ).
\]
Clearly, $\lambda \in \C$ is an eigenvalue of $M_t$ exactly when $\lambda/ L_{t,k}$ is a zero of $f_t$ (and the multiplicities coincide). By the properties  \ref{item:tame-a}, \ref{item:weaktame-b}  and \ref{item:tame-c} of $\innone{t}{\cdot\,, \cdot}$, the entries of the matrix $M_t \in \R^{n \times n}$ behave like
\begin{align*}
	M_t(i, i') &= O\Big(L_{t,\max\{k, k'\}} \Big), &&i \in I_k, \, i' \in I_{k'},\\
	M_t(i, i') &= L_{t,k} \Big ( N_{k} (i, i') + o(1) \Big), &&i, i' \in I_k,
\end{align*}
as $t$ tends to $\infty$. Moreover, since  $L_{t, k' +1} / L_{t,k'}= o(1)$ for all $k' \in [r]$, one readily obtains that
\[
f_t(\lambda)  = \det \begin{pmatrix}
N_{1} - \frac{L_{t,k}}{L_{t,1}} \lambda \id  & O(1) & \dots & O(1) & 0 & \dots & 0 \\
0 & N_{2} - \frac{L_{t,k}}{L_{t,2}} \lambda \id &  \dots & O(1) & 0 & \dots & 0  \\
0 & 0 & \dots & O(1) & 0 & \dots & 0 \\
\vdots & \vdots & \ddots &\vdots & \vdots & \ddots &\vdots   \\
0 & 0 & \dots & O(1) & 0 & \dots & 0 \\
0 & 0 & \dots & N_{k} - \lambda \id & 0 & \dots & 0 \\
0 & 0& \dots & 0 & - \lambda \id & \dots & 0 \\
\vdots & \vdots & \ddots & \vdots & \vdots & \ddots & \vdots\\
0 & 0 & \dots & 0 & 0 & \dots & - \lambda \id
\end{pmatrix} + o(1)
\]
as $t$ tends to $\infty$. It then follows that $f_t$ converges to the polynomial
\[
g(\lambda) \coloneqq \prod_{k' < k} \det(N_{k'}) \cdot \det \big (N_{k} - \lambda \id \big) \cdot (-\lambda)^{n - m_{k+1}}
\]
uniformly on compact subsets of $\C$ as $t$ tends to $\infty$. The roots of $g$ are the $n_k$ eigenvalues of $N_{k}$, counted with their respective multiplicities, and $\lambda = 0$, with multiplicity $n - m_{k+1}$. Thus, as $t$ tends to $\infty$, precisely $m_k$ roots of $f_t$ tend to infinity, $n_k$ roots converge to the eigenvalues of $N_{k}$, and the $(n-m_{k+1})$ remaining  ones converge to $\lambda = 0$. Taking into account this behavior for all $k \in [r]$, the relation \eqref{eq:limiting_eigenvalues} follows and the proof is complete.
\end{proof}

\begin{proof}[Proof of Theorem~\ref{thm:volume_degeneration}]
This follows by combining Theorem~\ref{thm:limiting_eigenvalues} with \eqref{eq:torus_volume_matrix} and \eqref{eq:torus_volume_matrix2}.
\end{proof}

\subsection{Metric degeneration} In this section, we study the degeneration of the flat tori $(\T, \dist_t)$, $t \in \R_+$, from the perspective of Gromov--Hausdorff convergence.

 For a metric space $(X,\rho)$, we denote, as in Section~\ref{sec:Hausdorff_convergence}, by $\rho_{\HA}(A,B)$ the Hausdorff distance between two non-empty subsets $A, B \subset X$
\[
\rho_{\HA}(A,B) \coloneqq \max \big \{ \sup_{a \in A} \rho(a, B), \, \sup_{b \in B} \rho(b, A) \big \}, \qquad \rho(a,B) = \inf_{b \in B}\rho(a,b), \quad  \rho(b,A) = \inf_{a \in A}\rho(a,b).
\]
Recall that the {\em Gromov--Hausdorff distance} between two metric spaces $X$ and $Y$ is defined by
\begin{align*}
\dist_{\GH}(X,Y) = \inf \big \{ \varrho_{\HA}(\iota_X(X), \iota_Y(Y)) \},
\end{align*}
where the infimum is taken over all metric spaces $(Z,\varrho)$ and isometric embeddings $\iota_X \colon X \hookrightarrow Z$ and $\iota_Y \colon Y \hookrightarrow Z$. A family of metric spaces $(X_t)_{t \in \R_+}$ {\em converges to a metric space $X$ in the Gromov--Hausdorff sense} if $\dist_{\GH}(X_t,X) \to 0$ for $t \to \infty$. In this situation, we also write $X_t \xrightarrow{\GH} X$ for $t \to \infty$.

\begin{thm} \label{thm:GHConvergenceAbstract}
Let $\innone{t}{\cdot\,, \cdot}$, $t \in \R_+$, be a family of scalar products which $\omega$-tamely degenerates to $\highinn{}{\cdot\,, \cdot}$ with parameters $\underline L_t =(L_{t,1}, \dots, L_{t,r})$.

Then,
\[
\Big (\T, \frac{1}{\sqrt{L_{t,1}}} \, \dist_t \Big ) \xrightarrow{\GH} \Big (\Theta_{1}, \dist_{1} \Big )
\]
as $t \to \infty$. More generally, for all $j=1, \dots, r$,
\[
\Big ( \T^j,\frac{1}{\sqrt{L_{t,j}}} \,\dist_t \rest{\T^j\times \T^j} \Big ) \xrightarrow{\GH} \Big ( \Theta_{j}, \dist_{j} \Big )
\]
as $t \to \infty$. Here, we view $\T^j$ as a closed subset $\T^j \subset \T$, endowed with the restriction $\dist_t \rest{\T^j\times \T^j}$ of $\dist_t$ to $\T^j\times \T^j$.
\end{thm}
We need the following result (see e.g.~\cite[Example 7.4.4]{Burago}). If $(\varrho_t)_{t \in \R_+}$ is a family of metrics on $X$ and $\rho$ is a semimetric on $X$ such that $\lim_{t \to \infty} \varrho_t = \rho$ uniformly on $X \times X$, then $(X, \varrho_t)$ converges to the quotient $\rquot{X}{\rho}$ in the Gromov--Hausdorff sense. The quotient metric space $\rquot{X}{\rho}$ is obtained by identifying points $x,y \in X$ with $\rho(x,y) = 0$. We will prove Theorem~\ref{thm:GHConvergenceAbstract} by applying this result to the distance functions $\dist_t$, $t \in \R_+$.

For $j = 1, \dots, r$, consider the semi-metric $\rho_j \colon \T^j \times \T^j \to [0, \infty)$ given by
\[
\rho_j(x,y) \coloneqq \dist_{j}(\proj_j(x), \proj_j(y)),  \qquad x,y \in \T^j,
\]
where $\proj_j \colon \T^j \to \Theta_j$ is the natural projection map from $\T^j = \rquot{\filter^j}{\L^j}$ onto $\Theta_{j} = \rquot{\grm{}{j} H}{\grm{}{j} \L}$. Note that the quotient $\T^j / \rho_j$ is isometric to $(\Theta_j, \dist_j)$.

\begin{lem} \label{lem:ConvergenceDistanceFunction}
Let $\innone{t}{\cdot\,, \cdot}$, $t \in \R_+$, be a family of scalar products which $\omega$-tamely degenerates to $\highinn{}{\cdot\,, \cdot}$ with parameters $\underline L_t =(L_{t,1}, \dots, L_{t,r})$, $t \in \R_+$. Then,
\[
\lim_{t \to \infty} \frac{1}{\sqrt{L_{t,1}}} \, \dist_t = \rho_1
\]
uniformly on $\T \times \T$. More generally, for all $j =1, \dots, r$,
\[
\lim_{t \to \infty} \frac{1}{\sqrt{L_{t,j}}} \,\dist_t \rest{\T^j\times \T^j} = \rho_j
\]
uniformly on $\T^j \times \T^j$.
\end{lem}
\begin{proof}
It suffices to treat the case $j=1$. Indeed, as follows from Lemma~\ref{lem:finiteness-discrete} and \eqref{eq:DistanceFunctionExplicit}
, for sufficiently large $t$, the restriction $\dist_t \rest{\T^j\times \T^j}$ of $\dist_t$ to $\T^j = \filter^j / \L^j$ (viewed as a closed subset $\T^j \subset\T$) coincides with the metric stemming from the scalar product $\innone{t}{\cdot\,, \cdot} \rest{\filter^j} \colon \filter^j \times \filter^j \to \R$. Moreover, the family of scalar products $\innone{t}{\cdot\,, \cdot}  \rest{\filter^j}$,  $t \in \R_+$, converges $\omega$-tamely to the inner product $[\cdot \, , \cdot] \colon \filter^j \times\filter^j \to \R^{r-j+1}$ given by $[\cdot \, , \cdot] = (\highinn{j}{\cdot\, ,\cdot}, \dots, \highinn{r}{\cdot\, ,\cdot})$. Thus, the result for $\T^j$, $j >1$ can be reduced to the case $j=1$.

Let $H = \aplus_{j=1}^r H_j$ be the almost orthogonal decomposition of $H$. We decompose $x \in H$ into $x= x_1 + \dots + x_r$ with $x_j \in H_j$, $j \in [r]$, and set $\norm{x}_1 \coloneqq \sqrt{\highinn{1}{x,x}}$ for $x \in H_1$. The distance between two points $u,v \in \T$ can then be rewritten as
\[
\rho_1(u, v) = \min_{\gamma \in \L} \norm{x_1 - y_1 - \gamma_1}_1,
\]
where  $x, y \in H$ are any fixed representatives of $u$ and $v$, respectively. The assumptions \ref{item:tame-a}, \ref{item:weaktame-b} and \ref{item:tame-c} imply that for every compact subset $B \subset H$ and every finite set $\S \subset \L$, 
\[
\frac{1}{\sqrt{L_{t,1}}} \min_{\gamma\in \S} \norm{x - y - \gamma}_t \to \min_{\gamma \in \S} \norm{x_1 - y_1 - \gamma_1}_1
\]
uniformly for $x,y\in B$. Fix a compact subset $B \subset H$ which covers $\T$ via the projection map $\proj \colon H \to \T$. By Lemma~\ref{lem:finiteness-discrete}, there exists a finite set $\S_B \subset \L$ such that
\[
\min_{\gamma \in \L} \norm{x - y - \gamma}_t = \min_{\gamma \in \S_B} \norm{x - y -\gamma}_t \quad \text{ and } \quad \min_{\gamma \in \L} \norm{x_1 -y_1 - \gamma_1}_1 = \min_{\gamma \in \S_B} \norm{x_1 - y_1 -\gamma_1}_1
\]
for all $x,y \in B$ and large $t$. Fixing representatives $x= x(u)$ and $y = y(v)$ in $B$ for all $u, v \in \T$,  we conclude that
\[
\left| \frac{1}{\sqrt{L_{t,1}}} \dist_t(u,v) - \rho_1(u,v)\right|  = \left| \frac{1}{\sqrt{L_{t,1}}} \min_{\gamma \in \S_B} \norm{x - y -\gamma}_t - \min_{\gamma \in \S_B} \norm{x_1 - y_1 -\gamma_1}_1\right| \to 0 
\]
uniformly for $u,v \in \T$ as $t \to \infty$. The proof is complete. \qedhere
\end{proof}
\begin{proof}[Proof of Theorem~\ref{thm:GHConvergenceAbstract}]
By Lemma~\ref{lem:ConvergenceDistanceFunction}, the metrics $\dist_t \rest{\T^j \times \T^j}$ converge to the semimetric $\rho_j$ uniformly on $\T^j \times \T^j$ for $t \to \infty$. This implies that $(\T, \dist_t \rest{\T^j \times \T^j})$ is Gromov--Hausdorff convergent to the quotient $\T^j / \rho_j$, which is isometric to $(\Theta_{j}, \dist_{j})$.
\end{proof}

\subsection{Metric degeneration for pullback families}
For pullback families, we obtain more precise results on the metric degeneration of the torus.

Let $\innone{t}{\cdot\,, \cdot}$, $t \in \R_+$, be a pullback family of $\highinn{}{\cdot\,, \cdot}$ with parameters $\underline L_t =(L_{t,1}, \dots, L_{t,r})$, $t \in \R_+$. Recall that
\[
\innone{t}{\cdot\,,\cdot}= \innone{\underline L_t}{\cdot\,,\cdot} = L_{t,1} \highinn{1}{\cdot\,,\cdot} + \dots + L_{t,r} \highinn{r}{\cdot\,,\cdot}, \qquad t\in \R_+,
\]
and $\lim_{t \to \infty}L_{t,j}/L_{t,j+1} = + \infty$ for all $j= 1, \dots, r-1$.

We introduce the following function
\begin{align*}
\dtor = (\dtor_1, \dots, \dtor_r) \colon \T\times \T \to \R^r \\
\dtor(p,q) \coloneqq \min_{\substack{x \in H\\ x = p-q \textrm{ in } \T}} \highinn{}{x,x}.
\end{align*}
It follows from Theorem~\ref{thm:voronoi-discrete} that the above minimum exists. Thus, $\zeta$ is well-defined. Note that $\zeta(p,q) \gexeq 0$ for all $p,  q \in  \T$ and $\zeta(p,q) =0$ exactly when $p=q$.
\begin{remark}
The definition of $\dtor \colon  \T\times \T \to \R^r$ is similar to the distance function $\dist_t \colon \T \times \T \to \R$ induced by the scalar product $\innone{t}{\cdot\,,\cdot}$, see \eqref{eq:DistanceFunctionExplicit}. Thus, one may interpret $\zeta$ as a higher rank version of the (square of the) distance function.
\end{remark}

Let $\Vor_\L(0)$ be the  Voronoi cell of the origin for $\highinn{}{\cdot\,,\cdot}$ and denote by $\interior{\Vor_\L(0)}$ its interior. The projection map $\proj \colon H \to \T$ is injective on $\interior{\Vor_\L(0)}$. In the following, we identify $\interior{\Vor_\L(0)}$ with a subset of $\T$. 

\begin{thm}\label{thm:metric-asymptotics-pullback} Let $\innone{t}{\cdot\,, \cdot}$, $t \in \R_+$, be a pullback family of $\highinn{}{\cdot\,, \cdot}$ with parameters $(\underline L_t)_{t \in \R_+}$.

\begin{itemize}
\item [(i)]
Consider a fixed pair of points $p,q \in \T$. Then,
\begin{equation} \label{eq:EqualityHigherRankMetric}
\dist_t(p,q)  = \sqrt{L_{t,1} \dtor_1(p,q) + \dots + L_{t,r} \dtor_r(p,q)}
\end{equation}
for all large $t \in \R_+$.
\item [(ii)]
Fix a compact subset $B \subset \interior{\Vor_\L(0)} \subset \T$ and define $K = \{(p,q) \in \T\times \T \st p-q \in B \}$. Then, $K \subset \T \times \T$ is compact and for all large $t \in \R_+$, we have 
\begin{equation} \label{eq:EqualityHigherRankMetric2}
\dist_t(p,q) = \sqrt{L_{t,1} \dtor_1(p,q) + \dots + L_{t,r} \dtor_r(p,q)}, \qquad \text{for all } (p,q) \in K.
\end{equation}
\end{itemize}
\end{thm} 
\begin{proof} 
In order to prove the claim in $(i)$, recall that
\[
\dist_t(p,q)^2 = \min_{\substack{x \in H\\ x = p-q \textrm{ in } \T}} \norm{x}_t^2 = \min_{x \in S} \norm{x}_t^2,
\]
where $S = \{x \in H \st x= p-q \text{ in } \T \}$. Fix an element $x \in S$ such that $\highinn{}{x,x} = \zeta(p,q)$. Then
\[
 \norm{x}_t^2 = \innone{t}{x,x} = L_{t,1} \highinn{1}{x,x} + \dots +L_{t,r} \highinn{r}{x,x} = L_{t,1} \zeta_1(p,q) + \dots +L_{t,r} \zeta_r(p,q),
\]
which is precisely the square of the right-hand side in \eqref{eq:EqualityHigherRankMetric}. Thus, it suffices to prove that
\[
\min_{y \in S} \norm{y}_t^2 = \norm{x}_t^2
\]
for all large $t \in \R_+$. Since the lattice $\L$ is admissible, $S$ is an admissible discrete subset of $H$. By Lemma~\ref{lem:finiteness-discrete},
there exists a finite subset $\S_0 \subset S$ such that $x \in S$ and
\[
\min_{y \in S} \norm{y}_t^2 = \min_{y \in \S_0} \norm{y}_t^2
\]
for all large $t \in \R_+$. It thus suffices to prove that for each fixed $y \in \S_0$ with $\highinn{}{y,y} \neq \zeta(p,q)$, we have
\[
 \norm{y}_t^2 > \norm{x}_t^2
\]
for all large $t \in \R_+$. Indeed, since $\S_0$ is finite, this implies that $\dist_t(p,q)^2 = \min_{y \in S} \norm{y}_t^2 = \min_{y \in \S_0} \norm{y}_t^2 = \norm{x}_t^2$ for all large $t \in \R_+$.

Fix $y \in S$ with $\highinn{}{y,y} \neq \zeta(p,q)$. Let $j \in [r]$ be the smallest index with $\highinn{j}{y,y} \neq \zeta_j(p,q)$. By the definition of $\zeta(p,q)$, we have $\zeta(p,q) \lexst \highinn{}{y,y}$ and hence $\zeta_j(p,q) < \highinn{j}{y,y}$. Therefore,
\begin{align*}
\innone{t}{y,y} - \innone{t}{x,x} &= \sum_{k=1}^r L_{t,k}(\highinn{k}{y,y} - \highinn{k}{x,x}) = \sum_{k=1}^r L_{t,k}(\highinn{k}{y,y} - \zeta_k(p,q)) \\
&= L_{t,j} \Big(\highinn{j}{y,y} - \zeta_j(p,q)  + \sum_{k=j+1}^r \frac{L_{t,k}}{L_{t,j}} (\highinn{k}{y,y} - \zeta_k(p,q)) \Big) >0,
\end{align*}
where the last estimate follows from the fact that $\lim_{t \to \infty} \frac{L_{t,k+1}}{L_{t,k}} = 0$ for all $k$.

It remains to prove the claim in $(ii)$. Since the addition $\T\times \T \to \T \colon (p,q) \mapsto p-q$, is continuous, the set $K \subset \T \times \T$ is compact. By assumption, $B$ lies in the interior of the Voronoi cell $\Vor_\L(0)$. It follows that for some $\epsilon\in (0,1)$, we have $B \subset U_\varepsilon$ for the set $U_\varepsilon = (1-\varepsilon) \overline \Vor_{S}(0)$ in Lemma~\ref{lem:inclusion_voronoi}.  As in Section~\ref{sec:Hausdorff_convergence}, let $W_t$ be the Voronoi cell of the origin for the scalar product $\innone{t}{\cdot\,,\cdot}$. Applying Lemma~\ref{lem:inclusion_voronoi}, we have $B \subset W_t$  for all large $t$.  

For a pair $(p,q) \in K$, choose a representative $x \in B$ of $p-q$.  Since $x$ lies in $W_t$ and $ \Vor_\L(0)$, we have
\[
 \dist_t(p,q)^2 = \min_{\gamma \in \L} \norm{x-\gamma}_t^2 = \norm{x}_t^2 \qquad \text{and} \qquad \zeta(p,q) =  \min_{\gamma \in \L} \highinn{}{x-\gamma,x-\gamma} = \highinn{}{x,x}.
\]
The second statement in the theorem now follows by observing that
\[
\dist_t(p,q)^2 = \norm{x}_t^2 = L_{t,1} \highinn{1}{x,x} + \dots + L_{t,r}\highinn{r}{x,x} =  L_{t,1} \zeta_1(p,q)+ \dots + L_{t,r}\zeta_r(p,q). \qedhere
\]
\end{proof}

\begin{remark} Note that the distance $\dist_t(p,q)$ between two points $p,q \in \T$ can become infinite as $t$ tends to infinity. In such a case, Theorem~\ref{thm:metric-asymptotics-pullback} gives the precise asymptotics.
\end{remark}

\begin{remark}
One can obtain similar results for families of scalar products which behave similar to pullback families. For instance, if $\innone{t}{\cdot \,, \cdot}$, $t \in \R_+$, is a tamely degenerating family of scalar products with parameters $(\underline L_t)_{t \in \R_+}$ for $\highinn{}{\cdot, \cdot}$ such that, additionally, $L_{t, r} \equiv 1$ and
\[
\lim_{t \to \infty} \left(\innone{t}{x,y} - \innone{\underline L_t}{x,y}\right) = 0
\]
for all $x,y \in H$, then, we get
\[
\lim_{t \to \infty} \left(\dist_t(p,q) - \sqrt{L_{t,1} \dtor_1(p,q) + \dots + L_{t,r} \dtor_r(p,q)}\right) = 0
\]
for all fixed pairs $(p,q) \in \T \times \T$. Since the proof of these results becomes more technical, we decided not to include them.
\end{remark}


\section{Geometric applications} \label{sec:tropical_curves}
In this section, we discuss geometric applications of our results in two situations. The first one concerns the multi-scale geometry of degenerating metric graphs, and the second one concerns degenerating Riemann surfaces. We treat the former in detail. The latter is described here informally because its treatment requires the introduction of specific tools for degenerations of Riemann surfaces, which is not the main objective of this paper. This will appear in our forthcoming work~\cite{AN-AG-hybrid}, see also Remark~\ref{rem:OutlookAGHybrid}.

\subsection{Metric graphs and tropical curves} \label{ss:degeneration-metric-graphs} In this section, we apply our results to metric graphs and their geometric limits, called tropical curves of higher rank.

\subsubsection{Metric graphs}
Let $G=(V, E)$ be a finite connected graph with vertex set $V$ and edge set $E$. Let $l \colon E \to (0,+\infty)$ be a length function, which assigns a positive real number $l(e)$ to every edge $e \in E$. The pair $(G, l)$ is called {\em a weighted graph}. To a weighted graph $(G,l)$, we associate a metric space $\mgr$ as follows. We insert an interval $\Ical_e= [0,l(e)]$ of length $l(e)$ between the two extremities $u$ and $v$ of any edge $e=\{u,v\}$ in the graph (with the two extremities of $\Ical_e$ identified with the two extremities $u$ and $v$ of $e$). The obtained space $\mgr$ carries a natural quotient topology which is metrizable by the {\em path metric} $\rho$: the distance $\rho(x,y)$ between $x,y \in \mgr$ is the arc length of the shortest path connecting them.

A {\em metric graph} is a compact metric space $\mgr$ which arises in this way for some weighted graph $(G, l)$. In this case, $\mgr$ is called the {\em metric realization} of $(G, l)$, and the pair $(G, l)$ is called a {\em model} of $\mgr$. If we want to emphasize that a metric graph $\mgr$ has a model with an underlying graph $G$, then we call it a {\em metric graph over $G$}.  

A graph $G$ is called \emph{essential} if either it has no vertices of degree two, or it is the graph consisting of one vertex and one loop edge. Any metric graph $\mgr$ has a model $(G,l)$ in which $G$ is essential. Such a model is called an \emph{essential model} for $\mgr$. All essential models of $\mgr$ have the same underlying essential graph $G$. Moreover, the essential model $(G,l)$ for a metric graph $\mgr$ is unique unless $G$ has non-trivial automorphisms, in which case, $\mgr$ has a finite number of essential models $(G, l)$. (These are obtained by applying the automorphisms $\varphi$ to one fixed model $(G,l)$, that is, considering edge length functions given by compositions $\ell \circ \varphi$.)

Let $H = H_1(\mgr, \R)$ and $\L =H_1(\mgr, \Z)$ be the first homology groups of $\mgr$ with real and integral coefficients, respectively.  Given a model $(G, l)$ of $\mgr$, we can identify $H = H_1(\mgr, \R)$ and $\L =H_1(\mgr, \Z)$ with the spaces of real-valued and integer-valued flows on $G$.  Let $\E$ be the set of size $2\abs E$ obtained by replacing each edge $e =\{u,v\}$ of $E$ by two oriented edges $uv$ and $vu$ (including loops for which $u=v$).  By an abuse of the notation, we also denote elements of $\mathbb E$ by $e$, and use $\bar e$ when referring to the same edge but with the opposite orientation. A {\em real-valued flow} on $G$ is a map $x\colon \E \to \R$ which verifies
\begin{align*}
x(e) = -x(\bar e) \qquad \forall\,\, e\in \E, \qquad \textrm{ and } \qquad \sum_{e =vu \in \E} x(e) =0\qquad \forall\,\, v\in V.
\end{align*}
If $x(e) \in \Z$ for all $e \in \E$, then $x$ is called {\em an integer-valued flow}. The spaces of real- and integer-valued flows on $G$ are denoted by $H_1(G, \R)$ and $H_1(G, \Z)$, respectively. We can then naturally identify $H_1(\mgr, \R) = H_1(G, \R)$ and $H_1(\mgr, \Z) = H_1(G, \Z)$. In the following, we tacitly make this identification and make it explicit, whenever there is a risk of confusion.

The dimension $\graphgenus = \dim(H_1(\mgr, \R))$ is called {\em the genus} of $\mgr$. Given a model $(G,l)$ of $\mgr$ with vertex set $V$ and edge set $E$, we have $\graphgenus = |E| - |V| +1$ by Euler's formula.

\smallskip

Recall the definition of polarization on metric graphs~\cite{BLN97, KS00, BR07, MZ08, CV10}.  Fix a model $(G, l)$ of $\mgr$. Fixing an orientation $\orient \colon E \to \E$, we identify $E$ with a subset of $\E$.  The {\em polarization} on $H$ is the scalar product $\innone{\mgr}{\cdot \, , \cdot} \colon H \times H \to \R$ defined by
\[
\innone{\mgr}{x, y} \coloneqq \sum_{e\in E}l(e)x(e)y(e) \qquad \forall x, y\in H.
\] 
It does not depend on the choices of the model and orientation. If the graph $G$ is fixed, then sometimes we denote $\innone{\mgr}{\cdot\,, \cdot}$ by $\innone{l}{\cdot\,, \cdot}$.

\smallskip

The {\em Jacobian} of $\mgr$ is the $\graphgenus$-dimensional torus
\[
\Jac(\mgr) \coloneqq  \rquot{H}{\L}=\rquot{H_1(\mgr,\R)}{H(\mgr, \Z)}
\]
We endow $\Jac(\mgr)$ with the distance function $\dist_\mgr \colon \Jac(\mgr)\times\Jac(\mgr)\to [0, \infty)$ induced by the polarization $\innone{\mgr}{\cdot\,, \cdot}$.

\subsubsection{Tropical curves} \label{ss:TropicalCurves}
In the following, we recall the definition of tropical curves in the sense of \cite{AN,AN2}. This framework was introduced in \cite{AN,AN2} for answering analytic questions on degenerating metric graphs and Riemann surfaces. Tropical curves in this sense should be viewed as multi-scale limits of metric graphs. 

\smallskip

Let $G=(V,E)$ be a finite, connected graph. We assume that $G$ is essential, which will be convenient when realizing tropical curves as limits of metric graphs. A {\em layering} of $G$ is an ordered sequence $\pi = (\pi_1,  \dots, \pi_r, \pi_\fin)$ of pairwise disjoint subsets $\pi_j \subset E$, $j \in \{1, \dots, r, \fin\}$, such that $\pi_j  \neq \varnothing$ for $j \in \{1,\dots, r\}$ and
\[
E = \bigsqcup_{j \in \{1, \dots, r, \fin \}} \pi_j.
\]
We stress that we allow $\pi_\fin = \varnothing$. Equivalently, a layering of $G$ is a pair $\pi = (\pi_\infty, \pi_\fin)$ of a subset $\pi_\fin \subset E$ and an ordered partition $\pi_\infty=(\pi_1, \dots, \pi_r)$ of $E \setminus \pi_\fin$.

The integer $r \in \{0\} \cup \N $ is called the {\em (infinitary) rank} of $\pi$. The sets $\pi_j$, $j \in \{1, \dots, r\}$, are the {\em infinite layers} of $\pi$. The last set $\pi_\fin$ is called the {\em finite layer} of $\pi$. The pair $(G, \pi)$ is called {\em a layered graph}. The meaning of this terminology will become clear below, when we relate layered graphs to degenerations of metric graphs. 

A {\em weighted layered graph} $(G, \pi, \ell)$ is a layered graph $(G,\pi)$, with $G$ essential, together with an edge length function $\ell \colon E \to (0, \infty)$ satisfying the {\em normalization condition}
\[
\sum_{e \in \pi_j} \ell(e) = 1, \qquad \text{ for all } j =1, \dots,  r.
\]
To a weighted layered graph $(G, \pi, \ell)$, we associate the metric realization $\curve$ of $(G, \ell)$ and the following decreasing sequence
\begin{equation} \label{eq:DecreasingSubgraphs}
\curve = \curve^1 \supset \curve^2 \supset \dots \supset \curve^r \supseteq \curve^{\fin}
\end{equation}
of metric subgraphs: for $j \in \{1, \dots, r, \fin\}$, the metric subgraph $\curve^j$ of $\curve$ is obtained by removing all interval edges $e \in \pi_1\dots \cup \pi_{j-1}$ from $\curve$.  Equivalently, $\curve^j$ is the metric realization of the graph $G^j = (V, E \setminus (\pi_1 \cup \dots \cup \pi_{j-1}))$ with the edge length function $\ell\rest{\pi_{j} \cup\dots \pi_r \cup  \pi_\fin}$.

{\em A tropical curve}  is a pair $(\curve, \curve^\bullet)$ of a metric graph $\curve$ and a decreasing sequence of metric subgraphs $\curve^\bullet = (\curve^j)_{j \in \{1, \dots, r, \fin\}}$ with $\curve = \curve^1$ that arises in this way for some triple $(G, \pi, \ell)$, with the graph $G$ essential.

\begin{remark} We remark that the definition of a tropical curve here differs slightly from the one given in \cite{AN2}, in the sense that the graph $G$ was not assumed to be essential there. The results in \cite{AN2} are motivated by applications to moduli spaces of Riemann surfaces. In that setting, we allowed to have vertices of degree two in order to treat the case where that vertex represents a Riemann surface of positive genus in the dual graph. In order to simplify the presentation, we do not consider vertex-weighted graphs here. Putting the essential assumption then allows to define a meaningful convergence of metric graphs toward tropical curves, see below.
\end{remark}

By an abuse of the notation, we simply write $\curve$ for the pair $(\curve, \curve^\bullet)$. In this case, the tropical curve $\curve$ is called {\em the realization} of the weighted layered graph $(G, \pi, \ell)$, and the triple $(G, \pi, \ell)$ is called a {\em model} of the tropical curve $\curve$. All models of a tropical curve $\curve$ have the same underlying essential graph $G$.  A tropical curve $\curve$ has a unique model unless the graph $G$ has a non-trivial automorphism group, in which case, $\curve$ has a finite number of models with the same underlying graph $G$. If we want to emphasize that a tropical curve $\curve$ has a model $(G,\pi, \ell)$ with underlying graph $G$, then we call it a {\em tropical curve over $G$}. 

\smallskip

The {\em (infinitary) rank of a tropical curve} $\curve$ is by definition the rank of any defining layering $\pi$, or equivalently, the number of metric graphs in the sequence $\curve^\bullet$ minus one. We note that with this terminology, metric graphs are precisely tropical curves of (infinitary) rank $r = 0$.

\smallskip

A tropical curve $\curve$ of rank $r$ defines a collection of metric graphs $\Gamma^1, \dots, \Gamma^r, \Gamma^\fin$, called the {\em graded minors} of $\curve$. The $j$th graded minor $\Gamma^j$ is the metric graph obtained by contracting all interval edges $e \in \pi_{j+1} \dots \cup \pi_{r} \cup \pi_\fin$ in $\curve^ j$. Equivalently, $\Gamma^j$ is the metric realization of the $j$th \emph{graded minor} $\grm{}{j} G$ of $G$ defined by  
\[
\grm{}{j} G\coloneqq G^j / (\pi_{j+1} \dots \cup \pi_{r} \cup \pi_\fin)
\] 
endowed with edge length function $\ell\rest{\pi_{j}}$. 

(For a finite graph $R = (V(R), E(R))$ and an edge set $F\subset E(R)$, the graph $R / F$ is the contraction of $R$ along $F$.)

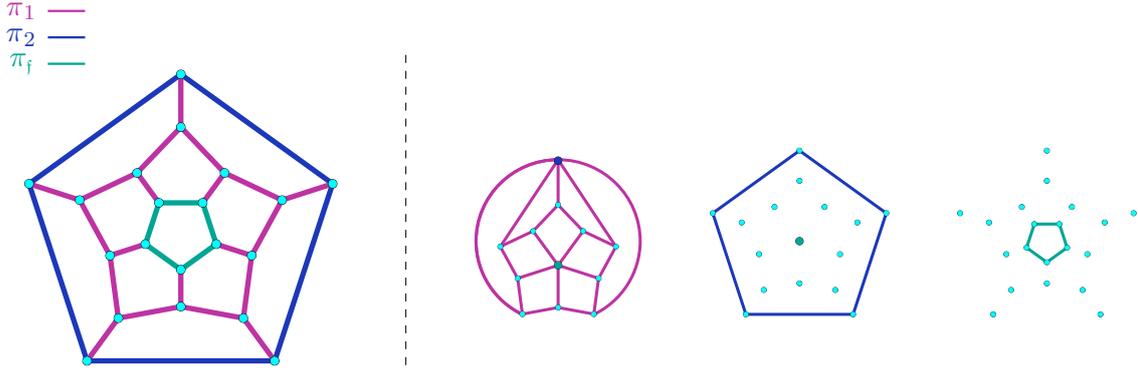
\begin{figure}
\centering
\begin{tikzpicture}[scale=.7]
\node[byzant] at (-3, 4.2) {$\pi_1$};
\draw[line width=0.30mm, byzant] (-2.5, 4.2) -- (-1.8, 4.2);
\node[pblue] at (-3, 3.7) {$\pi_2$};
\draw[line width=0.30mm, pblue] (-2.5, 3.7) -- (-1.8, 3.7);
\node[pgreen] at (-3, 3.2) {$\pi_\fin$};
\draw[line width=0.30mm, pgreen] (-2.5, 3.2) -- (-1.8, 3.2);

\foreach \x in {1,...,5}
{
\draw[line width=0.70mm, pblue] ({\x*72+18}:3) -- ({\x*72+90}:3);
\draw[line width=0.70mm, pgreen] ({\x*72-18}:.7) -- ({\x*72-90}:.7);

}
\foreach \x in {1,...,5}
{
\draw[line width=0.70mm, byzant] ({\x*72-90}:1.4) -- ({\x*72-90}:.7);
\draw[line width=0.70mm, byzant] ({\x*72+90}:3) -- ({\x*72+90}:2);
\draw[line width=0.70mm, byzant] ({\x*72+90}:2) -- ({\x*72+54}:1.4);
\draw[line width=0.70mm, byzant] ({\x*72+90}:2) -- ({\x*72+126}:1.4);
}
\foreach \x in {1,...,5}
{
\draw[line width=0.1mm] ({\x*72+90}:3) circle (0.8mm);
\filldraw[aqua] ({\x*72+90}:3) circle (0.7mm);
\draw[line width=0.1mm] ({\x*72+90}:2) circle (0.8mm);
\filldraw[aqua] ({\x*72+90}:2) circle (0.7mm);
\draw[line width=0.1mm] ({\x*72-90}:.7) circle (0.8mm);
\filldraw[aqua] ({\x*72-90}:.7) circle (0.7mm);
\draw[line width=0.1mm] ({\x*72-90}:1.4) circle (0.8mm);
\filldraw[aqua] ({\x*72-90}:1.4) circle (0.7mm);
}
\end{tikzpicture}
\qquad
\begin{tikzpicture}[scale=.7]
\draw[line width=0.1mm, dashed] (0,0) -- (90:6);
\end{tikzpicture}
\qquad
\begin{tikzpicture}[scale=.4]
\node at (0, -3) { };
\foreach \x in {1,...,5}
{
\draw[line width=0.40mm, byzant] ({\x*72-90}:1.4) -- (0,0);
\draw[line width=0.40mm, byzant] ({\x*72+90}:2) -- ({\x*72+54}:1.4);
\draw[line width=0.40mm, byzant] ({\x*72+90}:2) -- ({\x*72+126}:1.4);
}

\draw[line width=0.40mm, byzant] (90:3.5) -- (90:2);
\draw[line width=0.40mm, byzant] (90:3.5) -- (162:2);
\draw[line width=0.40mm, byzant] (90:3.5) -- (18:2);
\draw[line width=0.40mm, byzant] (90:3.5) arc (90:245:2.7);
\draw[line width=0.40mm, byzant] (90:3.5) arc (90:245:2.7);
\draw[line width=0.40mm, byzant] (90:3.5) arc (90:-65:2.7);

\draw[line width=0.1mm] (90:3.47) circle (1.2mm);
\filldraw[pblue] (90:3.47) circle (1.1mm);

\draw[line width=0.1mm] (0,0) circle (1.2mm);
\filldraw[pgreen] (0,0) circle (1.1mm);

\foreach \x in {1,...,5}
{
\draw[line width=0.1mm] ({\x*72+90}:2) circle (0.8mm);
\filldraw[aqua] ({\x*72+90}:2) circle (0.7mm);
\draw[line width=0.1mm] ({\x*72-90}:1.4) circle (0.8mm);
\filldraw[aqua] ({\x*72-90}:1.4) circle (0.7mm);
}
\end{tikzpicture}
\qquad
\begin{tikzpicture}[scale=.4]
\node at (0, -3.8) { };
\foreach \x in {1,...,5}
{
\draw[line width=0.40mm, pblue] ({\x*72+18}:3) -- ({\x*72+90}:3);
}
\draw[line width=0.1mm] (0,0) circle (1.2mm);
\filldraw[pgreen] (0,0) circle (1.1mm);

\foreach \x in {1,...,5}
{
\draw[line width=0.1mm] ({\x*72+90}:3) circle (0.8mm);
\filldraw[aqua] ({\x*72+90}:3) circle (0.7mm);
\draw[line width=0.1mm] ({\x*72+90}:2) circle (0.8mm);
\filldraw[aqua] ({\x*72+90}:2) circle (0.7mm);
\draw[line width=0.1mm] ({\x*72-90}:1.4) circle (0.8mm);
\filldraw[aqua] ({\x*72-90}:1.4) circle (0.7mm);
}
\end{tikzpicture}
\qquad
\begin{tikzpicture}[scale=.4]
\node at (0, -3.8) { };
\foreach \x in {1,...,5}
{
\draw[line width=0.40mm, pgreen] ({\x*72-18}:.7) -- ({\x*72-90}:.7);
}
\foreach \x in {1,...,5}
{
\draw[line width=0.1mm] ({\x*72+90}:3) circle (0.8mm);
\filldraw[aqua] ({\x*72+90}:3) circle (0.7mm);
\draw[line width=0.1mm] ({\x*72+90}:2) circle (0.8mm);
\filldraw[aqua] ({\x*72+90}:2) circle (0.7mm);
\draw[line width=0.1mm] ({\x*72-90}:.7) circle (0.8mm);
\filldraw[aqua] ({\x*72-90}:.7) circle (0.7mm);
\draw[line width=0.1mm] ({\x*72-90}:1.4) circle (0.8mm);
\filldraw[aqua] ({\x*72-90}:1.4) circle (0.7mm);
}
\end{tikzpicture}
\caption{A layered graph of rank two, and its graded minors.}
\label{fig:layered-graph}
\end{figure}

\subsubsection{Convergence of metric graphs to tropical curves} \label{ss:GraphConvergence}

Let $G =(V,E)$ be a finite, connected, essential graph. Consider a family $(G, l_t)$, $t \in \R_+$, of weighted graphs and a weighted layered graph $(G,  \pi, \ell)$.    
For $t \in \R_+$, define the total length of edges in $(G, l_t)$ as
\begin{equation} \label{eq:DefLength}
	L_t  \coloneqq \sum_{e \in E} l_t(e),
\end{equation}
and denote by $L_{t,j}$, $j\in [r]$, the total length of edges in $\pi_j$:
\begin{equation} \label{eq:DefGrLength}
 L_{t, j} \coloneqq   \sum_{e \in \pi_j} l_t(e).
\end{equation}
We also set
\begin{equation} \label{eq:DefGrLength}
 L_{t, \fin} \coloneqq  L_{t,r+1} \coloneqq  1.
\end{equation}

We say that the family $(G, l_t)$, $t \in \R_+$, {\em converges} to $(G, \pi, \ell)$ for $t \to \infty$, if the following conditions hold:
\begin{enumerate}[label=(\Alph*)]
\item \label{item:GraphConvergenceA}  for all $j\in [r]$, we have $\rquot{L_{t,j}}{L_{t,j+1}} \to \infty$ as $t \to\infty$,
\item \label{item:GraphConvergenceB} for each edge $e\in \pi_j$, $j\in[r]$, we have $\rquot{l_t(e)}{L_{t,j}} \to \ell(e)$ as $t \to \infty$.
\item \label{item:GraphConvergenceC}  for each edge $e \in \pi_\fin$, we have $l(e) \to \ell(e)$ as $t \to \infty$.
\end{enumerate}
Informally speaking, the weighted layered graph $(G, \pi, \ell)$ captures the asymptotic behavior of the edge lengths $l_t(e)$, $e \in E$. Note that (i) implies that $L_{t, j} \to \infty$ as $t \to \infty$ for all $j \in [r]$. Moreover, as $t \to \infty$, we have $L_{t,1} \gg L_{2,t} \gg \dots \gg L_{t,r} \gg 1$. Thus, (ii) implies that the lengths of all edges in the infinite layers $\pi_j$, $j \in [r]$, go to infinity, however, they grow with different rates $L_{t,j}$, $j \in [r]$.  Moreover, by (ii) and (iii), the precise asymptotic behavior of the edge lengths $l_t(e)$, $e\in E$, is described by the edge lengths $\ell(e)$, $e \in E$, of the weighted layered graph $(G,\pi, \ell)$. 

We say that a family $(\mgr_t)_{t \in\R_+}$ of metric graphs {\em converges} to a tropical curve $\curve$, if there exist models $(G, l_t)$, $t \in \R_+$, of $\mgr_t$, with $G$ essential, and a model $(G, \pi, \ell)$ of $\curve$ such that $(G, l_t)$, $t \in \R_+$, converges to $(G, \pi, \ell)$ in the above sense.

\subsubsection{Homology and polarization of tropical curves} \label{ss:tropical_polarization}
Notations as above, let $\curve$ be a tropical curve of rank $r$.  We define the {\em first homology} of $\curve$ with real and integer coefficients as $H_1(\curve, \R)$ and $H_1(\curve, \Z)$, respectively. As for metric graphs, upon choosing a model $(G, \pi, \ell)$ of $\curve$, we can further identify the homology with the space of flows on $G$. The {\em Jacobian} of $\curve$ denoted $\Jac(\curve)$ is the quotient
\[
\Jac(\curve) \coloneqq \rquot{  H_1(\curve, \R)}{H_1(\curve, \Z)}.
\]
For each $j \in \{1, \dots, r,  \fin \}$, the inclusion $\curve^j \subset \curve$ gives a canonical embedding
\[
H_1(\curve^j, \R) \hookrightarrow H_1(\curve,  \R).
\]
Thus, we obtain a non-increasing filtration
\[
H_1(\curve, \R) = H_1(\curve^1, \R) \supseteq H_1(\curve^2, \R) \supseteq \dots \supseteq H_1(\curve^r, \R) \supseteq H_1(\curve^\fin, \R) \supseteq \{0\}.
\]
The contraction map from $\curve^j$ to the graded minor $\Gamma^j$ induces a projection map $\proj_j \colon H_1(\curve^j, \R) \to H_1(\Gamma^j, \Z)$. It is surjective and has kernel given by $\ker(\proj_j) = H_1(\curve^{j+1}, \R) \subseteq H_1(\curve^j, \R)$.

\smallskip

In the following, we abbreviate $H= H_1(\curve,  \R)$ and $\L=H_1(\curve,  \Z)$.
\begin{defi}[Tropical polarization]
Consider $\R^{r+1}$ endowed with the lexicographic order $\preceq=\lexeq$. The {\em tropical polarization} on $H$ is the bilinear form $\highinn{}{\cdot\,,\cdot} \colon H \times H \to \R^{r+1}$ given by
\[
\highinn{}{x,y} \coloneqq \left(\sum_{e\in \pi_1} \ell(e)x(e)y(e), \dots, \sum_{e\in \pi_r} \ell(e)x(e)y(e), \sum_{e\in \pi_\fin} \ell(e)x(e)y(e)\right), \qquad x,y \in H. \qedhere
\]
\end{defi}

The following proposition summarizes the basic properties of the tropical polarization.

\begin{prop}
The tropical polarization $\highinn{}{\cdot\,,\cdot} \colon H \times H \to \R^{r+1}$ is an inner product on $H= H_1(\curve,  \R)$.  The associated filtration $\filter^\bullet$ is given by
\[
\filter^j = H_1(\curve^j, \R) \subset H_1(\curve, \R), \qquad j \in \{1,  \dots, r, \fin \}.
\]
\emph{By convention, we set $\filter^{r+1} \coloneqq \filter^\fin$ and $\filter^{\fin +1} \coloneqq \{0\}$.}

\noindent For $j \in  \{1, \dots, r, \fin\}$, the projection map $\proj_j \colon H_1(\mgr^j, \R) \to H_1(\Gamma^j, \R)$ induces a canonical isomorphism
\[
\grm{}{j} H = \filter^j / \filter^{j+1}  \cong H_1(\Gamma^j, \R).
\]
Under this isomorphism, the $j$th graded piece of $(H,\highinn{}{\cdot\,,\cdot})$ is isomorphic to
\[
\big ( \grm{}{j} H, \highinn{j}{\cdot \,, \cdot} \big) \cong \big ( H_1(\Gamma^j,\R),  \innone{\Gamma^j}{\cdot \, ,  \cdot} \big ).
\]
In particular, $\grm{}{j} H = \{0\}$ if and only if $\Gamma^j$ is a tree.
\end{prop}
\begin{proof}
All properties are direct consequences of the definition of $\highinn{}{\cdot\,,\cdot}$. We omit the formal proofs.
\end{proof}
Next, we describe the lifting operators $\proj_j^\ast$,  $j \in \{1, \dots,r, \fin \}$, introduced in Lemma~\ref{lem:LiftingLemma}. For $j \in \{1, \dots, r, \fin\}$, the induced lifting operator $\proj_j^\ast \colon H_1(\Gamma^j, \R) \to H_1(\curve,  \R)$ has the following form. For $\gamma \in H_1(\Gamma^j, \R)$, the image $\eta = \proj_j^\ast(\gamma)$ is the unique flow $\eta \colon \E \to \R$ on $G$ such that
\begin{itemize}
\item [(i)] $\eta(e) = 0$ for all edges $e \in \pi_1\cup \dots \cup \pi_{j-1}$ ,
\item [(ii)] $\eta(e) = \gamma(e)$ for all edges $e \in \pi_j$ (here we identify $\pi_j$ with the edge set of $\Gamma^j$),
\item [(iii)] and for all $k \in \{j+1, \dots, r, \fin\}$, we have
\begin{equation} \label{eq:Exactness}
\sum_{e\in \pi_k} \ell(e) \, \eta(e) \alpha(e) = 0,  \qquad \forall \alpha  \in H_1(\Gamma^j, \R).
\end{equation}
\end{itemize}

\begin{remark}
Using the Hodge decomposition of the space of one-forms on a metric graph, one can prove that \eqref{eq:Exactness} is equivalent to requiring that
\[
\gamma(e) = \frac{f_k(v) - f_k(u)}{\ell(e)}, \qquad e= uv \in \pi_k,
\]
for a function $f_k \colon V \to \R$ such that $f_k(x) = f_k(y)$ for all vertices $x, y \in V$ having the same contractions in $\Gamma^j$, via the contraction map $\curve^k \to \Gamma^k$.
\end{remark}

As is easily seen, the lattice $\L = H_1(\curve, \Z)$ is admissible for the inner product $\highinn{}{\cdot\,,\cdot}$. Moreover, under the identification $\grm{}{j} H = H_1(\Gamma^j, \R)$, we have
\[
\grm{}{j}\L= H_1(\Gamma^j, \Z), \qquad j \in \{1, \dots, r, \fin \}.
\]
By Corollary~\ref{cor:covering-lattices}, the Voronoi cells $\Vor_{\curve}(\gamma)$, $\gamma \in \L$, associated to the inner product $\highinn{}{\cdot \,, \cdot}$ and the admissible lattice $\L$, provide a decomposition of $H$. Moreover, applying Theorem~\ref{thm:VoronoiCellDecomposition-lattices}, we obtain the following result.

\begin{thm}[Almost orthogonal decomposition for the Voronoi cells of tropical curves] Notations as above, we have
\[
\overline \Vor_{\curve}(0) = V_1 \aplus \dots \aplus V_r\aplus V_\fin,
\]
where $V_j = \proj_j^\ast(\Vor_{\Gamma^{\grind{j}}}(0))$, $j \in \{1, \dots, r, \fin\}$ and $\Vor_{\Gamma^{\grind{j}}}(0)$ is the Voronoi cell in $H_1(\Gamma^j, \R)$ for the lattice $H_1(\Gamma^j, \Z)$. That is, the closure of the Voronoi cell  $\Vor_\curve(0)$ for the tropical curve $\curve$ is the Minkowski sum of the Voronoi cells of its graded minors, canonically lifted to $H$.
\end{thm}

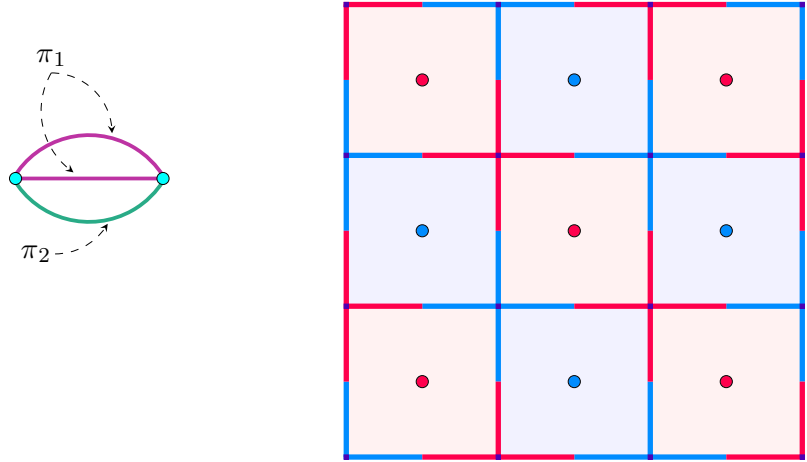
\begin{figure}[!t]
\centering
\begin{tikzpicture}
\draw[line width=.50mm, byzant] (0,2) -- (2,2);
\draw[line width=.50mm, byzant] (0,2) arc (150:30:1.15cm);
\draw[line width=.50mm, pomme] (0,2) arc (-150:-30:1.15cm);
\filldraw[aqua] (0.03,2) circle (0.7mm);
\draw[line width=0.1mm] (0.03,2) circle (0.8mm);
\filldraw[aqua] (1.97,2) circle (0.7mm);
\draw[line width=0.1mm] (1.97,2) circle (0.8mm);
\draw[-stealth, line width=.1mm, dashed] (.5, 3.4) arc (90:0:0.8cm);
\draw[-stealth, line width=.1mm, dashed] (.5, 3.4) arc (150:235:1cm);
\node at (.5, 3.6) {$\pi_1$};
\node at (.5, -1.6) { };
\node at (.3, 1) {$\pi_2$};
\draw[-stealth, line width=.1mm, dashed] (.55, 1) arc (270:330:.8cm);
\end{tikzpicture}
\hspace{2cm}
\begin{tikzpicture}
\fill[red!5!white] (0,0) -- (2,0) -- (2,2) -- (0,2) -- (0,0);
\fill[blue!5!white] (2,0) -- (4,0) -- (4,2) -- (2,2) -- (2,0);
\fill[red!5!white] (4,0) -- (6,0) -- (6,2) -- (4,2) -- (4,0);
\fill[blue!5!white] (0,2) -- (2,2) -- (2,4) -- (0,4) -- (0,2);
\fill[red!5!white] (2,2) -- (4,2) -- (4,4) -- (2,4) -- (2,2);
\fill[blue!5!white] (4,2) -- (6,2) -- (6,4) -- (4,4) -- (4,2);
\fill[red!5!white] (0,4) -- (2,4) -- (2,6) -- (0,6) -- (0,4);
\fill[blue!5!white] (2,4) -- (4,4) -- (4,6) -- (2,6) -- (2,4);
\fill[red!5!white] (4,4) -- (6,4) -- (6,6) -- (4,6) -- (4,4);
\draw[line width=0.70mm, lblue] (-0.035,0) -- (1,0);
\draw[line width=0.70mm, lblue] (3,0) -- (5,0);
\draw[line width=0.70mm, lblue] (1,2) -- (3,2);
\draw[line width=0.70mm, lblue] (5,2) -- (6,2);
\draw[line width=0.70mm, lblue] (0,4) -- (1,4);
\draw[line width=0.70mm, lblue] (3,4) -- (5,4);
\draw[line width=0.70mm, lblue] (1,6) -- (3,6);
\draw[line width=0.70mm, lblue] (5,6) -- (6.035,6);
\draw[line width=0.70mm, dred] (1,0) -- (3,0);
\draw[line width=0.70mm, dred] (5,0) -- (6.035,0);
\draw[line width=0.70mm, dred] (0,2) -- (1,2);
\draw[line width=0.70mm, dred] (3,2) -- (5,2);
\draw[line width=0.70mm, dred] (1,4) -- (3,4);
\draw[line width=0.70mm, dred] (5,4) -- (6,4);
\draw[line width=0.70mm, dred] (-0.035,6) -- (1,6);
\draw[line width=0.70mm, dred] (3,6) -- (5,6);
\draw[line width=0.70mm, lblue] (0,0) -- (0,1);
\draw[line width=0.70mm, lblue] (0,3) -- (0,5);
\draw[line width=0.70mm, lblue] (2,1) -- (2,3);
\draw[line width=0.70mm, lblue] (2,5) -- (2,6);
\draw[line width=0.70mm, lblue] (4,0) -- (4,1);
\draw[line width=0.70mm, lblue] (4,3) -- (4,5);
\draw[line width=0.70mm, lblue] (6,1) -- (6,3);
\draw[line width=0.70mm, lblue] (6,5) -- (6,6);
\draw[line width=0.70mm, dred] (0,1) -- (0,3);
\draw[line width=0.70mm, dred] (0,5) -- (0,6);
\draw[line width=0.70mm, dred] (2,0) -- (2,1);
\draw[line width=0.70mm, dred] (2,3) -- (2,5);
\draw[line width=0.70mm, dred] (4,1) -- (4,3);
\draw[line width=0.70mm, dred] (4,5) -- (4,6);
\draw[line width=0.70mm, dred] (6,0) -- (6,1);
\draw[line width=0.70mm, dred] (6,3) -- (6,5);
\filldraw[dred] (1,1) circle (0.7mm);
\draw[line width=0.1mm] (1,1) circle (0.8mm);
\filldraw[lblue] (3,1) circle (0.7mm);
\draw[line width=0.1mm] (3,1) circle (0.8mm);
\filldraw[dred] (5,1) circle (0.7mm);
\draw[line width=0.1mm] (5,1) circle (0.8mm);
\filldraw[lblue] (1,3) circle (0.7mm);
\draw[line width=0.1mm] (1,3) circle (0.8mm);
\filldraw[dred] (3,3) circle (0.7mm);
\draw[line width=0.1mm] (3,3) circle (0.8mm);
\filldraw[lblue] (5,3) circle (0.7mm);
\draw[line width=0.1mm] (5,3) circle (0.8mm);
\filldraw[dred] (1,5) circle (0.7mm);
\draw[line width=0.1mm] (1,5) circle (0.8mm);
\filldraw[lblue] (3,5) circle (0.7mm);
\draw[line width=0.1mm] (3,5) circle (0.8mm);
\filldraw[dred] (5,5) circle (0.7mm);
\draw[line width=0.1mm] (5,5) circle (0.8mm);
\foreach \x in {1,...,4}
{
\foreach \y in {1,...,4}
{
\filldraw[red!30!blue] ({2*\x-2.028},{2*\y-2.028}) -- ({2*\x-1.972},{2*\y-2.028}) -- ({2*\x-1.972},{2*\y-1.972}) -- ({2*\x-2.028},{2*\y-1.972}) -- ({2*\x-2.028},{2*\y-2.028}) -- ({2*\x-1.972},{2*\y-2.028});
}
}
\end{tikzpicture}
\caption{A tropical curve $\curve$ of rank two and its Voronoi tiling. The boundaries of Voronoi tiles are colored red or blue depending on whether they belong to the incident Voronoi cell with a red or blue center, respectively. }
\label{fig:voronoi-tiling}
\end{figure}

\subsubsection{Degeneration of the polarizations, Voronoi cells and Jacobians}
Let $G = (V,E)$ be a finite, connected, essential graph. Consider the spaces $H = H_1(G, \R)$ and $\L = H_1(G, \Z)$ of real- and integer-valued flows, respectively.

For every metric graph $\mgr$ over $G$, we fix a model $(G, l)$. Identifying $H_1(G, \R)$ with $H_1(\mgr, \R)$, the polarization on $H_1(\mgr, \R)$ induces a scalar product $\innone{\mgr}{\cdot \, , \cdot}$ on $H$. Analogously, for every tropical curve $\curve$ over $G$, we fix a model $(G, \pi, \ell)$ and obtain an inner product $\highinn{}{\cdot \, \cdot}$ on $H$. When considering a family $(\mgr_t)_{t \in \R_+}$ of metric graphs converging to a tropical curve $\curve$ over $G$, we additionally suppose that the weighted graphs $(G, l_t)$, $t \in \R_+$, converge to the weighted layered graph $(G, \pi, \ell)$.

In particular, for a family of metric graphs $\mgr_t$, $t \in \R_+$, over $G$, we obtain a family of scalar products  $\innone{\mgr_t}{\cdot \, , \cdot}$, $t \in \R_+$, on $H = H_1(G, \R)$. Our next results says that, if the metric graphs $\mgr_t$, $t \in \R_+$,  converge to a tropical curve $\curve$, then the scalar products $\innone{\mgr_t}{\cdot \, , \cdot}$, $t \in \R_+$, converge to the inner product $\highinn{}{\cdot \, , \cdot }$ in the sense of Section~\ref{ss:tameness}.

\begin{lem}
Let $(\mgr_t)_{t \in \R_+}$, be a family of metric graphs over $G$ which converges to a tropical curve $\curve$ over $G$. Consider the associated family of scalar products $\innone{\mgr_t}{\cdot \, , \cdot }$, $t \in \R_+$, on $H = H_1(G, \R)$. Then, $\innone{\mgr_t}{\cdot \, , \cdot }$, $t \in \R_+$, tamely degenerates to the tropical polarization $\highinn{}{\cdot \, , \cdot }$ with parameters $\underline L_t = (L_{t,1}, \dots, L_{t, r},  L_{t, \fin})$, $t \in \R_+$, given by
\begin{equation} \label{eq:GraphParameters}
L_{t,j} = \sum_{e \in \pi_j} l_t(e), \qquad j \in \{1, \dots, r \},
\end{equation}
and $L_{t, \fin} \equiv 1$.
\end{lem}

\begin{proof}
We need to show properties~\ref{item:tame-a}, \ref{item:tame-b} and \ref{item:tame-c} in the definition of tame degenerations. By assumption, the chosen models $(G,l_t)$ of $\mgr_t$, $t \in \R_+$, and $(G, \pi, \ell)$ of $\curve$ satisfy the conditions~\ref{item:GraphConvergenceA}, \ref{item:GraphConvergenceB} and \ref{item:GraphConvergenceC} in Section~\ref{ss:GraphConvergence}. Property \ref{item:tame-a} then follows immediately from \ref{item:GraphConvergenceA}.

We decompose $H = H_1\aplus\dots \aplus H_r\aplus H_\fin$, according to the almost orthogonal decomposition induced by the polarization of the tropical curve.

In order to prove property \ref{item:tame-b}, fix $i, j \in [r+1]$ with $i < j$ and let $\gamma \in H_i$ and $\gamma' \in H_j$. Since $\gamma=\proj_i^\ast(\eta)$ for some $\eta \in H_1(\Gamma^i, \R) \cong \grm{}{i}H$, we obtain from \eqref{eq:Exactness} that 
\[
 \sum_{e \in \pi_j} \ell(e) \, \gamma(e) \gamma'(e) =0. 
\]
Moreover, $H_j$ is contained in $\filter^j$, which is the first homology group $H_1(\curve^j,  \R) =H_1(G^j, \R)$, and hence $\gamma'(e) =0$ for all edges $e \in E\setminus (\pi_j \cup \dots \cup \pi_{r+1} \cup \pi_\fin)$ (again, by convention, $\filter^{r+1} = \filter^\fin$). Thus, taking into account \ref{item:GraphConvergenceA}, \ref{item:GraphConvergenceB} and \ref{item:GraphConvergenceC}, we conclude that, as $t \to \infty$,
\begin{align*}
\abs{\innone{\mgr_t}{\gamma , \gamma'}}&= \Big| \sum_{e \in E} l_t(e) \, \gamma(e) \gamma'(e) \Big| = \Big |  \sum_{e \in \pi_{j} \cup \dots \cup \pi_{r} \cup \pi_\fin} L_{t,j} (\ell(e) + o(1)) \gamma(e) \gamma'(e) \Big |\\
&= L_{t,j}  \Big| \sum_{e \in \pi_j} o(1) \gamma(e) \gamma'(e) + \sum_{k\in \{j+1, \dots, r, \fin\}} \sum_{e\in \pi_{k}} \frac{L_{t,k}}{L_{t,j}}(\ell(e) + o(1))\gamma(e) \gamma'(e)\Big |\\
&\le L_{t,j} \cdot  (\max_{e \in E} |\gamma(e)|)  \cdot  (\max_{e \in E} |\gamma'(e)|)  \cdot  o(1),
\end{align*}
where the $o(1)$-terms tend to zero uniformly for $\gamma \in H_i$ and $\gamma' \in H_j$.
This shows that property \ref{item:tame-b} holds.

It remains to prove property~\ref{item:tame-c}. Fix $i \in \{1, \dots, r, \fin \}$ and let $\gamma, \gamma'\in H_i$. Since $H_i \subset \filter^i$, and $\filter^i$ is the first homology group $H_1(\curve^i, \R)$, we obtain similar as above that 
\begin{align*}
\innone{\mgr_t}{\gamma , \gamma'} &= \sum_{e \in \pi_j \cup \dots \pi_r \cup \pi_\fin} l_t(e) \, \gamma(e) \gamma'(e) = L_{t,j} \sum_{k \in \{j, \dots, r, \fin\}} \frac{L_{t,k}}{L_{t,j}} \sum_{e \in \pi_k} (\ell(e) + o(1)) \, \gamma(e) \gamma'(e) \\
& = L_{t,j}\Big ( \highinn{j}{\gamma, \gamma'}   + \sum_{e \in \pi_j} o(1) \, \gamma(e) \gamma'(e) +  \sum_{k \in \{j+1, \dots, r, \fin\}} o(1) \sum_{e \in \pi_k} (\ell(e) + o(1)) \,  \gamma(e) \gamma'(e) \Big),
\end{align*}
where the $o(1)$-terms tend to zero uniformly for $\gamma, \gamma' \in H_i$. Since this implies that property~\ref{item:tame-c} holds, the proof is complete.
\end{proof}

Note that we may view the Voronoi cells $W_t = \Vor_{\innone{\mgr_t}{\cdot\, , \cdot}} (0)$ of the metric graph polarizations $\innone{\mgr_t}{\cdot \, , \cdot }$, $t \in \R_+$, and the Voronoi cell $ \Vor_{\curve}(0) =\Vor_{\highinn{}{\cdot \, , \cdot }}(0)$ of the tropical curve polarization $\highinn{}{\cdot \, , \cdot }$ as subsets of the same vector space $H= H_1(G, \R)$ with respect to the admissible lattice $\L= H_1(G, \Z)$. Applying Theorem~\ref{thm:Hausdorff_voronoi}, we obtain the following. 
\begin{cor}
Consider the lattice $\L = H_1(G, \Z)$ in $H= H_1(G, \R)$. Let $(\mgr_t)_{t \in \R_+}$, be a family of metric graphs over $G$ which converges to a tropical curve $\curve$ over $G$.  For $t \in \R_+$, let $W_t \subset H$ be the Voronoi cell of the origin $\gamma = 0$ for the polarization $\innone{\mgr_t}{\cdot \, , \cdot }$. 

Then, the Voronoi cells $W_{t}$, $t \in \R_+$, converge to the tropical Voronoi cell $\overline \Vor_{\curve}(0)$ in the Hausdorff sense.
\end{cor}

Consider the torus $\T = \rquot{H}{\L} = \rquot{H_1(G, \R)}{H_1(G, \Z)}$. For a metric graph $\mgr$ over $G$, we equip $\T$ with the distance function $\dist_\mgr \colon \T \times \T \to [0, \infty)$ induced from the scalar product $\innone{\mgr}{\cdot \,, \cdot}$.

Let $\curve$ be a tropical curve over $G$ and $(G, \pi, \ell)$ its fixed model. For $j \in \{1, \dots, r, \fin \}$, consider the subgraph $G^j = (V, E \setminus (\pi_1 \cup \dots \cup \pi_{j-1}))$ of $G$. Then we can naturally embed $H_1(G^j,  \R)$ into $H$ and $H_1(G, \Z)$ into $\L$. Setting
\[
\T^j = H_1(G^j, \R) / H_1(G^j, \Z), \qquad j \in \{1, \dots, r, \fin\},
\]
we obtain a non-increasing filtration of $\T$ by subtori:
\[
\T = \T^1 \supseteq \T^2 \supseteq \dots \supseteq \T^r \supseteq \T^\fin.
\]
Moreover, the contraction map $\curve^j \to \Gamma^j$ induces an isomorphism
\[
\rquot{\T^j}{\T^{j+1}} \cong\Jac(\Gamma^j), \qquad j \in \{1, \dots,r, \fin\}.
\]
(Here, by convention, $\T^{r+1} = \T^\fin$ and $\T^{\fin +1}  =  \{0\}$.) Applying Theorem~\ref{thm:GHConvergenceAbstract}, we obtain the following result.

\begin{cor}
Let $(\mgr_t)_{t \in \R_+}$ be a family of metric graphs which converges to $\curve$.
For $j \in \{1, \dots, r, \fin \}$, let  $L_{t,j}$ be the parameter defined in \eqref{eq:GraphParameters}. Then, as $t \to \infty$,
\[
\Big ( \T^j,\frac{1}{\sqrt{L_{t,j}}} \,\dist_{\mgr_t}  \rest{\T^j\times \T^j} \Big ) \xrightarrow{\GH} \Big ( \Jac(\Gamma^j), \dist_{\Gamma^j} \Big ).
\]
Here, we view $\T^j$ as a closed subset $\T^j \subseteq \T$, endowed with the restriction of $\dist_{\mgr_t}$ to $\T^j\times \T^j$. Moreover, $\dist_{\Gamma^j}$ is the distance function on $\Jac(\Gamma^j)$.
\end{cor}

\subsection{Degenerations of Riemann surfaces} \label{ss:DegenerationsRiemannSurfaces}
In our forthcoming work \cite{AN-AG-hybrid}, the results of Section~\ref{sec:metric_degeneration_general_tori} are applied to study the metric degeneration of Jacobians associated to degenerating Riemann surfaces. In this section, we informally discuss our results in a simplified setting.

\subsubsection{Riemann surfaces, Jacobians and polarization}
A Riemann surface $S$  is a compact one-dimensional complex manifold. We denote by $H_1(S, \R)$ and $H_1(S, \Z)$ the first homology of $S$ with coefficients in $\R$ and $\Z$, respectively. We have $\dim_\R (H_1(S, \R)) = 2g$, where $g$ is the genus of $S$, and $H_1(S, \Z)$ is a lattice of full rank in $H_1(S, \R)$. 

Denote by $\Omega$ and $\overline \Omega$ the spaces of holomorphic and anti-holomorphic one-forms on $S$, respectively, so that for $\omega \in \Omega$, the complex conjugate form $\overline \omega$ lives in $\overline \Omega$.  These are complex vector spaces of dimension $\dim_\C (\Omega)=\dim_\C (\overline \Omega)  = g$, and so, as real vector spaces, they have dimension $\dim_\R (\Omega)=\dim_\R (\overline \Omega)  = 2g$.

There exists a perfect pairing between $\Omega$ and $\overline \Omega$ defined by
\[
(\omega, \eta) \in \Omega \times \overline \Omega \quad \mapsto \quad \frac{i}{2} \int_S \omega \wedge \eta,
\]
which identifies $\overline\Omega$ with the complex dual $\Omega^*\coloneqq \Hom_\C(\Omega, \C)$.  On the other hand, as a real vector space, $\Omega^*$ can be identified with $H_1(S, \R)$ via the integration map
\[
\gamma \in H_1(S, \R) \quad \mapsto \quad  \Big (\omega \in \Omega \mapsto \int_\gamma \omega \Big),
\]
which induces  an $\R$-linear isomorphism between $H_1(S, \R)$ and $\Omega^*=\Hom_\C(\Omega, \C)$. 

The combination of these two identifications leads to a natural $\R$-linear isomorphism 
\[\Psi \colon H_1(S, \R) \to \overline \Omega,\]
more precisely defined as follows. For $\gamma \in H_1(S, \R)$, the image $\Psi(\gamma)$ is the unique anti-holomorphic one-form on $S$ such that
\[
\int_\gamma \omega = \frac{i}{2} \int_S \omega \wedge \Psi(\gamma), \qquad \forall \omega \in \Omega.
\]

The complex vector space $\overline \Omega$ has a natural pairing 
\[
(\alpha, \beta) \mapsto - \frac{i}{2} \int_S \alpha \wedge \overline \beta,
\]
which is a Hermitian sesquilinear form on $\overline \Omega$. Its real part
\[
 \innone{}{\alpha, \beta}  \coloneqq \Re \Big ( - \frac{i}{2} \int_S \alpha \wedge \overline \beta\Big ), \qquad \alpha, \beta \in \overline \Omega,
\]
is a scalar product on $\overline \Omega$, considered as a real vector space.

Altogether, we obtain a scalar product $\innone{S}{\cdot \, , \cdot} \colon H_1(S, \R) \times H_1(S, \R) \to \R $ given by
\[
\innone{S}{\gamma,\eta} = \innone{}{\Psi(\gamma), \Psi(\eta)}  , \qquad \gamma, \eta \in H_1(S, \R).
\]
We call it {\em the polarization} on $H_1(S, \R)$.

We define the {\em Jacobian} of the Riemann surface $S$ as the real $2g$-dimensional torus
\[
\Jac(S) \coloneqq \rquot{H_1(S, \R)}{H_1(S, \Z)}.
\]
The flat Riemannian metric induced by $\innone{S}{\cdot \, , \cdot}$ on $\Jac(S)$ is denoted by $\varphi_S$. Moreover, we let $\dist_S \colon \Jac(S) \times \Jac(S) \to [0, +\infty)$ be the associated distance function and $\vol(\Jac(S))$ the volume of $\Jac(S)$ for the Riemannian metric $\varphi_S$.

Note that since $\overline \Omega$ comes with an additional structure of a complex vector space, $\Psi$ can be used to give $H_1(S, \R)$ the structure of a complex vector space as well, leading to a complex structure on the Jacobian. However, this will not be important for us in what follows.

\subsubsection{Degenerating Riemann surfaces} 
\begin{figure}[ht!]
\centering
    \includegraphics[scale =0.45]{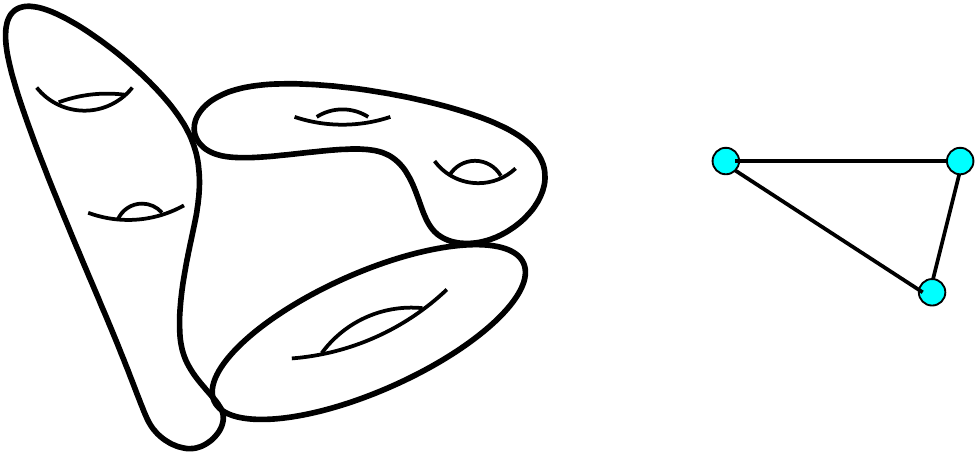}
\caption{A singular Riemann surface with three components. The corresponding graph is the cycle with three vertices.}
\label{fig:stable-surface}
\end{figure}

Let $G = (V,E)$ be a finite, connected graph. For each vertex $v \in V$, fix a smooth Riemann surface $C_v$. For each edge $e = uv$ in $G$, we fix two attachment points $p^e_u$ and $p^e_v$ on the Riemann surfaces $C_u$ and $C_v$, respectively (and additionally require that $p^e_v \neq  p^{e'}_v$ for $e \neq e'$). Identifying for each edge $e = uv$ the two points $p^ e_u$ and $p^e_v$ into a single point $p^e$, we obtain a singular Riemann surface $S$ with nodal singularities, which carries a natural analytic structure, see Figure~\ref{fig:stable-surface}.

\smallskip

What follows is a construction of a family of smooth Riemann surfaces $(\rsf_\varepsilon)_{0<\varepsilon<1}$ which degenerate to $S$ as $\varepsilon \to 0$. This is a simplified version of the so-called plumbing construction (see, e.g., \cite{HK14, Went91}).

Around each attachment point $p^e_v \in C_v$, we fix a local holomorphic coordinate $z^e_v$ and consider the unit disc 
\[
D^e_v = \left \{ z^e_v \st |z^e_v| \le 1 \right \} \subset S.
\]
 The origin $z^e_v = 0$ corresponds to the attachment point $p^e_v$. For each $0 < \varepsilon <1$, we now remove a disc of radius $\sqrt{\varepsilon}$ from $C_v$. More formally, we consider 
\[
 D^e_{v, \varepsilon} = \left\{ z^e_v \st \sqrt{\varepsilon} \le |z^e_v| \le 1 \right\}\subset S
\]
 with its inner boundary circle denoted by
 \[
 \mathbb{S}^e_{v, \varepsilon} \coloneqq \left\{ z^e_v \st |z^e_v| =  \sqrt{\varepsilon}\right\}.
 \] 
 For each edge $e = uv$, we then glue the two boundary cycles $\mathbb{S}^e_{v, \varepsilon}$ and $\mathbb{S}^e_{u, \varepsilon}$ together to a circle $\mathbb{S}^e_{\varepsilon}$. More precisely,  we identify each $z^e_u  \in\mathbb{S}^e_{v, \varepsilon}$ with the point 
\[
z^e_v = \frac{\varepsilon}{z^e_u}  \in \mathbb{S}^e_{v, \varepsilon}.
\]

By this procedure, we obtain a smooth surface $\rsf_\varepsilon$ which moreover carries a natural complex structure (that we do not explicit here), leading to a family of smooth Riemann surfaces $\rsf_\varepsilon$ parametrized by $\varepsilon \in (0,1)$. Moreover, for $\varepsilon \to 0$, the circles $\mathbb{S}^e_\varepsilon$, $e\in E$, shrink and $\rsf_\varepsilon$ degenerates to the singular Riemann surface $S$.
\begin{figure}[!h]
\centering
   \includegraphics[scale =0.45]{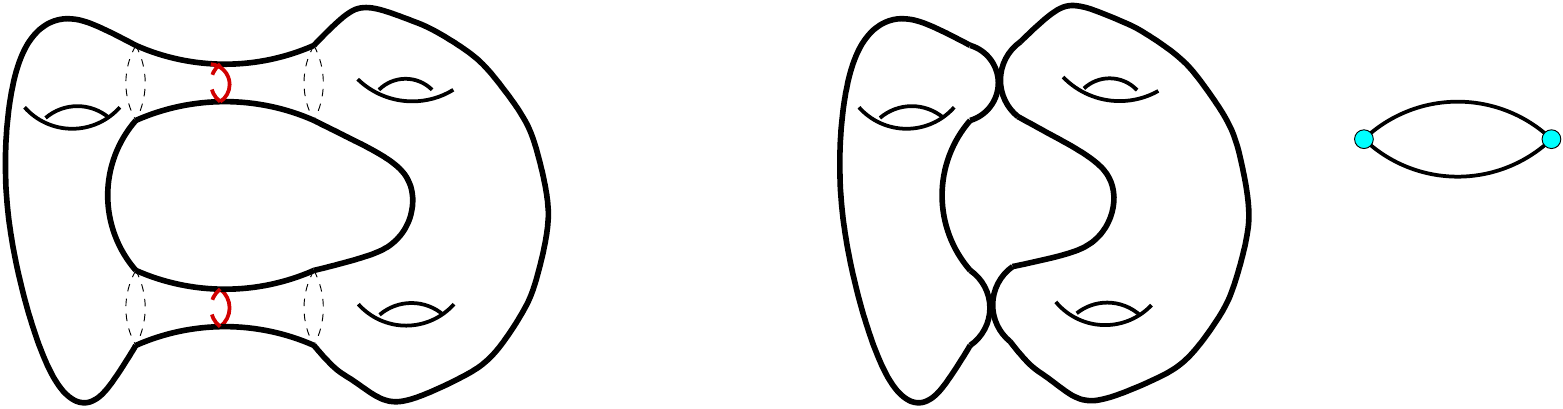}
\caption{A singular Riemann surface with two components and its associated graph on the right. A smooth Riemann surface $\rsf_\varepsilon$ in the plumbing family is depicted in the left. The cycles $\mathbb{S}^e_{\varepsilon}$ are shown in red. }
\label{fig:plumbing-family}
\end{figure}
\smallskip

\subsubsection{Jacobians of degenerating Riemann surfaces} 

For $0 < \varepsilon<1$, consider the Jacobian $\Jac(\rsf_\varepsilon) = \rquot{H_1(\rsf_\varepsilon,\R)}{H_1(\rsf_\varepsilon,\Z)}$ endowed with its Riemannian metric $\varphi_{\rsf_\varepsilon}$ and the distance function $\dist_{\rsf_\varepsilon}$. We are interested in the metric behavior of the flat tori $\Jac(\rsf_\varepsilon)$ as $\varepsilon$ tends to zero.

\smallskip

\begin{figure}[!h]
\centering
   \includegraphics[scale =0.45]{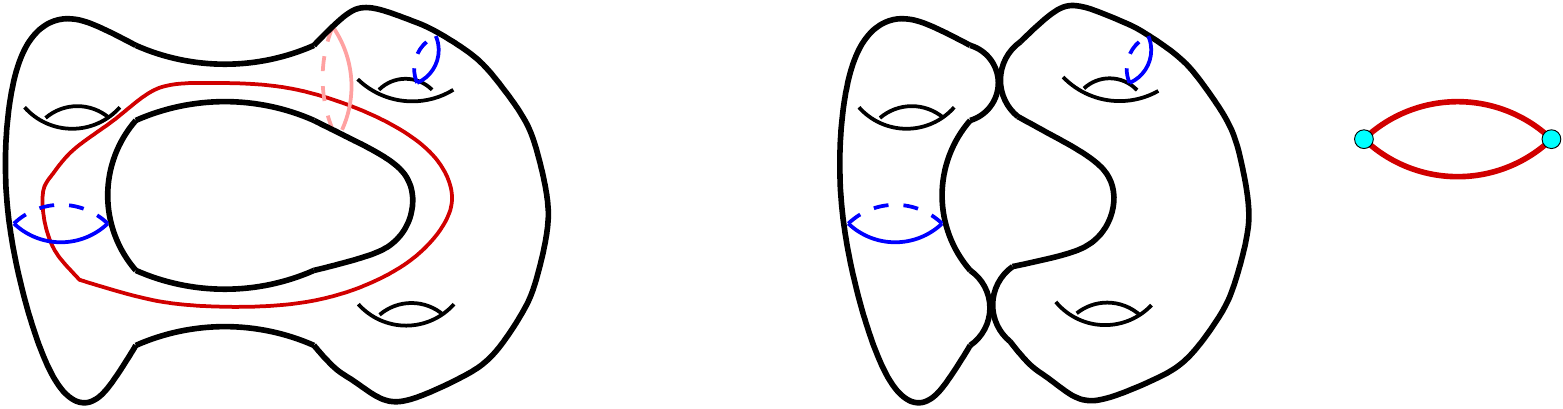}
\caption{The three type of cycles appearing in $\rsf_\varepsilon$ are shown in three colors. The cycle in pink is vanishing in the limit, it is an element of $\filter^3_\varepsilon$. The two cycles in blue live on the two components of $S$, they belong to $\filter^2_\varepsilon \setminus \filter^3_\varepsilon$. The red cycle lives above the red cycle in the graph, it belongs to $\filter^1_\varepsilon \setminus \filter^2_\varepsilon$.}
\label{fig:plumbing-cycles}
\end{figure}

Before formulating our results, we need some considerations on the homology $H_1(\rsf_\varepsilon,\R)$. Using the shape of the Riemann surface $\rsf_\varepsilon$, we can define a natural non-increasing filtration
\begin{equation}\label{eq:FiltrationSurface}
\filter^\bullet_{\varepsilon}\colon \quad  H_1(\rsf_\varepsilon,\R)  = \filter^1_\varepsilon  \supseteq \filter^2_\varepsilon \supseteq \filter^3_\varepsilon \supset \{0\}.
\end{equation}
Roughly speaking, $\filter^3_\varepsilon$ corresponds to the cycles originating from the attachment cycles $\mathbb{S}^e_\varepsilon$, $\filter^2_\varepsilon$ corresponds additionally to the cycles originating from the Riemann surfaces $C_v$, $v \in V$, and finally, $\filter^1_\varepsilon$ corresponds additionally to the cycles originating from the graph $G$ (see Figure~\ref{fig:plumbing-cycles}).  More formally, we proceed as follows.

\begin{figure}[!h]
\centering
\scalebox{.45}{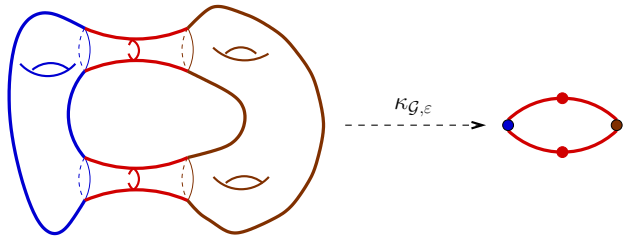}
\caption{The contraction map from $\rsf_{\varepsilon}$ to the metric graph $\mgr$. The middle cycles $\S_{e, \varepsilon}$ are mapped to the middle points of the edges, shown in red.}
\label{fig:plumbing-contraction}
\end{figure}

As is clear from the shape of the singular Riemann surface $S$, there exists a natural contraction map $\kappa_{\smallcc, \varepsilon} \colon \rsf_\varepsilon \to S$ (see Figure~\ref{fig:plumbing-family}). Namely, we contract the circles $\mathbb{S}^e_\varepsilon$ in $\rsf_\varepsilon$ to the singular points $p_e$, $e \in E$, and ``stretch" the remaining parts of $\rsf_\varepsilon$ to cover the rest of $S$. On the other hand, we can also define a natural contraction map $\kappa_{\mgr, \varepsilon} \colon \rsf_\varepsilon \to \mgr$  (see Figure~\ref{fig:plumbing-contraction}), where $\mgr$ is the metric graph obtained from $G$ and the length function $l \colon E \to (0, \infty)$ with $l(e) = 1 $. In order to define $\kappa_{\mgr, \varepsilon}$, we contract the circle $\mathbb{S}^e_\varepsilon$ of an edge $e =uv$ to the middle point of the corresponding interval edge in $\mgr$, contract the annuli $D^e_{v, \varepsilon}$ and $D^e_{u, \varepsilon}$ to cover the remaining parts of intervals, and map the remaining parts of $\rsf_\varepsilon$ to the vertices of $\mgr$.

The above contraction maps induce maps on the level of the homology, which we denote by $\proj_{\smallcc, \varepsilon} \colon H_1(\rsf_\varepsilon, \R) \to H_1(S, \R)$ and $\proj_{\mgr, \varepsilon} \colon H_1(\rsf_\varepsilon, \R) \to H_1(S,\R)$. We then define 
\[
\filter^3_\varepsilon = \ker(\proj_{\smallcc, \varepsilon}).
\]
and call it the {\em space of vanishing cycles}. Concretely, each circle $\mathbb{S}^e_{\varepsilon} \subset \rsf_\varepsilon$, $e \in E$, in $\rsf_\varepsilon$ gives rise to an element $a_{e,\varepsilon} \in H_1(\rsf_\varepsilon, \R)$. We then have
\begin{equation} \label{eq:SurfaceF3}
\filter^3_\varepsilon = \ker(\proj_{\smallcc, \varepsilon}) =\operatorname{span} \big( \{a_{e, \varepsilon}| e \in E \} \big).
\end{equation}
That is, $\filter^3_\varepsilon$ corresponds precisely to the glueing cycles $\mathbb{S}_{e, \varepsilon}$, $e \in E$, in the plumbing construction.  We stress that it might happen that the cycles $a_{e, \varepsilon}$ are not linearly independent, that is, $\sum_{e \in E} \lambda_e a_{e,\varepsilon} = 0$ in $H_1(\rsf_\varepsilon, \R)$ for some non-trivial coefficients $\lambda_e \in \R$,  $e \in E$. It turns out that $\filter^3_\varepsilon$ has dimension equal to $\graphgenus  \coloneqq \dim(H_1(\mgr, \R))$. The latter is called the genus of $\mgr$ and, by Euler's formula, can be computed as 
\[
\graphgenus = \dim(H_1(\mgr, \R)) =  |E| - |V|+1. 
\]

Note that for each smooth Riemann surface $C_v$, $v \in V$, we can naturally embed $H_1(C_v,\R)$ into $H_1(\rsf_\varepsilon,\R)$. More precisely, for each element $a_v \in H_1(C_v, \R)$, we choose a representing cycle which does not intersect the discs $D^e_v$, $e \sim v$, and view it as a cycle inside $\rsf_\varepsilon$. This allows to associate to each $a_v \in H_1(C_v, \R)$ an element $a_{v, \varepsilon} \in H_1(\rsf_\varepsilon,\R)$. Note that this choice is not canonical, since for different representing cycles, the elements $a_{v, \varepsilon}$ can differ by some combination of the $a_{e,\varepsilon}$s. Nevertheless, we obtain a non-canonical embedding $H_1(C_v, \R) \hookrightarrow  H_1(\rsf_\varepsilon,\R)$. We then define the subspace $\filter^2_\varepsilon \subset H_1(\rsf_\varepsilon, \R)$ in the filtration~\eqref{eq:FiltrationSurface} as
\begin{equation} \label{eq:SurfaceF2}
\filter^2_\varepsilon =\operatorname{span}\Big( \left\{a_{e, \varepsilon}\st e \in E \right\} \cup \bigcup_{v \in V} \left\{a_{v, \varepsilon}\st a_v \in H_1(C_v, \R) \right\} \Big).
\end{equation}
We have that $\dim(\filter^2_\varepsilon) = \graphgenus + \sum_{v \in V} 2g(C_v)$, where $g(C_v)$ denotes the genus of $C_v$. Moreover, $\filter^2_\varepsilon =\ker(\proj_{\mgr, \varepsilon}) $, where $\proj_{\mgr, \varepsilon} \colon H_1(\rsf_\varepsilon, \R) \to H_1(S,\R)$ is the projection induced by the contraction map $\kappa_{\mgr, \varepsilon} \colon \rsf_\varepsilon \to \mgr$.

\smallskip
Finally, we define the last set in the filtration~\eqref{eq:FiltrationSurface} by $\filter^1_\varepsilon =H_1(\rsf_\varepsilon,\R)$. One can give a precise description of $\filter^1_\varepsilon \setminus\filter^2_\varepsilon$. Namely, for a cycle $a_\mgr \in H_1(\mgr, \R)$, we can choose an element $a_{\mgr, \varepsilon}$ with $\proj_{\mgr, \varepsilon}(a_{\mgr, \varepsilon}) = a_\mgr$ (by ``following the graph cycle in the Riemann surface $\rsf_\varepsilon$"). This choice is only unique up to elements of $\filter^2_\varepsilon = \ker(\proj_{\mgr, \varepsilon})$. Thus, we obtain a non-canonical embedding $H_1(\mgr, \R) \hookrightarrow H_1(\rsf_\varepsilon, \R)$. Altogether, we have
\[
\filter^1_\varepsilon = H_1(\rsf_\varepsilon, \R) =\operatorname{span} \Big( \left\{a_{e, \varepsilon}\st e \in E \right\} \cup \bigcup_{v \in V} \left\{a_{v, \varepsilon}\st a_v \in H_1(C_v, \R) \right\} \cup \left\{ a_{\mgr, \varepsilon} \st  a_\mgr \in H_1(\mgr, \R) \right\}\Big) .
\]
In particular, the Riemann surface $\rsf_\varepsilon$ has genus
\[
g = \frac{\dim(H_1(\rsf_\varepsilon, \R))}{2} = \frac{2 \graphgenus + \sum_{v \in V} 2 g(C_v)}{2} = \graphgenus + \sum_{v \in V} g(C_v).
\]

\smallskip
In order to formulate our results on the degeneration of $\Jac(\rsf_\varepsilon)$, we need the following ingredients.  Setting $\L^j_\varepsilon \coloneqq \filter^j_\varepsilon \cap H_1(\rsf_\varepsilon, \Z)$, $j \in \{1,2,3\}$, we obtain a lattice of full rank in $\filter^j_\varepsilon$. Consider the associated tori
\[
 \T^j_\varepsilon \coloneqq \rquot{\filter^j_\varepsilon}{\L^j_\varepsilon},  \qquad j \in \{1,2,3\}.
\]
$\T^2_\varepsilon$ and $\T^3_\varepsilon$ have dimension  $\graphgenus + 2 \sum_{v \in V} g(C_v)$ and $\graphgenus$, respectively.
Altogether, we obtain a non-increasing filtration of $\Jac(\rsf_\varepsilon)$ by subtori:
\[
\Jac(\rsf_\varepsilon) = \T^1_\varepsilon \supseteq \T^2_\varepsilon \supseteq \T^3_\varepsilon \supset \{0\}.
\]
For each vertex $v \in V$, let $\Jac(C_v)$ be the Jacobian of the smooth Riemann surface $C_v$ and $\dist_{C_v} \colon \Jac(C_v) \times \Jac(C_v) \to [0, \infty)$ the distance function induced by the polarization $\innone{C_v}{\cdot \, , \cdot}$ on $H_1(C_v, \R)$.

 Let $\Jac(\mgr) = \rquot{H_1(\mgr, \R)}{H_1(\mgr, \Z)}$ be the Jacobian of the metric graph $\mgr$ and $\dist_\mgr \colon \Jac(\mgr) \times \Jac(\mgr) \to [0, \infty)$ the distance function induced by the polarization $\innone{\mgr}{\cdot \,, \cdot}$ on $H_1(\mgr, \R)$.

Consider the dual spaces $H_1(\mgr, \R)^\ast = \Hom_\R(H_1(\mgr, \R), \R)$ and $H_1(\mgr, \Z)^\ast = \Hom_\Z(H_1(\mgr, \Z), \Z)$. The polarization $\innone{\mgr}{\cdot \,, \cdot}$ on $H_1(\mgr, \R)$ defines an isomorphism $\Psi \colon H_1(\mgr, \R) \to H_1(\mgr,  \R)^\ast$ by setting
\[
(\Psi\gamma)(\eta) = \innone{\mgr}{\eta, \gamma}, \qquad \gamma, \eta\in H_1(\mgr, \R).
\]
Thus, we obtain a scalar product $\innone{\mgr, \ast}{\cdot \,, \cdot} \colon H_1(\mgr,  \R)^\ast \times H_1(\mgr,  \R)^\ast \to \R$ on $H_1(\mgr, \R)^\ast$ given by
\[
\innone{\mgr, \ast}{x,y}\coloneqq \innone{\mgr}{\Psi^{-1}(x), \Psi^{-1}(y)}, \qquad x,y \in H_1(\mgr,  \R)^\ast.
\]
Consider the dual torus
\[
\Jac^\ast(\mgr) \coloneqq \rquot{H_1(\mgr, \R)^\ast}{H_1(\mgr, \Z)^\ast}.
\]
We endow $\Jac^\ast(\mgr) $ with the distance function $\dist_{\mgr}^\ast \colon \Jac^\ast(\mgr)  \times \Jac^\ast(\mgr)  \to [0, \infty)$ induced from the scalar product $\innone{\mgr, \ast}{\cdot \,, \cdot}$.

\begin{thm} \label{thm:RiemannSurfaceDegeneration}
The volume of the torus $\Jac(\rsf_\varepsilon)$ converges to
\[
\lim_{\varepsilon \to 0} \vol(\Jac(\rsf_\varepsilon))  = \prod_{v \in V} \vol(\Jac(C_v)).
\]
Moreover, as $\varepsilon \to 0$,
\begin{align*}
\Big (\Jac(\rsf_\varepsilon), \,  \frac{1}{\sqrt{|\log(\varepsilon)|}}\, d_{\rsf_\varepsilon} \Big) \qquad &\xrightarrow{\GH} \qquad \Big ( \Jac(\mgr), \dist_{\mgr} \Big ) \\[2mm]
\Big ( \T^2_\varepsilon, \,  d_{\rsf_\varepsilon} \rest{\T^2_\varepsilon \times \T^2_\varepsilon} \Big)\qquad &\xrightarrow{\GH} \qquad \Big (\bigoplus_{v \in V} \Jac(C_v), \bigoplus_{v \in V} \dist_{C_v} \Big) \\[2mm]
 \Big (\T^3_\varepsilon, \, \sqrt{|\log(\varepsilon)|}\, d_{\rsf_\varepsilon}\rest{\T^3_\varepsilon \times \T^3_\varepsilon} \Big) \qquad &\xrightarrow{\GH} \qquad  \Big (\Jac^\ast(\mgr), \dist_{\mgr}^\ast \Big)
\end{align*}
in the Gromov--Hausdorff sense. Here, we view $\T^j_\varepsilon$ as a closed subset of $\Jac(\rsf_\varepsilon)$, endowed with the restriction of $\dist_{\rsf_\varepsilon}$ to $ \T^j_\varepsilon \times  \T^j_\varepsilon$.
\end{thm}
The notation $ (\bigoplus_{v \in V} \Jac(C_v), \bigoplus_{v \in V} \dist_{C_v} \Big)$ here refers to the product metric.

\begin{remark}
Note that, as $\varepsilon \to 0$, the degeneration parameter $\sqrt{|\log(\varepsilon)|}$ tends to $+ \infty$. In particular, the diameter of the Jacobians $\Jac(\rsf_t)$ goes to $+ \infty$ under degeneration. On the other hand, the volume remains bounded. This is explained by the fact that the Jacobian becomes ``infinitely long" in some directions, and at the same, ``infinitely thin" in other ones.
\end{remark}

\subsubsection{Tame degenerations of polarization on Riemann surfaces} We informally describe here the proof of Theorem~\ref{thm:RiemannSurfaceDegeneration} given in \cite{AN-AG-hybrid}. Fix some $\varepsilon_0 >0$ and consider the vector space $H = H_1(\rsf_{\varepsilon_0}, \R)$. For every other $\varepsilon \in (0,1)$, we can define a natural diffeomorphism $f_\varepsilon \colon \rsf_{\varepsilon_0} \to  \rsf_\varepsilon$ by mapping the circles $\mathbb{S}^e_{\varepsilon_0}$ to the circles $\mathbb{S}^e_{\varepsilon}$ and ``stretching" the remaining parts of $\rsf_{\varepsilon_0}$. This induces an isomorphism between $H_1(\rsf_\varepsilon, \R)$ and $H = H_1(\rsf_{\varepsilon_0}, \R)$. Thus the scalar products $\innone{\rsf_\varepsilon}{\cdot \,,  \cdot}$ on $H_1(\rsf_\varepsilon,  \R)$ define a family of scalar products $\innone{\varepsilon}{\cdot \, , \cdot}$,  $\varepsilon \in (0,1)$, on $H$.

\smallskip
Consider the bilinear form $\highinn{}{\cdot \,  , \cdot} = (\highinn{1}{\cdot \,  , \cdot}, \highinn{2}{\cdot \,  , \cdot}, \highinn{3}{\cdot \,  , \cdot}) \colon H \times H \to \R^{3}$ defined as follows. The first component $\highinn{1}{\cdot \,  , \cdot}$ is defined by
\[
\highinn{1}{\gamma, \eta} \coloneqq \innone{\mgr}{\proj_{\mgr, \varepsilon_0} (\gamma), \proj_{\mgr, \varepsilon_0} (\eta)}, \qquad \gamma, \eta \in H.
\]
Note that $\oplus_v H_1(C_v,  \R)$ embeds into $H_1(S,\R)$ and the image of $\filter^2_{\varepsilon_0}$ under the projection map $\proj_{\smallcc, \varepsilon_0}$ lies in this subspace (see~\eqref{eq:SurfaceF2}). Using the polarizations $\innone{C_v}{\cdot \, , \cdot}$ of the Riemann surfaces $C_v$, $v \in V$, we define the second component $\highinn{2}{\cdot \,  , \cdot}$ of the bilinear form $\highinn{}{\cdot \,  , \cdot}$ by
\[
\highinn{2}{\gamma, \eta} \coloneqq \begin{cases} \sum_{v \in V} \innone{C_v}{\proj_{\smallcc, \varepsilon_0}(\gamma) \, , \, \proj_{\smallcc, \varepsilon_0}(\eta) },  & \gamma, \eta \in \filter^{2}_{\varepsilon_0}, \\
 0, & \gamma\notin \filter^{2}_{\varepsilon_0} \text{ or }\eta\notin \filter^{2}_{\varepsilon_0}.
\end{cases}
\]
In order to define the third component $\highinn{2}{\cdot \,  , \cdot}$ of the bilinear form $\highinn{}{\cdot \,  , \cdot}$, we introduce the following projection map $\proj_{\ast, \varepsilon_0}\colon \filter^3_{\varepsilon_0} \to H_1(\mgr, \R)^\ast$. Given $\gamma \in \filter^3_{\varepsilon_0}$, the linear functional $\proj_{\ast, \varepsilon_0} (\gamma) \in H_1(\mgr, \R)^\ast$ is defined by
\[
\big( \proj_{\ast, \varepsilon_0} (\gamma) \big) (a_\mgr)= \innone{\operatorname{int}}{\gamma, a_{\mgr, {\varepsilon_0}}}, \qquad a_\mgr \in H_1(\mgr, \R),
\]
where $\innone{\operatorname{int}}{\cdot \, , \cdot }$ is the intersection pairing on $H_1(\rsf_{\varepsilon_0}, \R)$. As above, for $a_\mgr \in H_1(\mgr, \R)$ we have chosen a lift $a_{\mgr, \varepsilon_0}$, that is, an element $a_{\mgr, \varepsilon_0}  \in H_1(\rsf_{\varepsilon_0})$ satisfying $\proj_{\mgr, \varepsilon}(a_{\mgr, \varepsilon_0}) = a_\mgr$. Note that two different such lifts differ by an element of $\filter^2_{\varepsilon_0} = \ker(\proj_{\mgr, \varepsilon_0})$. Moreover, we have $\innone{\operatorname{int}}{\gamma, \eta } =  0$ for all $\gamma \in \filter^3_{\varepsilon_0}$ and $\eta \in \filter^2_{\varepsilon_0}$ (consider \eqref{eq:SurfaceF3} and \eqref{eq:SurfaceF2}, taking into account the locations of the cycles on $\rsf_{\varepsilon}$ shown in Figure~\ref{fig:plumbing-cycles}). This ensures that the projection map $\proj_{\ast, \varepsilon_0}\colon \filter^3_{\varepsilon_0} \to H_1(\mgr, \R)^\ast$ is indeed well-defined.

We define the third component $\highinn{3}{\cdot \,  , \cdot}$ of $\highinn{}{\cdot \,  , \cdot}$ as
\[
\highinn{3}{\gamma, \eta} \coloneqq \begin{cases}  \innone{\mgr, \ast}{\proj_{\ast, \varepsilon_0} (\gamma) \,, \, \proj_{\ast, \varepsilon_0}  (\eta)},  & \gamma, \eta \in \filter^{3}_{\varepsilon_0}, \\
0, & \gamma\notin \filter^{3}_{\varepsilon_0} \text{ or }\eta\notin \filter^{3}_{\varepsilon_0}.
\end{cases}
\]

The following result is proven in \cite{AN-AG-hybrid}. Theorem~\ref{thm:RiemannSurfaceDegeneration} then follows by applying Theorem~\ref{thm:GHConvergenceAbstract}.

\begin{prop} \label{prop:DegeneratingSurfaceWeakTameness}
The bilinear form $\highinn{}{\cdot \,, \cdot} \colon H \times H \to \R^3$ is an inner product on $H = H_1(\rsf_{\varepsilon_0}, \R)$. Its filtration $\filter^\bullet_{\highinn{}{\cdot \,, \cdot}}$ coincides with $\filter^\bullet_{\varepsilon_0}$. Moreover, as $\varepsilon \to 0$, the scalar products $\innone{\varepsilon}{\cdot \,,  \cdot}$ degenerate $\omega$-tamely to $\highinn{}{\cdot \,,  \cdot}$ with parameters 
\[
\underline L_\varepsilon = (|\log\varepsilon|^{1/2}, 1, |\log\varepsilon|^{-1/2}).
\]
\end{prop}

\begin{remark} \label{rem:OutlookAGHybrid}
The results in \cite{AN-AG-hybrid} go beyond Proposition~\ref{prop:DegeneratingSurfaceWeakTameness} and Theorem~\ref{thm:RiemannSurfaceDegeneration} on several levels. On the one hand, we prove a more precise asymptotic for the degenerating scalar products $\innone{\varepsilon}{\cdot \,, \cdot}$, which implies the $\omega$-tameness. These results are based on our study of degenerations of solutions to Poisson equations in~\cite{AN2}, and are more elaborate to formulate here.

On the other hand, \cite{AN-AG-hybrid} deals with families of Riemann surfaces which depend on several parameters. The above degeneration parameter $\varepsilon >0$ is replaced by independent complex parameters $\varepsilon_e$, one for each edge $e \in E$.  In this case, additional multi-scale effects occur, which are then handled by using the notion of tropical curves from Section~\ref{ss:TropicalCurves}. In this more general framework, the tori $\Jac(\mgr)$ and $\Jac^\ast(\mgr)$ in Theorem~\ref{thm:RiemannSurfaceDegeneration} are replaced by several tori $\Jac(\Gamma^1), \dots, \Jac(\Gamma^r)$ and $\Jac^\ast(\Gamma^1), \dots, \Jac^\ast(\Gamma^r)$ for the graded minors $\Gamma^1, \dots, \Gamma^r$ of a tropical curve $\curve$, each corresponding to a different scale. The study of multi-parameter degenerations of Riemann surfaces is linked to the geometric properties of $\mgg{g}$, the moduli space of Riemann surfaces of genus $g$, and its Deligne--Mumford compactification $\mgbarg{g}$. In fact, the analysis of the analytic behavior of Riemann surfaces close to the boundary of the Deligne--Mumford compactification requires to consider degenerations depending on $3g-3$ complex parameters (recall that $\mgg g$ has complex dimension $3g-3$).

Finally, we mention that Jacobians of degenerating Riemann surfaces can be treated more naturally in the framework of hybrid curves. Hybrid curves are geometric objects which mix both tropical curves and Riemann surfaces. This framework has been established in \cite{AN, AN2} and allows to understand conceptually  several asymptotic properties of degenerating Riemann surfaces. In order to keep the exposition streamlined, we do not develop the hybrid point of view here and refer to our forthcoming work \cite{AN-AG-hybrid} for more details.
\end{remark}

\section{Further discussions and open problems} \label{sec:Discussion}

\subsection{Higher rank Hermitian inner products} \label{ss:HermitianInnerProducts} We briefly discuss how to adapt the set-up of the paper to Hermitian inner products. Let $H$ be a complex vector space of finite dimension and let $r$ be a positive integer. Consider a function $\highinn{}{\cdot \,, \cdot} \colon H \times H \to \C^r$ with coordinates $\highinn{j}{\cdot\,,\cdot}$, $j \in [r]$. We assume that $\highinn{}{\cdot \,, \cdot}$ is Hermitian, that is,
\[\highinn{}{x , y} = \overline{\highinn{}{y , x}}\qquad \forall x,y \in H.\]
Moreover, we assume that $\highinn{}{\cdot \,, \cdot}$ is sesquilinear, that is, it is linear in the first factor. It follows then that $\highinn{}{\cdot \,, \cdot}$ is antilinear in the second factor. 
 
Note that it follows from the conjugate-symmetric property that for any $x \in H$, $\highinn{}{x,x}$ belongs to $\R^r$, so that it makes sense to talk about positivity by saying that $\highinn{}{x \,, x} \succ 0$ for $x\neq 0$.

The form $\highinn{}{\cdot \,, \cdot}$ defines a non-increasing filtration
\[
\filter^\bullet\colon \qquad \filter^1\coloneqq H \supseteq \filter^2 \supseteq \filter^3 \supseteq \dots \supseteq \filter^r \supseteq \filter^{r+1} \coloneqq (0)
\]
on $H$ by setting
\begin{align*}
\filter^j  \coloneqq \left\{ x \in H \st \,\highinn{}{x \,, y}_i = 0 \text{ for all $i <j$ and all $y \in H$}\right\}, \qquad j \in [r].
\end{align*}

We say that $\highinn{}{\cdot\,, \cdot}$ is a (higher rank) \emph{inner product} 
if for each $j \in [r]$, the induced form $\highinn{j}{\cdot\,, \cdot}$ on $\filter^j / \filter^{j+1}$ is a Hermitian inner product. It is easy to see that this is equivalent to requiring $\highinn{j}{\gamma, \gamma} >0$ for all $\gamma \in \filter^j \setminus \filter^{j+1}$.

All the results of this paper can be extended to this setting.

\subsection{Topological characterization of admissible discrete sets} \label{ss:HybridTopology} Let $H$ be a real vector space of finite dimension endowed with an inner product $\highinn{}{\cdot\,,\cdot} \colon H\times H \to \R^r$. Denote by $\qf \colon H \to \R^r$ the corresponding quadratic form, $\qf(x) =\highinn{}{x,x}$ for all $x\in H$. 

For $r$ positive reals $\rho_1, \dots, \rho_r>0$, we define the \emph{polydisk} $B_\rho=B_{\rho_1, \dots, \rho_r}(0)$ with center $z = 0$ and radius $\rho=(\rho_1, \dots, \rho_r) \in \R_+^r$ by 
\begin{align*}
B_{\rho}(0)\coloneqq &\left\{x\in H \st 0<\qf_1(x) <\rho_1\right\} \sqcup \left\{x\in \filter^2 \st 0<\qf_2(x) <\rho_2\right\} \sqcup \dots \\
&\sqcup \left\{x\in \filter^r \st 0<\qf_r(x) <\rho_r\right\}\sqcup\{0\}. 
\end{align*} 
 More generally, we define the polydisk with center $z\in H$ and radius $\rho \in \R_+^r$ by
\[B_{\rho}(z)\coloneqq B_{\rho}(0) + z.\]

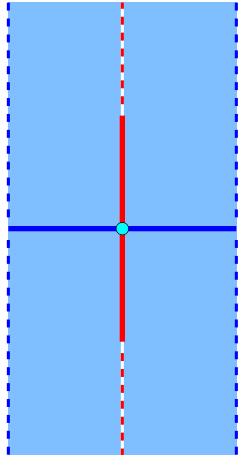
\begin{figure}[!t]
\centering
\begin{tikzpicture}
\fill[azure!50!white] (0,0) -- (3,0) -- (3,6) -- (0,6) -- (0,0);
\draw[line width=0.40mm, blue, dashed] (3,2.85) -- (3,0);
\draw[line width=0.40mm, blue, dashed] (3,3.15) -- (3,6);
\draw[line width=0.40mm, blue, dashed] (0,2.85) -- (0,0);
\draw[line width=0.40mm, blue, dashed] (0,3.15) -- (0,6);
\draw[line width=0.70mm, blue] (0,3) -- (3,3);
\draw[line width=0.70mm, red] (1.5,1.5) -- (1.5,4.5);
\draw[line width=0.40mm, white] (1.5,1.5) -- (1.5,0);
\draw[line width=0.40mm, white] (1.5,4.5) -- (1.5,6);
\draw[line width=0.40mm, red, dashed] (1.5,4.65) -- (1.5,6);
\draw[line width=0.40mm, red, dashed] (1.5,1.35) -- (1.5,0);
\draw[line width=0.1mm] (1.5,3) circle (0.8mm);
\filldraw[aqua] (1.5,3) circle (0.7mm);
\end{tikzpicture}
\caption{A polydisk in $H=\R^2$ endowed with inner product from Example~\ref{ex:HigherRankEuclidean}. It is the union of three sets, in blue $\bigl\{x\in H\st \qf_1(x) \in(0, \rho_1)\bigr\}$, in red $\bigl\{x\in \filter^2 \st \qf_2(x)\in (0, \rho_2)\bigr\}$, and the center of the polydisk in azure.}
\label{fig:polydisk}
\end{figure}

Figure~\ref{fig:polydisk} depicts a polydisk in $\R^2$ endowed with the Euclidean inner product of rank two from Example~\ref{ex:HigherRankEuclidean}.

\begin{thm} The polydisks $B_{\rho}(z)$, $\rho\in \R_+^r$ and $z\in H$, form the base of a topology on $H$. \end{thm}
\begin{proof} We decompose $H$ as the direct sum $H=H_1\oplus \dots \oplus H_r$ according to the almost orthogonal decomposition. More generally, for any such decomposition into a direct sum, we define a topology whose base of open sets are given by translations of the sets of the form 
\[B = \bigl(B_1^\times \oplus H_2 \oplus \dots \oplus H_r\bigr)\, \bigsqcup \, \bigl(B_2^\times \oplus H_3 \oplus \dots \oplus H_r\bigr) \,\bigsqcup \, \dots\, \bigsqcup\, B_r^\times \,\bigsqcup\, \{0\} \]
 where $B_1, \dots, B_r$ are open sets around the origin in $H_1, \dots, H_r$ and $B_j^\times$ denotes the punctured open set $B_j \setminus\{0\}$. Obviously, these sets and their translations $z+B$, $z\in H$, cover $H$. One verifies directly that for any point $x\in B$, for $B$ as above, there exists a set $B'$ as above such that $x+B'$ is contained in $B$. This shows that the sets $z +B$, $B$ as above and $z \in H$, form the base of a unique topology on $H$. The claim follows.
 \end{proof}
Note that this topology is a hybrid of the (pullback of the) order topology on $\R^r$ (induced by the lexicographic order $\lexeq$), and the Euclidean topologies of the graded pieces $\grm{}{j}H$, $j=1, \dots, r$. It does not have the Hausdorff property. Moreover, it is coarser than the Euclidean topology on $H$, since the identity map $H \to H$ from $H$ with Euclidean topology to $H$ with the new topology is continuous. It follows that $H$ has fewer discrete sets in this topology than discrete sets in Euclidean topology. We have the following theorem.
\begin{thm}[Topological characterization of admissible discrete sets] \label{thm:AdmissibleDiscreteTopological} Let $(H,\highinn{}{\cdot\,,\cdot})$ be an inner product space and $S \subset H$. Then, the following statements are equivalent... 
\begin{itemize} 
\item $S \subset H$ is an admissible discrete set.

\item $S$ is discrete for the above defined topology.

\end{itemize}

\end{thm}
We omit the proof of Theorem~\ref{thm:AdmissibleDiscreteTopological}.  Note that one obtains an alternative proof of Proposition~\ref{prop:discrete-discrete} by combining Theorem~\ref{thm:AdmissibleDiscreteTopological} with the above remarks.

 \subsection{Compact Hausdorff convergence} \label{ss:CompactHausdorffConvergence} 
The aim of this section is to provide an alternative description of the compact Hausdorff convergence from Section~\ref{sec:Hausdorff_convergence}.

\smallskip

Let $H$ be a real finite-dimensional vector space equipped with a norm $\refnorm{\cdot} \colon H \to [0, \infty)$. As in Section~\ref{sec:Hausdorff_convergence}, we denote the distance between a point $a \in H$ and a non-empty subset $X \subset H$  by $\dist(a, X) = \inf_{x \in X} \refnorm{a-x}$. The ball of radius $r >0$ with center $x \in H$ is denoted by $B_r(x)$. The Hausdorff distance $\hdist(X, Y)$ between two non-empty subsets $X, Y \subset H$ is given by
\[
\hdist(X,Y) = \max\Bigl\{ \sup_{x \in X} \dist(x,Y),  \, \sup_{y \in Y} \dist(y,X) \Bigr\}.
\]
Moreover, we set $\hdist(X,\varnothing) \coloneqq \hdist(\varnothing, X) = +\infty$ for non-empty $X \subset H$ and  $\hdist(\varnothing,\varnothing) \coloneqq 0$.
 A sequence $(Z_t)_{t \in \R_+}$ of non-empty subsets $Z_t\subset H$ is said to \emph{compactly converge in the Hausdorff distance} to a non-empty set $Z \subset H$ if the distance functions $\dist(y, Z_t)$, $t \in \R_+$, converge to the distance function $\dist(y, Z)$ uniformly on compact subsets of $H$ for $t \to \infty$.

We prove the following theorem.
 \begin{thm} \label{thm:CompactHausdorffConvergence} Let $Z$ and $Z_t$, $t\in \R_+$, be subsets of $H$. Then, the following are equivalent.
 \begin{enumerate}
 \item $(Z_t)_{t \in \R_+}$ compactly converges to $Z$ in the Hausdorff metric.
 \item There exists an open covering $H = \bigcup_{k\in \N} B_k$ by bounded open sets $B_k \subset H$ such that $\hdist(Z_t \cap B_k, Z\cap B_k) \to 0$ as $t \to \infty$ for every $k \in \N$. 
 \end{enumerate}
 \end{thm}
\begin{proof}
We only prove that (1) implies (2).  The proof of the other direction is more direct and we omit it.

If $Z_t$ converges compactly to $Z$, then for any bounded open set $U$ and any $\varepsilon >0$, there exists a positive $N=N(U, \varepsilon)$ such that
\begin{align} \label{eq:recrusive-key-open}
\forall \,\,t\geq N \qquad \sup_{z\in U \cap Z}\dist(z, Z_t) , \sup_{z\in  U\cap Z_t} \dist(z, Z) < \varepsilon.
\end{align}

Let $x \in H$ be a point. We show that there exists an open bounded set $B$ containing $x$ such that $\hdist(Z_t \cap B, Z\cap B) \to 0$ as $t \to \infty$. This will prove the result.

 We fix a sequence $\varepsilon_1, \varepsilon_2, \dots$ of positive reals such that $\sum_{j=1}^\infty \varepsilon_j <\infty$.

 Let $U_1$ be an open ball around $x$ such that $U_1\cap Z$ is non-empty. From the statement (1) in the theorem, we deduce that for all large $t$, $U_1 \cap Z_t$ is non-empty. Let $N_1$ be the positive number given in Equation~\eqref{eq:recrusive-key-open} for the open set $U=U_1$ and $\varepsilon =\varepsilon_1$.
 
 Proceeding by induction, we construct an increasing sequence of open sets $U_1 \subset U_2 \subset \dots$ as follows. Having constructed $U_k$, $k\in \N$, we set $N_k>N_{k-1}$ to be a positive number which verifies~\eqref{eq:recrusive-key-open} for the set $U=U_k$ and $\varepsilon = \varepsilon_k$. We define then $U_{k+1}$ by
\[U_{k+1} \coloneqq U_k \cup\bigcup_{z \in Z \cap U_k} B_{\varepsilon_k}(z) \cup \bigcup_{t\geq N_k} \bigcup_{z \in Z_t \cap U_k} B_{\epsilon_k}(z).\]
Let $B = \bigcup_{k\in \N} U_k$. By construction, $B$ is a bounded, open subset of $H$. We show that 
\begin{equation} \label{eq:hconverge}\hdist(Z_t \cap B, Z\cap B) \to 0, \,\, \textrm{as } \,\, t\to \infty, 
\end{equation}
which finishes the proof.  We need the following result.
\begin{claim} \label{claim1} Let $l>k$ be a pair of positive integers. Assume that $y \in Z\cap U_l$ or $y\in Z_t \cap U_l$ for some $t \geq N_k$. Then, there exists $z \in U_{l-1}$ such that $\dist(y, z) \leq \varepsilon_{l-1}$ and either, $z \in Z \cap U_{l-1}$, or, $z \in Z_s \cap U_{l-1}$ for some $s\geq N_k$. 
\end{claim}
\begin{proof} Suppose first that $y \in Z\cap U_l$. By the definition of $U_l$, we have either $y\in U_{l-1}$, or $y \in B_{\epsilon_{l-1}}(z)$ for $z\in Z\cap U_{l-1}$ or $z \in  B_{\epsilon_{l-1}}(z)$ for $z \in Z_s\cap U_{l-1}$ for some $s\geq N_{l-1} \geq N_k$, from which the result follows.

Suppose now that $y \in Z_t\cap U_l$. Again, we have either, $y\in U_{l-1}$, or $y \in B_{\epsilon_{l-1}}(z)$ for $z\in Z\cap U_{l-1}$, or $z \in  B_{\epsilon_{l-1}}(z)$ for $z \in Z_s\cap U_{l-1}$ for $s\geq N_{l-1} \geq N_k$. In the first case, we set $s=t$. In all cases, we have the result. 
\end{proof}

We now prove \eqref{eq:hconverge}.  We need to show that for any $\varepsilon>0$, there exists $T$ such that for all $t \geq T$, we have
\[\sup_{y\in Z\cap B} \dist(y, Z_t\cap B) , \sup_{y \in Z_t \cap B} (y, Z\cap B) \leq \varepsilon. \]

 We choose $k$ large enough such that $3\sum_{j \geq k}\varepsilon_{j} < \varepsilon$. 

Let first $y$ be a point in $Z\cap B$. There exists an integer $l \geq k+1$ such that $y \in Z \cap U_l$. Applying Claim~\ref{claim1} $(l-k)$ times, we infer the existence of a point $z$ either in $Z \cap U_k$ or in $Z_s\cap U_k$, for some $s \geq N_k$, such that $\dist(y,z) \leq \varepsilon_{k}+ \dots + \varepsilon_{l-1}$. In the second case, by the choice of $N_k$, there exists a point $z'$ in $Z \cap B_{\varepsilon_k}(z)$. We infer that in either case, there exists a point $w \in Z\cap U_{k+1}$ such that  
$\dist(y,w) \leq 2 \sum_{j \geq k}\varepsilon_{j}$.

 Since $B_{\varepsilon_{k+1}}(w)$ contains a point of $Z_t \cap U_{k+1}$ for any $t\geq N_{k+1}$,  we conclude setting $T \coloneqq N_{k+1}$ that
$\dist(y, Z_t \cap B) \leq 3\sum_{j \geq k}\varepsilon_{j} \leq \varepsilon$ for any $t \geq T$, as required.

We prove that for $t \ge T = N_{k+1}$ the inequality $\sup_{y \in Z_t \cap B} (y, Z\cap B) \leq \varepsilon$ holds as well. Let $y$ be a point in $Z_t \cap B$ for some $t\geq T$. There exists an integer $l \geq k+1$ such that $y \in Z_t\cap U_l$. Applying again Claim~\eqref{claim1} $(l-k-1)$ times, we infer the existence of a point $z$ either in $Z \cap U_{k+1}$ or in $Z_s\cap U_{k+1}$, for some $s \geq N_{k+1}$, such that $\dist(y,z) \leq \varepsilon_{k+1}+ \dots + \varepsilon_{l-1}$. In the second case, by the choice of $N_k$, the ball $B_{\varepsilon_{k+1}}(z)$ contains a point $w$ of $Z$. Moreover, by construction, $w$ belongs to $Z \cap U_{k+2} \subset Z \cap B$. This means that in either case, we have $\dist(y, Z \cap B) \leq 2 \sum_{j \geq k+1}\varepsilon_{j} \leq \varepsilon$, as required.
\end{proof}

\subsection{An example of a metric collapse} \label{ss:PathologyTorus} 
In the Section~\ref{sec:metric_degeneration_general_tori}, we obtained results on the metric degeneration of tori associated to tamely degenerating families of scalar products. In this section, we present an example of a degenerating family of scalar products which looks quite similar to a tame degeneration, but where Theorem~\ref{thm:GHConvergenceAbstract} fails.

\smallskip

Let $H = \R^2$ and consider the standard lattice $\L = \Z^2$  in $H$. We denote elements $x \in H$ by $x = (x_1, x_2)$.  For $t\in \R_+$, consider the scalar bilinear form
\[
\innone{t}{x,y} = x_1 y_1 + b_t (x_1 y_2 + x_2 y_1) + t^{-1} x_2 y_2,\qquad x,y \in H,
\]
where $b_t$ is a real parameter with $0< b_t < t^{-1/2}$ such that $\lim_{t \to \infty} t \cdot b_t^2= 1$. The Gram matrix  of $\innone{t}{\cdot\,,\cdot}$ with respect to the standard basis $e_1, e_2 \in \Z^2$ is given by
\[
M_t = \begin{pmatrix} 1 & b_t \\
b_t & t^{-1} \end{pmatrix}
\]
and has strictly positive determinant. It follows that $\innone{t}{\cdot\,, \cdot}$ is a scalar product on $H$.

The scalar products $\innone{t}{\cdot \,, \cdot}$, $t \in \R_+$, can be realized as a pullback family as follows. Consider the bilinear form
\begin{align*}
\highinn{}{\cdot\,,\cdot} &\colon H \times H \to \R^3\\
\highinn{}{x\,, y} &\coloneqq \Big ( x_1y_1, x_1 y_2 + x_2 y_1,  x_2 y_2  \Big ), \qquad x,y \in H.
\end{align*}
 The scalar product $\innone{t}{\cdot\,, \cdot}$ coincides with the pullback $\innone{\underline L_t}{\cdot\,, \cdot}$ for the vector 
\[\underline L_t= (L_{t,1}, L_{t,2}, L_{t,3}) = (1, b_t, t^{-1}).\]
   We also have that $L_{t,1} / L_{t,2} \to \infty$ and $L_{t,2} /L_{t,3} \to \infty$ for $t \to \infty$. 
  
The bilinear form $\highinn{}{\cdot\,, \cdot}$ is easily seen to be positive, that is, $\highinn{}{x,x} \gex 0$ for $x \in H \setminus \{0\}$. Indeed, for a non-zero element $x \in H$, if $x_1\neq 0$, then we have $\highinn{1}{x,x}>0$, and if $x_1=0$, then $\highinn{1}{x,x} = \highinn{2}{x,x}=0$ and $\highinn{3}{x,x}>0$. 

The filtration $\filter^\bullet$ on $H$ induced by $\highinn{}{\cdot \,, \cdot}$ is given by
\[\filter^1 = H \supset \filter^2 =\{0\}\times \R \supset \filter^3=(0).\]
The restriction of $\highinn{1}{\cdot, \cdot}$ to the graded piece $\grm{}{1}H$ is the bilinear form $\highinn{1}{x,y} = x_1y_1$ which is definite. However, the restriction of $\highinn{1}{\cdot, \cdot}$ to $\grm{}{2}H = \{0\}\times \R$ is zero. It follows that $\highinn{}{\cdot,\cdot}$ is not an inner product.

We equip the torus $\T = \R^2 / \Z^2$ with the metric $\dist_t \colon \T \times \T \to [0, +\infty)$ induced by the inner product $\innone{t}{\cdot\,, \cdot}$ on $H$. Comparing with Theorem~\ref{thm:GHConvergenceAbstract} and taking into account the shape of the Gram matrix $M_t$, one might expect that $(\T, d_t)$ converges to the torus $\rquot{\R}{\Z}$  with its standard metric when $t$ tends to infinity.  However, it turns out that this drastically fails. In fact, the torus shrinks to a single point in the limit. Thus, although the products $\innone{t}{\cdot\,, \cdot}$, $t \in\R_+$, resemble the pullback family of an inner product, the flat tori show a very different behavior. 

\begin{prop}\label{prop:metriccollapse}
As $t$ tends to infinity, the flat tori $(\T, \dist_t)$, $t \in \R_+$, converge in the Gromov--Hausdorff sense to the metric space $X = \{ \bullet \}$ consisting only of one point.
\end{prop}
\begin{proof}

Fix a bounded subset $B \subset H$ such that $B$ covers the torus $\T$ via the projection map $H \to \T$. It suffices to prove that, as $t$ tends to infinity, the function $f_t(x) = \min_{\gamma \in \L} \norm{x -\gamma}_t$ goes to zero uniformly on $B$.  Clearly, this implies that $\lim_{t \to \infty} \dist_t(x,0) = 0$ uniformly for $x \in \T$ and thus $(\T, \dist_t)$ shrinks to a point in the limit $t \to \infty$.

Given $x = (x_1, x_2)$ in $B$, consider the lattice point $\gamma(x) = (0, \lceil \frac {x_1}{b_t} \rceil)$. Then,
\begin{align*}
f_t(x)^2  &\le \normsq{x -\gamma(x)}_t = x_1^2 + t^{-1} \gamma_2^2 - 2 b_t x_1 \gamma_2 + 2b_t x_1 x_2 + t^{-1}(x_2^2 - 2 x_2 \gamma_2) \\
&= x_1^2 + t^{-1} \gamma_2^2  - 2 b_t x_1 \gamma_2 + R_t(x),
\end{align*}
where the remainder term $R_t(x)$ tends to zero uniformly for $x \in B$ when $t$ tends to infinity.
Using the explicit estimate $(|x_1|/b_t) -1 \le |\gamma_2| \le (|x_1|/b_t) + 1$, we arrive at
\[
f_t(x) \le x_1^2 \Big(\frac{1}{t b_t^2} -1\Big) + \frac{2|x_1|}{t b_t}  + \frac{1}{t} + 2b_t |x_1| + R_t(x).
\]
Since $b_t \to 0$ and $t b_t^2 \to1$ for $t \to \infty$, it follows that $f_t(x) \to 0$ uniformly for $x \in B$.
\end{proof}

\begin{figure}[!h]
\centering
   \scalebox{.35}{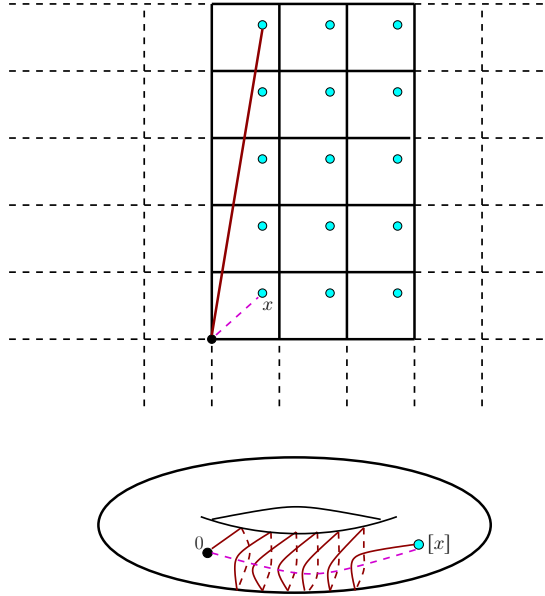}
\caption{An illustration of the collapse of the flat tori $(\T, \dist_t)$, $t \in \R_+$, in Proposition~\ref{prop:metriccollapse}. The blue points are all the predecessors in the universal cover $H$ of the point $[x]$ in the torus. The red line gives the shortest path joining the origin to the point $[x]$. As $t$ tends to infinity, the red path becomes more and more winding and its lengths tends to zero. Thus, the distance between the origin and $[x]$ goes to zero and the flat torus collapses to the origin. On the other hand, the length of the violet path remains bounded from below.}
\label{fig:metric-collapse}
\end{figure}

In the above example, the flat tori $(\T, \dist_t)$ undergo a collapse in the limit $t \to\infty$. Note that dimension drops appear naturally in the framework of Gromov--Hausdorff convergence of Riemannian manifolds, see e.g. Fukaya's survey paper~\cite{Fuk90}, and the works~\cite{Gro78, Ruh82, CG86, Fuk87, CG90, GTZ}. From this perspective, $\omega$-tamely degenerating scalar products provide a situation, where the dimension drop happens in a controlled way. Namely, by Theorem~\ref{thm:GHConvergenceAbstract}, after suitable renormalization, the dimension reduces from $\dim(H)$ to $\dim(\grm{}{1}{H})$. On the other hand, in the above example, one might naively expect a dimension drop of one, whereas actually the dimension decreases by two. 

It would be interesting to understand this phenomenon in a systematic way. We are therefore led to the following question,  which, to the best of our knowledge, does not seem to have been studied in the literature.

\begin{question}
Consider a real torus $\T = \rquot{H}{\L}$ of dimension $n$, defined by a finite dimensional real vector space $H$ and a lattice $\L \subset H$ of full rank. For $t \in \R_+$, let $\innone{t}{\cdot \,, \cdot}$ be a scalar product on $H$ and $\dist_t \colon \T \times \T \to [0, \infty)$ the associated distance function on the torus $\T$.

Consider the sequence of flat tori $(\T, \dist_t)$, $t \in \R_+$ and suppose that $(\T, \dist_t)$ converges in the Gromov--Hausdorff sense to a torus $X$ of dimension $m \le n$. How can the dimension drop $d= n-m$ be described in terms of the scalar products $\innone{t}{\cdot \,, \cdot}$?
\end{question}

\subsection{Voronoi cells of metric graphs and tropical curves} In this section, we discuss the structure of the Voronoi cells from Section~\ref{ss:degeneration-metric-graphs}. We fix a finite connected essential graph $G = (V,E)$.

\subsubsection{Voronoi cell of a metric graph} Let $l \colon E \to (0, +\infty)$ be an edge length function and denote by $\mgr$ the  metric graph associated to the pair $(G,l)$. Consider $H=H_1(\mgr, \R)$ endowed with the polarization $\innone{l}{\cdot\,,\cdot} \colon H\times H \to\R$.   Let $\L =H_1(\mgr, \Z)$ be the lattice of integer points in $H$. The Voronoi cell $W_l =  \Vor_{\L, \innone{l}{\cdot\,, \cdot} }(0)$ associated to $\L$ and $\innone{l}{\cdot\,, \cdot}$ has the following explicit description. 

A {\em circuit} in $\mgr$ is an element $\gamma \in \L$ which verifies $\gamma(e) =\pm 1$ for any edge $e \in \E$. Equivalently, a circuit of $\mgr$ is obtained by taking a cycle in $G$ and orienting it in one of the two possible ways to obtain an oriented cycle. Let $\S$ be the set of circuits of $\mgr$. Note also that $\S$ is a finite subset of $\L$.
\begin{thm}\label{thm:voronoi_graph} 
The Voronoi cell $W_l$ is the set of all points $x \in H$ which verify the inequalities
\begin{align*}
2 \innone{l}{x, \gamma} \leq \innone{l}{\gamma, \gamma} \qquad \textrm{ for any circuit } \gamma \in \S.  
\end{align*}
\end{thm}
The proof of this theorem is based on the fact that the circuits of the graph $G$ form a positive set of generators for the lattice $\L$, see  Section~\ref{ss:TropicalVoronoiCell} below. Using this theorem, it is possible to deduce the following description of the face structure of $W_l$.

Recall that a \emph{spanning subgraph} $R$ of $G$ is a subgraph which has the same vertex set $V$ and an edge set $F \subset E$. We say that $R$ is \emph{bridgeless} if there is no edge in $F$ whose removal from $R$ increases the number of connected components of $R$. A \emph{strongly connected orientation} of a bridgeless subgraph $R$ is an orientation such that each pair of vertices $u$ and $v$ in the same connected component of $R$ can be connected by oriented paths, one from $u$ to $v$ and another from $v$ to $u$.
 
\begin{thm} \label{thm:face-voronoi-graphs} There is a bijection between the faces of $W_l$ and strongly connected orientations of bridgeless spanning subgraphs of $G$. Under this bijection, a face $\tau$ is included in a face $\sigma$ in $W_l$ if and only if the subgraph $R_\sigma$ corresponding to $\sigma$ is included in the subgraph $R_\tau$ corresponding to $\tau$, and the orientation on $R_\sigma$ coincides with the one induced from the orientation on $R_\tau$. 
\end{thm}
\begin{proof} See~\cite{Ami10} for a proof. Alternatively, one can proceed by showing first that the normal fan of the Voronoi cell $W_l$ corresponds to the arrangement of hyperplanes in $H = H_1(\mgr, \R) \subset \R^E$ induced by the coordinate hyperplanes in $\R^E$, and then use the results of~\cite[Section 8]{GZ83} to deduce the bijection. 
\end{proof}

\subsubsection{Tropical Voronoi cell} \label{ss:TropicalVoronoiCell}

Consider a tropical curve $\curve$ associated to a weighted layered graph $(G, \pi, \ell)$, as in Section~\ref{ss:TropicalCurves}. Let $\highinn{}{\cdot\,,\cdot}\colon H_1(\curve, \R) \times  H_1(\curve, \R) \to \R^r$ be the corresponding polarization, and denote by $\Vor_{\curve}(0) \subset H_1(\curve, \R)$ the Voronoi cell with respect to the lattice $\L = H_1(\curve, \Z)$. We have the following description.

\begin{thm}\label{thm:voronoi_tropical_curve}
The Voronoi cell $\Vor_\curve(0)$ is the set of all points $x \in H_1(\curve, \R)$ which verify the inequalities 
\begin{align*}
2 \highinn{}{x, \gamma} \preceq \highinn{}{\gamma, \gamma} \qquad \textrm{ for any circuit } \gamma \in \S.  
\end{align*}
\end{thm}

The proof of this theorem is based on the following definition.

\begin{defi}[Positively generating set] \label{defi:positive-generating} Let $(H, \highinn{}{\cdot\,,\cdot})$ be an inner product space and $\L\subset H$ an admissible lattice. A subset $\S \subset \L$ is said to \emph{positively generate} $\L$ if it verifies the following property: Each element $\eta \in \L$ can be written as a positive linear combination $\sum_{\gamma\in A} a_\gamma \gamma$ for a subset $A \subset \S$ and positive integer coefficients $a_\gamma$ such that moreover, $\highinn{}{\gamma,\gamma'} \succeq 0$ for any pair $\gamma, \gamma' \in A$.
\end{defi}

\begin{lem}\label{lem:voronoi_positive_generating_set}  Let $(H, \highinn{}{\cdot\,,\cdot})$ be an inner product space and $\L\subset H$ an admissible lattice. Suppose that $\S$ positively generates $\L$. Then, the Voronoi cell of the origin $\Vor_\L(0)$ is given by 
\[
\Vor_\L(0) = \left\{x \in H \st 2\highinn{}{x,\gamma} \leq \highinn{}{\gamma ,\gamma } \text{ for any $\gamma \in \S$}\right\}.
\]
\end{lem}
\begin{proof} It suffices to show that the above inequalities imply all the inequalities
\[
2\highinn{}{x,\eta} \leq \highinn{}{\eta,\eta}, \qquad \eta\in \L.
\]
We write $\eta = \sum_{\gamma\in A} a_\gamma \gamma$ with $A \subset\S$, $a_\gamma >0$ positive integer, and $\highinn{}{\gamma,\gamma'} \succeq 0$ for any pair $\gamma, \gamma' \in A$. Then,
\begin{align*}
2\highinn{}{x,\eta} &= 2\sum_{\gamma \in A} a_\gamma \highinn{}{x,\gamma} \preceq \sum_{\gamma \in A} a_\gamma \highinn{}{\gamma,\gamma}\\
& \preceq \sum_{\gamma \in A} a^2_\gamma \highinn{}{\gamma, \gamma}\preceq \highinn{}{\sum_{\gamma \in A} a_\gamma \gamma,\sum_{\gamma \in A} a_\gamma \gamma} = \highinn{}{\gamma ,\gamma }. 
\end{align*}
In the last inequality, we used that $a_\gamma>0$ and $\highinn{}{\gamma,\gamma'} \succeq 0$ for any $\gamma, \gamma'\in A$.
\end{proof}

\begin{lem}\label{lem:circuits_positive_generating_set} Let $\curve$ be a tropical curve, $H = H_1(\curve, \R)$ and $\L = H_1(\curve, \Z)$. Let $\highinn{}{\cdot\,,\cdot}$ be the polarization in $\curve$. Then, the set of circuits $\S \subset \L$ positively generates $\L$.
\end{lem}

\begin{proof} This is a well-known property for graphs, see e.g.~\cite[Lemma 8.3]{GZ83}.  For $\eta \in \L\setminus \{0\}$, define $\supppm{\eta}$ as the set of oriented edges $e\in \E$ such that $\pm \eta(e) >0$. The oriented graph $(V, \suppp{\eta})$ has an oriented cycle $\gamma_1$. Let $\eta_1 =\eta- \gamma_1$. We have $\suppp{\eta_1} \subset \suppp{\eta}$. If $\eta_1 \neq0$, then we find  an oriented cycle $\gamma_2$ in the oriented graph $(V, \suppp{\eta_1})$. Proceeding this way, we find a decomposition $\eta = \gamma_1 + \dots + \gamma_k$ such that $\suppp{\gamma_j} \subseteq \suppp{\eta}$ and $\suppm{\gamma_j} \subseteq \suppm{\eta}$. The definition of polarization implies that $\highinn{}{\gamma_i, \gamma_j}\succeq 0$.  This implies that $\S$ is a positive generating set.
\end{proof}

\begin{proof}[Proof of Theorems~\ref{thm:voronoi_graph} and~\ref{thm:voronoi_tropical_curve}] Follows by combining Lemma~\ref{lem:voronoi_positive_generating_set} and Lemma~\ref{lem:circuits_positive_generating_set}. 
\end{proof}
We have the following higher rank analog of Theorem~\ref{thm:face-voronoi-graphs}.

\begin{thm}[Face structure of the closure of the tropical Voronoi cell] There is a bijection between the faces of $\overline \Vor_{\curve}(0)$ and tuples $\underline D=(D_1, D_2, \dots, D_r, D_\fin)$ consisting for each $j=1, \dots, r, \fin$ of a strongly connected orientation $D_j$ of a bridgeless spanning subgraph $R_{j}$ of the graded minor $\grm{}{j}(G)$. 

Let $R_\sigma$ be the spanning subgraph of $G$ obtained by taking the union of the edges appearing in $R_j$s oriented according to $D_j$s. Under the above bijection, a face $\tau$ is included in a face $\sigma$ in $\overline \Vor_{\curve}(0)$ if and only if the subgraph $R_\sigma$ corresponding to $\sigma$ is included in the subgraph $R_\tau$ associated to $\tau$, and the orientation on $R_\sigma$ coincides with the one induced from the orientation on $R_\tau$. 
\end{thm}

It would be interesting to obtain a description of the boundary points of the Voronoi cell $\Vor_{\curve}(0)$. In addition, it would be nice to understand which parts of the Voronoi cells $\Vor_{\mgr_t}(0)$ converge to a given face of $\overline \Vor_{\curve}(0)$, when a family of metric graphs $(\mgr_t)_{t \in \R_+}$ converges to $\curve$.

\subsubsection{Positive generation property for admissible lattices}  Let $\L\subset H$ be an admissible lattice in an inner product space $(H,\highinn{}{\cdot\,,\cdot})$. 
\begin{question} Does $\L$ have a finite, positively generating set $\S$ in the sense of Definition~\ref{defi:positive-generating}? 
\end{question}
If the rank of the inner product is one, it is possible to show that this is the case. To see this, without loss of generality, we can assume that the vector space generated by $\L$ is $H$. We then take a complete fan $\Sigma$ in $H$ which we suppose to be rational with respect to $\L$, and such that the scalar product of any two non-zero vectors lying in the same cone is positive. By Gordon's lemma (see \cite{AB86, Ful93}), for each cone $\sigma \in \Sigma$, the semigroup $\sigma \cap \L$ is finitely generated. We choose a finite set of generators $\S_\sigma$ for $\sigma \in \Sigma$. Then, the set $\S\coloneqq \bigcup_{\sigma\in\Sigma}\S_\sigma$ is finite and positively generates $\L$.

\subsection{Inner products with value in general ordered vector spaces}\label{ss:InnerProductsGeneralValueGroup}
So far in this paper we have studied inner products with values in the ordered vector space $(\R^r, \lexeq)$. In this section, we briefly explain how the set-up generalizes to inner products with values in an arbitrary totally ordered vector space $(\Lambda, \preceq)$.

Recall the following basic result (see~\cite[Section XV.2]{Bir-Lattice} and \cite{HW52}).
\begin{thm} \label{thm:Birkhoff}
Let $(\Lambda, \preceq)$ be a totally ordered real vector space of dimension $r \in \N$. Then $(\Lambda, \preceq)$ is isomorphic to $(\R^r, \lexeq)$.
\end{thm}

We will use the following lemma to compare different isomorphisms from $(\Lambda, \preceq)$ to $(\R^r, \lexeq)$. The proof is elementary and we omit it.
\begin{lem} \label{lem:LexicographicAutomorphisms}
Let $r \in \N$ and consider the ordered vector space $(\R^r, \lexeq)$. Then, the isomorphisms $\psi \colon (\R^r, \lexeq) \to(\R^r, \lexeq)$ are in bijection with matrices $B = (b_{ij})_{1 \le i,j \le r} \in \R^{r \times r}$ which are lower triangular and satisfy $b_{jj}>0$, $j \in [r]$. Here, $B \in \R^{r \times r}$ is identified with the linear map $\psi_B (x) = B x$, $x \in \R^r$. Moreover, let
\[
\Lambda^\bullet  \colon  \qquad \Lambda^1 \coloneq \Lambda\supset \Lambda^2 \supset \dots \supset \Lambda^{r} \supset \Lambda^{r+1} \coloneq (0)
\]
be the decreasing filtration on $(\R^r, \lexeq)$ defined in \eqref{eq:FiltrationRr}, that is,
\[
\Lambda^j \coloneqq \left\{a = (a_1, \dots, a_r)  \in \R^r \st a_1 = \dots =a_{j-1}=0\right\}, \qquad j \in [r+1].
\]
Then, the subspaces $\Lambda^j \subset H$, $j \in[r+1]$, are invariant under every isomorphism of $(\R^r, \lexeq)$. 

\end{lem}

\smallskip

Let $\Lambda$ be a real vector space of dimension $r \in \N$ endowed with an order $\preceq$. Fix an isomorphism $\varphi \colon (\Lambda, \preceq) \to (\R^r, \lexeq)$.  Let $H$ be a real vector space of finite dimension.

\begin{defi} \label{defi:GeneralVectorSpaceInnerProduct}
An inner product on $H$ with values in $(\Lambda, \preceq)$ is a symmetric bilinear form $\highinn{}{\cdot \,, \cdot} \colon H \times H \to \Lambda$ such that the bilinear form $[\cdot \,, \cdot] = \varphi(\highinn{}{\cdot \,, \cdot})$ is an inner product on $H$ with values in $\R^r$ (in the sense of Definition~\ref{defi:InnerProduct}).
\end{defi}
Lemma~\ref{lem:LexicographicAutomorphisms} ensures that Definition~\ref{defi:GeneralVectorSpaceInnerProduct} does not depend on the choice of the isomorphism $\varphi \colon (\Lambda, \preceq) \to (\R^r, \lexeq)$. Indeed, let $\{\cdot \,, \cdot\} = \widetilde \varphi (\highinn{}{\cdot \,, \cdot})$ be obtained from another isomorphism $\widetilde \varphi \colon (\Lambda, \preceq) \to (\R^r, \lexeq)$. The composition $\psi = \widetilde \varphi \circ \varphi^{-1}$ is an isomorphism of $(\R^r, \lexeq)$. By Lemma~\ref{lem:LexicographicAutomorphisms}, $\psi$ is given by a lower triangular matrix $B \in \R^{r \times r}$ with strictly positive diagonal entries $b_{jj} >0$, $j \in [r]$. 
Since the subspaces $\Lambda^j \subset \R^r$, $j \in [r]$, are invariant under $\psi$, the bilinear forms $[\cdot \, , \cdot]$ and $\{\cdot \, , \cdot \}$ induce the same filtration $\filter^\bullet = (\filter^j)_{j=1}^r$ on $H$. For all $j \in [r]$, we moreover have
\[
\{x,y\}_j = b_{jj} [x,y]_j, \qquad x,y \in \filter^j.
\]
Taking into account Definition~\ref{defi:InnerProduct}, it follows that $[\cdot \,, \cdot] \colon H \times H \to \R^r$ is an inner product with values in $\R^r$ if and only if $\{ \cdot\, , \cdot \} \colon H \times H \to \R^r$ is so.

\smallskip

Analogously, using an isomorphism $\varphi \colon (\Lambda, \preceq) \to (\R^r, \lexeq)$, one can introduce the filtration $\filter^\bullet$, the graded pieces $\grm{}{j}H = \rquot{\filter^j}{\filter^{j+1}}$, $j \in [r]$, the lifting operators $\proj_j^\ast \colon \grm{}{j}H \to H$, $j \in [r]$, and the almost orthogonal decomposition $H = \aplus_{j=1}^r H_j$ for an inner product $\highinn{}{\cdot \, , \cdot}$ with values in $\Lambda$. These objects do not depend on the choice of isomorphism $\varphi \colon (\Lambda, \preceq) \to (\R^r, \lexeq)$.  Note however that the scalar products $\highinn{j}{\cdot \,, \cdot} = (\varphi(\highinn{}{\cdot \,, \cdot}))_j$ on the graded minors $\grm{}{j}H$, $j \in [r]$, depend on $\varphi$ and change by a multiplicative constant $\lambda_j>0$ when changing $\varphi$.

\smallskip

Concerning our results on higher rank Voronoi decompositions, note that the definition of the Voronoi cell \eqref{eq:DefinitionVoronoiCell} directly generalizes to an inner product $\highinn{}{\cdot \, , \cdot} \colon H \times H \to \Lambda$. Fixing an isomorphism $\varphi \colon \Lambda \to \R^r$, one immediately obtains the results in Section~\ref{sec:HigherRankVoronoiDecompositions}, Section~\ref{sec:StructureVoronoiCells} and Section~\ref{sec:AdmissibleLattices} for $\highinn{}{\cdot \, , \cdot}$.

\smallskip

Finally, we discuss briefly how the notions of pullback families and tamely degenerating families depend on the choice of isomorphism. Let $\highinn{}{\cdot \, , \cdot} \colon H \times H \to \Lambda$ be an inner product with values in $\Lambda$. As above, consider two inner products $[\cdot \,, \cdot] = \varphi (\highinn{}{\cdot \,, \cdot})$ and $\{\cdot \,, \cdot\} = \tilde \varphi(\highinn{}{\cdot \,, \cdot})$ with values in $\R^r$ obtained from isomorphisms $\varphi \colon (\Lambda, \preceq) \to (\R^r, \lexeq)$ and $\tilde \varphi \colon (\Lambda, \preceq) \to (\R^r, \lexeq)$, respectively. Let $B \in \R^{r \times r}$ be the matrix describing the isomorphism $\psi = \tilde \varphi \circ \varphi^{-1}$ of $(\R^r, \lexeq)$.

For a vector $\underline L = (L_1, \dots, L_r) \in \R_+^r $, the pullback bilinear forms for $[\cdot \,, \cdot]$ and $\{\cdot \,, \cdot\}$ (see Definition~\ref{defi:PullbackBilinearForm}) are related by
\[ \{\cdot \,, \cdot\}_{\underline L} = [\cdot \, , \cdot]_{B^\transpose \cdot \underline L}, \qquad
[\cdot \,, \cdot]_{\underline L} = [\cdot \, , \cdot]_{(B^{-1})^\transpose \cdot \underline L}.
\]
Thus, $[\cdot \,, \cdot]$ and $\{\cdot \, , \cdot \}$ have the same pullback families (see Section~\ref{ss:PullbackFamilies} for the definition). (Note in particular that the parameters $\underline L_t \in \R_+^r$, $t \in \R_+$ of a pullback family, satisfy $B^\transpose \cdot \underline L_t \in \R_+^r$ and $(B^{-1})^\transpose \cdot \underline L_t \in\R_+^r$ for all $t$ large.)

A similar relation holds for tamely degenerating families of scalar products. More precisely, let $\innone{}{\cdot \,, \cdot}$, $t \in \R_+$, be a family of scalar products on $H$. Then $\innone{}{\cdot \, ,\cdot}$, $t \in \R_+$, degenerates tamely ($\omega$-tamely) to $[\cdot \,, \cdot]$ with parameters $\underline L_t \in \R_+^r$, if and only if $\innone{}{\cdot \,, \cdot}$, $t \in \R_+$, degenerates tamely ($\omega$-tamely) to $\{\cdot \,, \cdot\}$ with parameters $\widetilde{\underline L}_t \in \R_+^ r$ given by $\tilde L_{t,j} = \frac{1}{b_{jj}} L_{t,j}$, $j \in [r]$.

Using these relations, one can formulate analogs of the results in Section~\ref{sec:Hausdorff_convergence} and Section~\ref{sec:metric_degeneration_general_tori} for inner products with values in general ordered vector spaces $(\Lambda, \preceq)$.

\subsection{Tame equivalence of higher rank inner products} \label{ss:TameEquivalence}
The results in the preceding sections suggest to view tame degenerations as a kind of convergence of scalar products to inner products of higher rank. In this section, we discuss the uniqueness of limits in this context. 

\smallskip

 Let $\highinn{}{\cdot\,,\cdot} \colon H \times H \to \R^r$ and $\highinnprime{}{\cdot\,,\cdot}\colon H \times H \to \R^r$ be inner products on a finite dimensional real vector space $H$. We call $\highinn{}{\cdot\,,\cdot}$ and $\highinnprime{}{\cdot\,,\cdot}$ \emph{tame equivalent} (resp. \emph{$\omega$-tame equivalent}), if they have the same tamely degenerating families (resp. $\omega$-tamely degenerating families) with the same parameters. That is, a family $\innone{t}{\cdot \, , \cdot}$, $t \in \R_+$, of scalar products degenerates tamely (resp. $\omega$-tamely) to $\highinn{}{\cdot\,,\cdot}$ with parameters $\underline L_t$, $t \in \R_+$, if and only if it degenerates tamely (resp. $\omega$-tamely) to $\highinnprime{}{\cdot\,,\cdot}$ with parameters $\underline L_t$, $t \in \R_+$. Clearly, ($\omega$-)tame equivalence defines an equivalence relation on the space of inner products on $H$ with values in $\R^r$.

\smallskip

The next theorem characterizes tame equivalent and $\omega$-tame equivalent inner products.

\begin{thm} \label{thm:TameEquivalence} Let $\highinn{}{\cdot\,,\cdot} \colon H \times H \to \R^r$ and $\highinnprime{}{\cdot\,,\cdot}\colon H \times H \to \R^r$ be inner products. Denote by $\filter^\bullet$ and ${\filter'}^\bullet$ the associated filtrations, and by $H_j$ and $H'_j$, $j \in [r]$, the spaces in the associated almost orthogonal decompositions $H=H_1 \aplus \dots \aplus H_r$ and $H=H'_1\aplus \dots \aplus H_r'$. Then:

\begin{itemize}
\item [(i)] $\highinn{}{\cdot\,,\cdot}$ and $\highinnprime{}{\cdot\,,\cdot}$ are $\omega$-tame equivalent if and only if  $\filter^j = {\filter'}^j$ for all $j \in [r]$, and
\[
\highinn{j}{x,y} = \highinnprime{j}{x,y } \qquad \forall x,y \in \filter^j.
\]
\item [(ii)] $\highinn{}{\cdot\,,\cdot}$ and $\highinnprime{}{\cdot\,,\cdot}$ are tame equivalent if and only if they are $\omega$-tame equivalent, and additionally, we have $H_j = H'_j$ for all $j \in [r]$.
\end{itemize}
\end{thm}

\begin{proof} In order to prove the claim in $(i)$, assume first that the two filtrations $\filter^\bullet$ and ${\filter'}^\bullet$ coincide and the equality $\highinn{j}{x,y} = \highinnprime{j}{x,y }$ holds for all $x,y \in \filter^j = {\filter'}^j$. Let $\innone{t}{\cdot\,, \cdot}$, $t \in \R_+$, be a family of scalar products which $\omega$-tamely degenerates to $\highinn{}{\cdot \, , \cdot}$ with parameters $\underline L_t$, $t \in \R_+$. We verify the properties \ref{item:tame-a}, \ref{item:weaktame-b}, and \ref{item:tame-c} in Definition~\ref{def:TameDegenerations} for $\highinnprime{}{\cdot\,,\cdot}$ and the parameters $\underline L_t$, $t \in \R_+$. Clearly, \ref{item:tame-a} and \ref{item:weaktame-b} are automatically satisfied.

Fix $j \in [r]$. In the following, we will prove that for given $\varepsilon >0$, we have
\begin{equation} \label{eq:TameProof}
\abs{L_{t,j}^{-1} \innone{t}{x, y} -  \highinn{j}{x, y}} \leq \varepsilon \refnorm{x}\cdot\refnorm{y}, \qquad \forall x,y \in \filter^j \text{ and $t$ large}.
\end{equation}
This then implies property \ref{item:tame-c} for $\highinnprime{}{\cdot\,,\cdot}$, since $H'_j \subset {\filter'}^j = \filter^j$ and $\highinn{j}{\cdot \, ,\cdot} = \highinnprime{j}{\cdot \,, \cdot}$ on $\filter^j = {\filter'}^j$ by assumption.

Writing $x \in \filter^j$ as $x = \sum_{k \ge j} x_k$ with $x_k \in H_k$, we can apply property \ref{item:weaktame-b} for $\highinn{}{\cdot\,,\cdot}$ to get
\begin{align*}
L_{t,j}^{-1}  | \innone{t}{x,y} - \innone{t}{x_j,y_j}| &= L_{t,j}^{-1}  | \innone{t}{x_j,y-y_j} + \innone{t}{x-x_j,y}| \\
& \le C \frac{L_{t,j+1}}{L_{t,j}} \Big ( \refnorm{x_j} (\sum_{k=j+1}^r \refnorm{y_k}) + (\sum_{k=j+1}^r \refnorm{x_k}) \refnorm{y} \Big )
\end{align*}
for $x,y \in \filter^j$ and $t$ large. Recalling that $\highinn{j}{x,y} = \highinn{j}{x_j, y_j}$ for $x,y \in \filter^j$, and applying the properties \ref{item:tame-a}, \ref{item:weaktame-b} and \ref{item:tame-c} for $\highinn{}{\cdot \,, \cdot}$, we obtain the following. For every $\varepsilon >0$, the estimate
\[
  | L_{t,j}^{-1}\innone{t}{x,y} - \highinn{j}{x,y}| \le | L_{t,j}^{-1}\innone{t}{x_j,y_j} - \highinn{j}{x_j,y_j}| + L_{t,j}^{-1}  | \innone{t}{x,y} - \innone{t}{x_j,y_j}| \le \varepsilon (\sum_{k =j}^r \refnorm{x_k}  ) (\sum_{k =j}^r \refnorm{y_k}  ) 
\]
holds for all $x , y \in \filter^j$ and $t$ large. On the other hand, since all norms on $H$ are equivalent, there exists a constant $D >0$ with $\sum_{k=1}^r \refnorm{z_k} \le D \refnorm{z}$ for all $z \in H$. Thus, \eqref{eq:TameProof} holds and we have proven the implication ``$\Leftarrow$" in $(i)$.

To prove the converse implication ``$\Rightarrow$" in (i), fix a common $\omega$-tamely degenerating family $\innone{t}{\cdot \, , \cdot}$, $t \in \R_+$, for $\highinn{}{\cdot \, , \cdot }$ and $\highinnprime{}{\cdot \, , \cdot }$ with parameters $\underline L_t \in \R_+^r$. As follows from \eqref{eq:TameProof} and property \ref{item:tame-a}, the filtration $\filter^\bullet = (\filter^j)_{j=1}^r$ can be expressed as
\[
\filter^j = \{ x \in H \st \lim_{t \to \infty} L_{t,i}^{-1} \innone{t}{x,x} = 0 \text{ for all } i = 1, \dots, j-1 \}.
\]
Since this expression only involves the scalar products $\innone{t}{\cdot \, , \cdot}$ and parameters $\underline L_t$, $t \in \R_+$, we infer that $\filter^j  = {\filter'}^j$ for $j \in [r]$.  The equality $\highinn{j}{\cdot\,,\cdot} = \highinnprime{j}{\cdot\,,\cdot}$ on $\filter^j = {\filter'}^j$ follows from \eqref{eq:TameProof}.

The implication ``$\Leftarrow$" in (ii) is obvious, since the properties \ref{item:tame-a}, \ref{item:tame-b}, and \ref{item:tame-c} of tame degenerations are formulated only in terms of the spaces $H_j$ and components $\highinn{j}{\cdot\,,\cdot}$.

To prove the implication ``$\Rightarrow$" in (ii),  fix a common tamely degenerating family $\innone{t}{\cdot \, , \cdot}$, $t \in \R_+$, for $\highinn{}{\cdot \, , \cdot }$ and $\highinnprime{}{\cdot \, , \cdot }$ with parameters $\underline L_t \in \R_+^r$. By the above argument, $\filter^j = {\filter'}^j$ and $\highinn{j}{\cdot \, , \cdot } = \highinnprime{j}{\cdot \, , \cdot }$ on $\filter^j = {\filter'}^j$ for all $j \in [r]$. In order to complete the proof, we show that
\begin{equation} \label{eq:HjExpressionDegeneration}
H_j = \{ x \in \filter^j \st   \lim_{t \to \infty} L_{t,i}^{-1} \innone{t}{x,y} = 0 \text{ for all $y \in \filter^i$ and $i \in \{j+1, \dots, r\}$} \}.
\end{equation}
Since this expression only depends on the scalar products $\innone{t}{\cdot \, , \cdot}$ and parameters $\underline L_t$, $t \in \R_+$, it follows that $H_j  = H_j'$ for all $j \in [r]$. If $x \in H_j$, then the properties \ref{item:tame-a} and \ref{item:tame-b} of $\highinn{}{\cdot \, , \cdot }$ imply that $x$ belongs to the right hand side of \eqref{eq:HjExpressionDegeneration}. Conversely, suppose that $x \in H$ belongs to the right hand side of \eqref{eq:HjExpressionDegeneration}. Writing $x = \sum_{k\ge j} x_k$ for $x_k \in H_k$, we have to prove that $x_k = 0$ for all $k >j$. By assumption, we have
\[
\lim_{t \to \infty} L_{t,k}^{-1} \innone{t}{x,x_k} = 0.
\]
On the other hand, properties \ref{item:tame-a}, \ref{item:tame-b} and \ref{item:tame-c} of $\highinn{}{\cdot \, , \cdot }$ imply that
\[
\lim_{t \to \infty} L_{t,k}^{-1} \innone{t}{x,x_k} =\lim_{t \to \infty} \sum_{i =j}^{k-1} L_{t,k}^{-1} \innone{t}{x_i,x_k} + L_{t,k}^{-1} \innone{t}{x_k,x_k} + \sum_{i =k+1}^{r} L_{t,k}^{-1} \innone{t}{x_i,x_k} = \highinn{k}{x_k, x_k}.
\]
Since $\highinn{k}{\cdot \,, \cdot} = 0$ is a scalar product on $H_k$, we obtain $x_k = 0$ for $k >j$, as desired.
\end{proof}

Theorem~\ref{thm:TameEquivalence} suggests a way to define a partial compactification of the space of scalar products on a fixed vector space $H$. We discuss this in the next section.

Note also that tame equivalent inner products $\highinn{}{\cdot\,,\cdot}$ and $\highinnprime{}{\cdot\,,\cdot}$ on $H$ have the same admissible discrete sets $S \subset H$. Moreover, by Theorem~\ref{thm:Hausdorff_voronoi}, the associated Voronoi cells have the same closure, that is, $\overline \Vor_{S}(\gamma) = \overline{\Vor_{S}'}(\gamma)$ for all $\gamma$ in an admissible discrete subset $S \subset H$. On the other hand, the non-closed Voronoi cells can differ, $ \Vor_{S}(\gamma) \neq \Vor_{S}'(\gamma)$. It would be certainly interesting to obtain a more accurate description of the Voronoi cells $\Vor_{S}(\gamma)$ for all inner products $\highinn{}{\cdot\,,\cdot}$ in a fixed tame equivalence class.

\subsection{Tame compactification of the cone of scalar products} \label{sec:tame-compactifications}
In this last section, we briefly discuss how higher rank inner products provide partial compactifications of the cone of scalar products. A more through discussion of this topic is beyond the scope of this paper and is postponed to a future work.

Let $H$ be a real vector space of finite dimension $n$, and consider the cone of scalar products on $H$. Fixing a basis for $H$, we can identify this cone with the cone $\mathscr S_n^+$ of symmetric positive definite matrices of order $n$,
\[\mathscr S^+_n \coloneqq \bigl\{A \in M_n(\R) \st A^\transpose=A \quad \textrm{and} \quad A \succ 0\bigr\}.\]
The cone $\mathscr S^+_n$ is an open subset of the space of symmetric matrices $\mathscr S_n$, which is a real vector space of dimension $n(n+1)/2$.

Using higher rank inner products, and tame and $\omega$-tame degenerations, we can define two partial compactifications of $\mathscr S^+_n$. We add tame, resp. $\omega$-tame, equivalence classes of higher rank inner products to the cone $\mathscr S^+_n$ and endow naturally the new space with a relevant topology such that tame, resp. $\omega$-tame, degenerating families lead to convergent sequences. Note that our Theorem~\ref{thm:TameEquivalence} provides a description of tame and $\omega$-tame equivalence classes. We call these partial compactifications  the \emph{tame} and \emph{$\omega$-tame (partial) compactifications} of $\mathscr S^+_n$, respectively. Both of these spaces are Hausdorff. 

These compactifications are in the esprit of the family of compactifications produced in our previous work~\cite{AN2} for polyhedral cones (although, the cone $\mathscr S^+_n$ is not polyhedral).  They refine the data of the filtration on the underlying vector space, induced by the multi-scale behavior of the vector lengths in the limit, by additional coordinates which provide the data of the higher rank inner product.

Classically, spaces such as the cone $\mathscr S^+_n$, or its Hermitian analogue, can be compactified in multiple ways, as instances of compactifications produced by Tits, Baily-Borel, Borel-Serre, Satake, Furstenberg, Martin, and Karpelevic. (For example, Tits' compactification of the cone $\mathscr S^+_n$ is obtained by endowing it with a natural Riemannian metric that turns it into a Hadamard space, and adding the ideal boundary.) A uniform presentation of these compactifications is provided in the work by Borel and Ji~\cite{BJ06, BJ07}. 

Tame compactifications provide a new way of (partially) compactifying this cone, as a limit of toroidal compactifications enriched by ordered conical blow-ups, as we do in producing higher rank compactifications of polyhedral cones in~\cite{AN2}. Going beyond the case of positive definite matrices, we expect that a generalization of the ideas developed here would produce multi-scale tame compactifications of more general Hadamard, symmetric and locally symmetric spaces. This will be studied in our future work.

\bibliographystyle{alpha}
\bibliography{bibliography}

\end{document}